\documentclass[3p,final,13pt]{elsarticle}

\pdfoutput=1 

\usepackage{lineno}
\modulolinenumbers[1]

\usepackage[pdfstartview=FitH,colorlinks=true,citecolor=blue,bookmarks=false]{hyperref}
\hypersetup{hypertexnames=false}

 \makeatletter
\def\ps@pprintTitle{%
 \let\@oddhead\@empty
 \let\@evenhead\@empty
 \def\@oddfoot{\centerline{\thepage}}%
 \let\@evenfoot\@oddfoot}
\makeatother  

\makeatletter
\providecommand{\doi}[1]{%
  \begingroup
    \let\bibinfo\@secondoftwo
    \urlstyle{rm}%
    \href{http://dx.doi.org/#1}{%
      doi:\discretionary{}{}{}%
      \nolinkurl{#1}%
    }%
  \endgroup
}
\makeatother
\usepackage{xcolor}
\definecolor{green1}{rgb}{0,0.314,0} 
\definecolor{green2}{rgb}{0,0.498,0} 
\definecolor{brown1}{rgb}{0.545,0.270,0.0745} 
\usepackage{bbm}
\usepackage{mathtools}
\usepackage{amsmath}
\usepackage{amssymb}
\usepackage{amsthm}
 \usepackage{mathrsfs}
\usepackage{graphicx}%
\usepackage{graphics}%
\usepackage{subcaption}
\usepackage{overpic}%

\usepackage{wasysym}
\usepackage{relsize}

\newtheorem{proposition}{Proposition}[section]
\newtheorem{remark}{Remark}[section]
\newtheorem{definition}{Definition}[section]

\usepackage{appendix}

\graphicspath{ {./Figures/} }

\usepackage{accents}
\newcommand\ubar[1]{\underaccent{\bar}{#1}}

\renewcommand{\leq}{\leqslant}
\renewcommand{\geq}{\geqslant}

\begin{document}
\begin{frontmatter}
\title{%
Response of an oscillatory delay differential equation to a periodic stimulus\\
}%

\author[mcgill]{Daniel C. De Souza}
\ead{daniel.desouza@mail.mcgill.ca}

\author[mcgill2]{Michael C. Mackey}
\ead{michael.mackey@mcgill.ca}

\address[mcgill]{Department of Mathematics and Statistics, McGill University,  Montreal, Quebec, H3A 0B9, Canada}
\address[mcgill2]{Departments of Physiology, Physics and Mathematics, McGill University, 3655 Promenade Sir William Osler, Montreal, QC H3G 1Y6, Canada}

\date{\today}

\begin{abstract}
Periodic hematological diseases such as cyclical neutropenia or cyclical thrombocytopenia, with their characteristic oscillations of circulating neutrophils or platelets, may pose grave problems for patients.  Likewise, periodically administered chemotherapy has the unintended side effect of establishing periodic fluctuations in circulating white cells, red cell precursors and/or platelets. These fluctuations, either spontaneous or induced,  often have serious consequences for the patient (e.g. neutropenia, anemia, or thrombocytopenia respectively) which exogenously administered cytokines can partially correct. The question of when and how to administer these drugs is a difficult one for clinicians and not easily answered.  In this paper we use a simple model consisting of a delay differential equation with a piecewise linear nonlinearity, that has a periodic solution, to model the effect of a periodic disease or periodic chemotherapy.  We then examine the response of this toy model to both single and periodic perturbations, meant to mimic the drug administration, as a function of the drug dose and the duration and frequency of its administration to best determine how to avoid side effects.
\end{abstract}


\end{frontmatter}


\section{Introduction}\label{sec.introduction}

\par Hematopoiesis is the term for the  process of blood cell formation. Normally maintained at a homeostatic level within certain bounds (that vary between cell types), the numbers of circulating white cells, red cells, and platelets typically do not display any evidence of oscillatory dynamics.  However, there are many hematopoietic diseases (so called periodic diseases~\citep{Glass_1988}) in which cycling of one or more circulating blood cell types is seen.  Examples include cyclical thrombocytopenia (CT)~\citep{Langlois_2017}, cyclical neutropenia (CN)~\citep{Colijn_2005b} and periodic chronic myelogenous leukemia (PCML)~\citep{Colijn_2005a}, and there have been numerous mathematical modeling studies of these disorders aimed at their understanding and treatment~\citep{Craig_2016,Craig_2015,Foley_2009a}.

While cyclicity in circulating blood cell numbers is relatively rare in disease states, induced cyclicity of one or more circulating hematopoietic cell types as a byproduct of periodically administered chemotherapy is all too common.  This cycling is most often encountered in the neutrophils (with a concomitant risk of infection when neutrophil numbers fall to sufficiently low levels), but  also may be observed in the platelets (with an accompanying increased risk of bleeding and stroke at the low point of the platelet cycle, or thrombosis at the high point) as well as rarely in the erythrocytes (red blood cells, with accompanying anemia at the low point of the cycle).  The commonality of this cycling with its attendant side effects (infection, bleeding, anemia) is one of the primary reasons leading to an interruption and/or cessation of chemotherapy.

 In mammals hematopoiesis starts in the bone marrow with the proliferation and subsequent differentiation of hematopoietic stem cells (HSCs) into one of the three major cell lines, and ends  with the release of  mature blood cells into the circulation.
 Although all mature blood cells have the HSCs as their common origin, the control of their production is only partially understood~\citep{Beuter_2003}.  However, the broad outline is clear.

 The numbers of circulating blood cells are controlled by a delayed negative feedback mediated by cytokines, such as granulocyte colony-stimulating factor (G-CSF) for the white blood cells, thrombopoietin (TPO) for the platelets, and erythropoietin (EPO) for the red blood cells~\citep{Mackey_2017}. The periodic administration of chemotherapy, or the existence of  hematological disorders like CN, PCML, or CT,  may lead the level of peripheral blood cells to exhibit oscillations that are more or less regular~\citep{Craig_2015, Foley_2009a, Mackey_2017}. There is a vast literature of mathematical models that propose how to control chemotherapy side effects or understand periodic hematological disorders, see for example~\citet[Chapter 8]{Beuter_2003}, \citet{Foley_2009a} and~\citet{Menjouet_2016}.

A scalar delay differential equation (DDE), with a linear piecewise constant negative feedback nonlinearity, which captures the essence of the negative delayed feedback mechanism involved in the control of blood cells by cytokines, was analyzed in~\citet{Mackey_2017}.  This equation has an oscillatory solution (to mimic an inherent oscillation due to a periodic disease or induced, for example, by the administration of chemotherapy), and was used as a toy model to examine analytically the effect of a {\it single} perturbation on the limit cycle.  The single perturbation was applied at different points in the limit cycle to mimic the delivery of a cytokine (G-CSF, TPO, or EPO) and to examine the subsequent effect on the model oscillation.  The single perturbation amplitude and duration were related to the dose and time of administration of an exogenous cytokine,  as it is known that the timing of the administration of a cytokine can be crucial in its clinical effect~\citep{Mackey_2017}.  In~\citet{Mackey_2017} the results were limited to an examination of a single stimulus as the authors were unable to deal with the clinically more interesting case of a periodic stimulus.
 Here, we study the effect of a {\it periodic perturbation} on the dynamics of this delay differential equation.

This paper is organized as follows. Section~\ref{sec.model.disc} presents the model background and its DDE with discontinuous (Heaviside step function) delayed feedback and summarizes  some fundamental results from~\citet{Mackey_2017} concerning the response of the limit cycle to a single stimulus.  Section~\ref{sec.single} extends the analysis from~\citet{Mackey_2017} of the response of the periodic solution to a single stimulus.
We further define the concept of resetting time for a single perturbation and investigate its maximum and minimum values as function of the perturbation phase. We also describe a special case where changes in the phase and amplitude of the perturbation leads to solutions close to unstable limit cycles.

Section~\ref{sec.per} analyzes the response of the DDE 
to a periodic perturbation.
We define sufficient conditions to obtain a periodic solution such that all local minima are positive (clinically important), and also show examples of solutions with different numbers of local minima and maxima.  The proofs of all the results stated in the remarks and propositions are presented in Appendix~\ref{sec.app}.
Section~\ref{sec:treatment} considers our modeling results in the context of cytokine administration, and carries out a detailed comparison with a more comprehensive model of \cite{zhuge2012neutrophil}.  The penultimate Section~\ref{sec.num.exp.per.pert} presents a variety of bifurcation results in  our model system, while Section~\ref{sec.conc} gives a brief summary and prospects for further work.  

\section{Model Background}\label{sec.model.disc}

Here we consider the simple mathematical model used in~\citet{Mackey_2017} to describe the dynamics of a circulating blood cell population
$x(t)$ in which the cell death rate of circulating cells is denoted by $\gamma$ and their production rate is described by a delayed negative feedback mechanism $f(x(t-\tau))$. The delay $\tau$ captures the physiologically known delay due to cellular division, differentiation and maturation. The dynamics of $x(t)$ is taken to be described by (see~\citet{Mackey_2017})
\begin{equation}\label{dxdt0}
x^{\prime}(t) = -\gamma x(t) + f(x(t-\tau)),
\end{equation}
where the piecewise constant nonlinearity $f$ is of the form
\begin{equation}\label{ftau0}
f(x(t-\tau)) \coloneqq \left\{\begin{array}{ll}
b_{L}\qquad\mbox{for}\qquad x(t-\tau)<\theta,\\
b_{U}\qquad\mbox{for}\qquad x(t-\tau)\geq \theta,
\end{array}\right.
\end{equation}
with $\gamma>0$, $\theta>0$, $b_{L}>b_{U}>0$, $b_{U}\neq \gamma\theta$. To solve the initial value problem~\eqref{dxdt0} we must specify the initial function $\varphi:[-\tau,0]\to \mathbb{R}$
\begin{equation}\label{varphi}
x(t)=\varphi(t)\quad\mbox{for}\quad -\tau\leq t\leq 0.
\end{equation}

Denote the solution of~\eqref{dxdt0} with initial function~\eqref{varphi} by $x(t,\varphi)$ and the zeros of the solution $x(t,\varphi)$ for $t\in[-\tau,\infty)$ as the set of all $z_{j}$, with $j\in\mathbb{N}$, such that $x(z_{j},\varphi)=0$. Here we only consider initial (history) functions~\eqref{varphi} that are continuous $C([-\tau,0],\mathbb{R})$ with a finite number of zeros. Denote by $Z_{0}\subset C([-\tau,0],\mathbb{R})$ a set of history functions which has at most one zero on $[-\tau,0]$ and changes sign at this zero. Given a history function $\varphi_{0}\in Z_{0}$, it follows from~\citet[Section 3]{Mackey_2017} that $(x(t,\varphi_{0})-\theta)$ has a strictly increasing sequence of zeros $z_{j}=z_{j}(\varphi_{0})$ in $(0,\infty)$, $j\in\mathbb{N}$, such that
\begin{equation}\nonumber
z_{j}+\tau<z_{j+1}\quad\mbox{for all}\quad j\in\mathbb{N}.
\end{equation}

\begin{figure}[!htbp]
\centering
\includegraphics[width=1.0\textwidth]{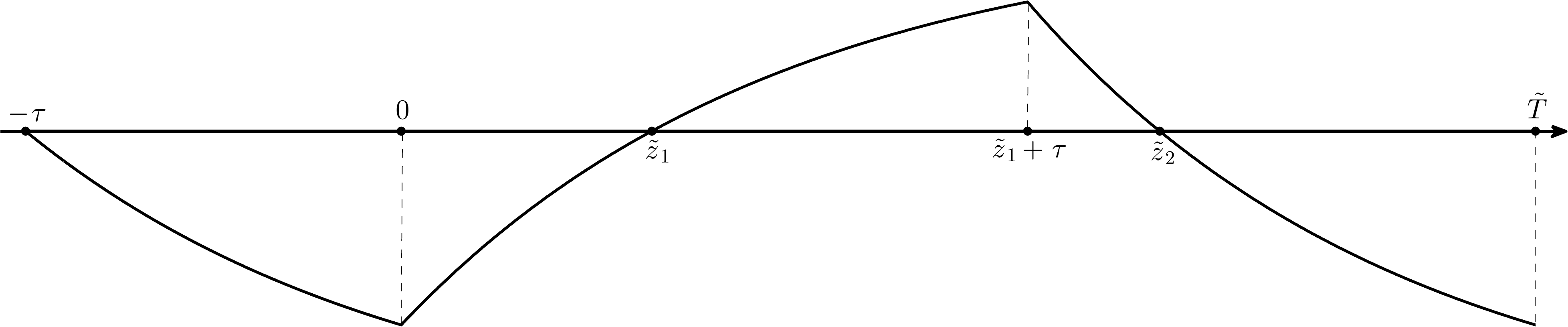}
\caption{A schematic representation of the limit cycle~\eqref{tildex} on the interval $[-\tau,\tilde{T}]$.}
\label{fig_limit_cycle}
\end{figure}

The model given by~\eqref{dxdt0} and~\eqref{ftau0} contains five parameters $\{\gamma,\tau,b_{U},b_{L},\theta\}$. Using the change of variables
\begin{equation}\label{transform}
\left\{\begin{array}{ll}
x(t) \to x(t) + \theta,\\
t \to \gamma t,\\
\tau \to \gamma \tau,\\
b_{L} \coloneqq \gamma(\theta+\beta_{L}),\\
b_{U} \coloneqq \gamma(\theta-\beta_{U}),
\end{array}\right.
\end{equation}
we can rewrite~\eqref{dxdt0}-\eqref{ftau0} as a function of three parameters $\{\tau,\beta_{U},\beta_{L}\}$, namely
\begin{equation}\label{dxdt}
x^{\prime}(t) = - x(t) + f(x(t-\tau)),
\end{equation}
where
\begin{equation}\label{ftau}
f(x(t-\tau)) \coloneqq \left\{\begin{array}{ll}
\phantom{-}\beta_{L}\qquad\mbox{for}\qquad x(t-\tau)<0,\\
-\beta_{U}\qquad\mbox{for}\qquad x(t-\tau)\geq 0,
\end{array}\right.
\end{equation}
with $-\beta_{U}<0<\beta_{L}$, $\tau>0$. Eq.~\eqref{dxdt} with $f(x(t-\tau))$ defined by~\eqref{ftau} captures the negative delayed feedback mechanism involved in the control of hematopoietic cells by cytokines.
For the initial function
\begin{equation}\nonumber
\tilde{x}(t)\coloneqq-\beta_{U}+\beta_{U}\mathrm{e}^{-(t+\tau)} \quad\mbox{for}\quad -\tau\leq t\leq 0,
\end{equation}
the solution of Eq.~\eqref{dxdt} is a limit cycle $\tilde{x}(t)$, and is given by~\eqref{tildex}~\citep{Mackey_2017} and is presented in Figure~\ref{fig_limit_cycle}.
\begin{equation}\label{tildex}
\tilde{x}(t) \coloneqq \left\{\begin{array}{ll}
\phantom{-}\beta_{L} + (\ubar{x}-\beta_{L})\mathrm{e}^{-t} \quad&\mbox{for}\quad 0\leq t\leq\tilde{z}_{1}+\tau,\\
 -\beta_{U} + (\bar{x}+\beta_{U})\mathrm{e}^{-(t-\tilde{z}_{1}-\tau)} \quad&\mbox{for}\quad \tilde{z}_{1}+\tau\leq t\leq\tilde{T}.
\end{array}\right.
\end{equation}
where the minimum $(\ubar{x})$ and maximum $(\bar{x})$ of $\tilde{x}(t)$ are given by \citep{Mackey_2017}
\begin{equation}\label{max.min}
\ubar{x}\coloneqq-\beta_{U}(1-\mathrm{e}^{-\tau}),\qquad \bar{x}\coloneqq\beta_{L}(1-\mathrm{e}^{-\tau}),
\end{equation}
$\tilde{T}$ is the period of the limit cycle~\eqref{tildex} and $\tilde{z}_{j}$, with $j\in\mathbb{N}_{>0}\coloneqq\{x\in\mathbb{N} \mid x>0\}$ represents the zeros of $\tilde{x}(t)$.
For $j>2$ we have $\tilde{z}_{j}=\tilde{z}_{j-2}+\tilde{T}$. The period $\tilde{T}$ and the zeros $\tilde{z}_{1}$, $\tilde{z}_{2}$ are given by~\eqref{z1z2T}~\citep{Mackey_2017}.
\begin{equation}\label{z1z2T}
\tilde{z}_{1}\coloneqq\ln\frac{\beta_{L}-\ubar{x}}{\beta_{L}},\qquad \tilde{z}_{2}\coloneqq\tilde{z}_{1}+\tau+\ln\frac{\beta_{U}+\bar{x}}{\beta_{U}}, \qquad  \tilde{T}\coloneqq\tilde{z}_{2}+\tau.
\end{equation}
At $t=0$ where the limit cycle has a minimum point  and at the maximum of the limit cycle when $t_{max}\coloneqq\tilde{z}_{1}+\tau$ the derivative of $\tilde{x}(t)$ is undefined. From~\eqref{z1z2T} it is clear that $\tilde{z}_{j+1}>\tilde{z}_{j}+\tau$ for all $j\in\mathbb{N}_{>0}$.

For every history function $\varphi$ restricted to the set $\varphi_{0}$ the solution $x(t,\varphi_{0})$ of~\eqref{dxdt}-\eqref{ftau} converges to the limit cycle~\eqref{tildex} within a finite time and with a new phase, and this limit cycle is stable~\citep[Theorem 3.3]{Mackey_2017}. The DDE~\eqref{dxdt}-\eqref{ftau} admits infinitely many other periodic orbits, but are all unstable~\citep[Remark 3.2 and 3.3]{Mackey_2017}.
%
\section{Single Stimulus}\label{sec.single}
In this section we consider the response of Eq.~\eqref{dxdt} to a single pulse perturbation.
We summarize some notation from~\citet{Mackey_2017}, 
defining the \textit{resetting time} for the pulse-like perturbation and distinguishing it from the concept of \textit{cycle length map} defined in~\citet{Mackey_2017}. We show that the resetting time is always less than the cycle length map, and that the minimum resetting time is equal to or greater than the stimulus duration. We also show how a single perturbation can lead to an infinite resetting time because of an unstable limit cycle.  The influence of the amplitude, phase and time duration of the perturbation on the cycle length map is also examined.

\par Consider a single pulse-like perturbation of amplitude $a > 0$ and duration $\sigma\in (0,\tau]$ which starts at $t=\Delta\in [0,\tilde{T})$, where $\tilde{T}$ is the period of the periodic solution of Eq.~\eqref{dxdt} with $f(x(t-\tau))$ given by the discontinuous function~\eqref{ftau}. For $\Delta\leq t\leq\Delta+\sigma$  Eq.~\eqref{dxdt} becomes
\begin{equation}\label{eq.dde.pert}
x^{\prime}(t) = -x(t) + f(x(t-\tau)) + a.
\end{equation}
\citet{Mackey_2017} examined the response of the limit cycle~\eqref{tildex} to a single pulse-like stimulus with positive amplitude $a$, as defined in Eq.~\eqref{eq.dde.pert}. They calculated the minima, maxima and period of the perturbed limit cycle for different values of the starting time and duration of the single pulse-like perturbation. (The single stimulus of amplitude $a$ and duration $\sigma$ can be related to the dose and temporal duration of the administration of exogenous cytokines in an attempt to regulate the peripheral blood cells~\citep{Mackey_2017}).

The response of the limit cycle solution~\eqref{tildex} to a perturbation with amplitude $a$ can be calculated piecewise by considering the phase which the stimulus begins $\Delta\in [0,\tilde{T})$ and ends $\Delta+\sigma$. We distinguish the perturbation for different values of $\{\Delta,\sigma,a\}$ as in~\citet{Mackey_2017}.

To use the same notation, we classify the stimulus distinguishing whether it begins in a rising phase \textbf{R} ($\dot{x}>0$), falling phase \textbf{F} ($\dot{x}\leq 0$), with negative value \textbf{N} ($x<0$), or non-negative value \textbf{P} ($x\geq 0$); and whether it ends in a phase \textbf{R} or \textbf{F}, with value \textbf{N} or \textbf{P}. Precisely, including the points where the derivative is undefined, we use the nomenclature \textbf{R} for $\Delta\in[0,t_{max})$, \textbf{F} for $\Delta\in[t_{max},\tilde{T})$, \textbf{P} for $\Delta\in[\tilde{z}_{1},\tilde{z}_{2}]$, and \textbf{N} for $\Delta\in[0,\tilde{z_{1}})\cup (\tilde{z}_{2},\tilde{T}]$. For example, the sequence of letters \textbf{RPFN} denotes a pulse which {\it begins} in the {\it rising} phase with {\it positive} value and {\it ends} in a {\it falling} phase with {\it negative} value.  We also denote the solution of the DDE with perturbation~\eqref{eq.dde.pert} by $x^{(\Delta)}(t)$, and its zeros by $z_{\Delta,j}$, with $j\in\mathbb{N}_{>0}$, as in~\citet{Mackey_2017}. The zeros of $x^{(\Delta)}(t)$ form an increasing sequence with $z_{\Delta,j}+\tau<z_{\Delta,j+1}\quad\mbox{for all}\quad j\in\mathbb{N}_{>0}$~\citep{Mackey_2017}. With this classification we have the following sequence of possibilities~\citep{Mackey_2017}:
\begin{equation}\label{cases}
\left\{\begin{array}{ll}
\textbf{RNRN}, \textbf{RNRP}, \textbf{RPRP}, \textbf{RPFP}, \textbf{RPFN}, \\
\textbf{FPFP}, \textbf{FPFN}, \textbf{FNFP}, \textbf{FNFN}, \textbf{FNRN}, \textbf{FNRP}.
\end{array}\right.
\end{equation}

Each subcase~\eqref{cases} is defined on a $\Delta$-subinterval $\mathit{I}\subset[0,\tilde{T})$. For each parameter vector $(\tau,\beta_{U},\beta_{L},\sigma,a)$ there is a set of subcases~\eqref{cases}.  The union of all their $\Delta$-subintervals is equal to $[0,\tilde{T})$ and their intersection is the empty set. We denote the $\Delta$-subintervals as $\mathit{I}$ with a subscript with the respective sequence of letters, e.g. for the case \textbf{RNRN} the $\Delta$-subinterval is denoted by $\mathit{I_{RNRN}}$.

For the subcase \textbf{FNFP}, the time required for the perturbed solution $x^{(\Delta)}(t)$ to return to the limit cycle $\tilde{x}$ may not be finite everywhere due to a rapidly oscillating periodic solution~\citep{Mackey_2017}. This is the most complex case in~\eqref{cases} and will be analyzed here in detail. In this case the pulse starts in the falling phase of the limit cycle $\tilde{x}$ with $x^{(\Delta)}(\Delta)<0$, which implies $\Delta\in(\tilde{z}_{2},\tilde{T})$. The pulse also ends in the falling phase, i.e. $\Delta\in(t_{max},\tilde{T}-\sigma)$, and with positive value $x^{(\Delta)}(\Delta+\sigma)>0$. The positivity condition $x^{(\Delta)}(\Delta+\sigma)>0$ implies $\Delta<\delta_{2}$, where $\delta_{2}$ is a constant defined by $x^{(\Delta)}(\delta_{2}+\sigma)=0$, which yields
\begin{equation}\label{delta2}
\delta_{2} \coloneqq \tilde{z}_{2} - \sigma + \ln\left(\frac{\beta_{U}}{\beta_{U}-a(1-\mathrm{e}^{-\sigma})}\right) \quad\mbox{for}\quad \beta_{U}>a(1-\mathrm{e}^{-\sigma}),
\end{equation}
the same constant $\delta_{2}$ defined in~\citet[Eq. (5.11)]{Mackey_2017}. Thus the $\Delta$-subinterval in the case \textbf{FNFP} is of the form
\begin{equation}\nonumber
\mathit{I_{FNFP}} = (\tilde{z}_{2},\tilde{T}-\sigma)\cap (-\infty,\delta_{2}).
\end{equation}

\begin{figure}[t]
\centering
\includegraphics[width=1.0\textwidth]{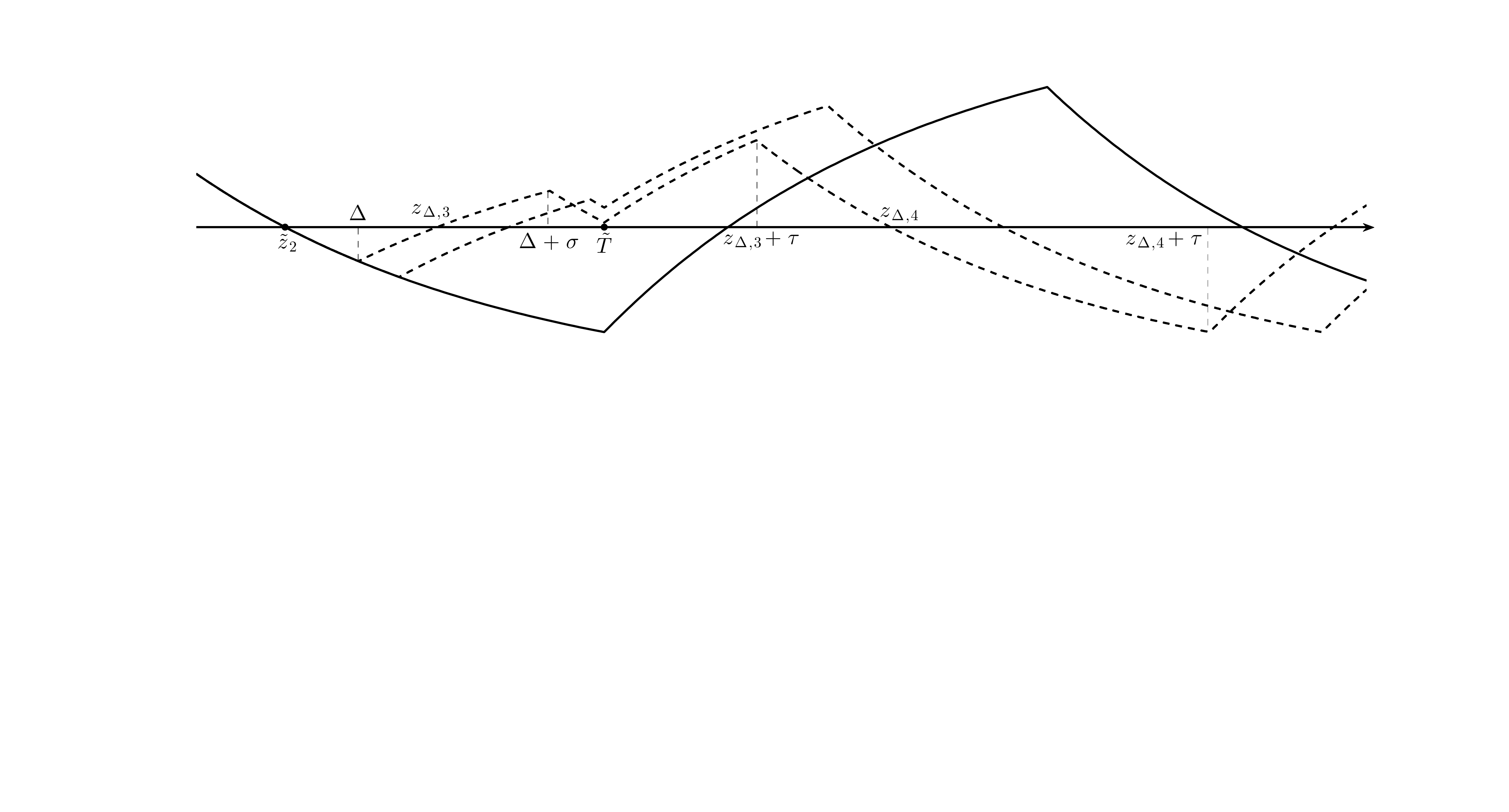}
\caption{A schematic representation of the unperturbed limit cycle~\eqref{dxdt} (solid line) and two perturbed solutions (dashed lines) for the case \textbf{FNFP1}. Here the parameters are $\tau=1$, $\sigma=0.6$, $\beta_{U}=0.3$, $\beta_{L}=0.4$, $a=0.52$, and $\Delta\in\{2.23, 2.36\}$.}
\label{fig_FNFP2}
\end{figure}

\begin{figure}[t]
\centering
\includegraphics[width=1.0\textwidth]{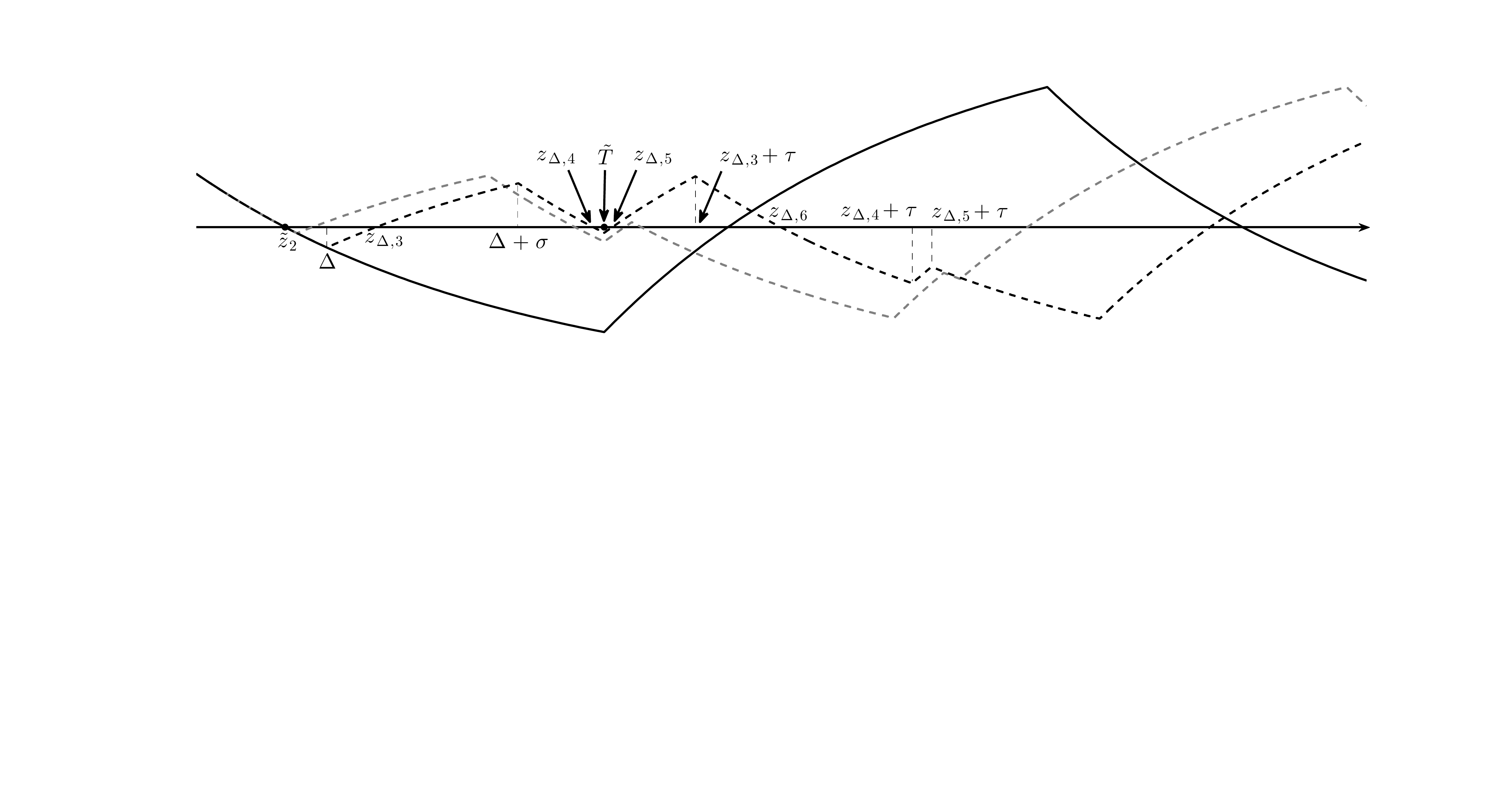}
\caption{As in Figure~\ref{fig_FNFP2}, but for the case \textbf{FNFP2} and \textbf{FNFP4} with $\Delta\in\{2.04, 2.13\}$.}
\label{fig_FNFP1}
\end{figure}

In Remark~\ref{rem.FNFP} we show that the inequality $a>\beta_{U}$ holds for the case \textbf{FNFP} and we distinguish  between two subcases, $\textbf{FNFP1}$ and $\textbf{FNFP2}$. For the case $\textbf{FNFP1}$ we have $x^{(\Delta)}(\tilde{T})>0$ as in the perturbed solutions shown in Figure~\ref{fig_FNFP2}, while for the case $\textbf{FNFP2}$ we have $x^{(\Delta)}(\tilde{T})\leq 0$ as in the examples displayed in Figure~\ref{fig_FNFP1}.
\begin{remark}\label{rem.FNFP}
\textit{%
The inequality $a>\beta_{U}$ holds for the case \textbf{FNFP} and we can distinguish two subcases by incorporating further conditions as follows:
\begin{itemize}
\item\textbf{FNFP1} If $x^{(\Delta)}(\tilde{T})\geq 0$, then
\begin{equation}\label{FNFP1}
\Delta\in(\tilde{z}_{2},\tilde{T}-\sigma)\cap(-\infty,\delta_{2})\cap[\delta_{4},\infty)=\mathit{I_{FNFP1}};
\end{equation}
\item\textbf{FNFP2} If $x^{(\Delta)}(\tilde{T})<0$, then
\begin{equation}\label{FNFP2}
\Delta\in(\tilde{z}_{2},\tilde{T}-\sigma)\cap(-\infty,\delta_{2})\cap(-\infty,\delta_{4})=\mathit{I_{FNFP2}};
\end{equation}
\end{itemize}
\noindent where the constant $\delta_{4}$ is such that  $x^{(\Delta)}(\tilde{T})=0$ for $\Delta=\delta_{4}$ and is given by
\begin{equation}\label{delta4}
\delta_{4} = \tilde{z}_{2} + \ln\frac{\beta_{U}(\mathrm{e}^{\tau}-1)}{a(\mathrm{e}^{\sigma}-1)}.
\end{equation}
}%
(Remember that all proofs are presented in Appendix~\ref{sec.app}.)
\end{remark}

\begin{figure}[t]
\centering
\includegraphics[width=1\textwidth]{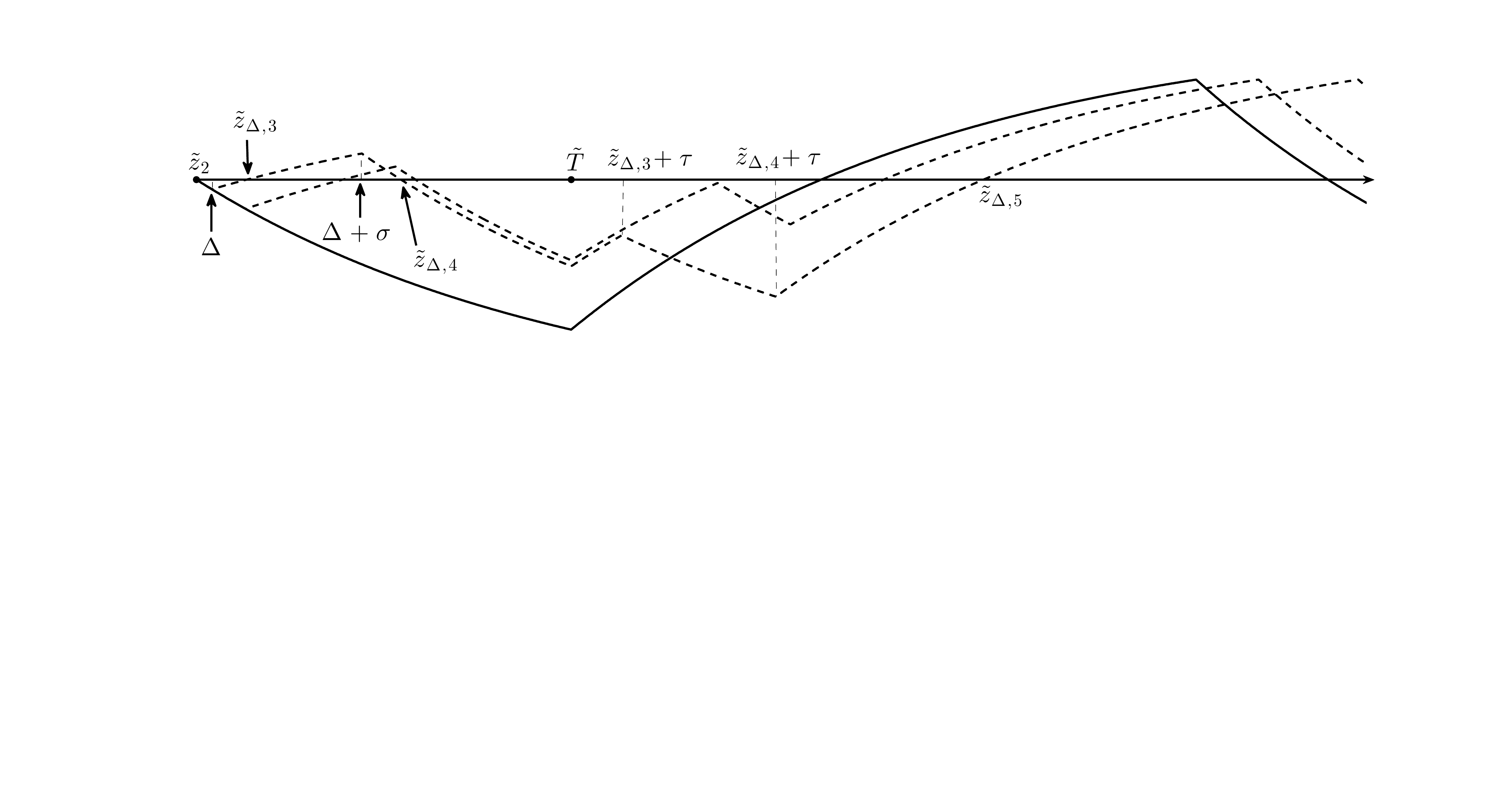}
\caption{A schematic representation of two solutions for Case \textbf{FNFP3}, where the parameters are $\tau=1$, $\sigma=0.4$, $\beta_{U}=0.6$, $\beta_{L}=0.4$, $a=0.85$, and $\Delta\in\{2.06, 2.15\}$.}
\label{fig_FNFP24}
\end{figure}

\begin{figure}[!b]
\centering
\includegraphics[width=1.0\textwidth]{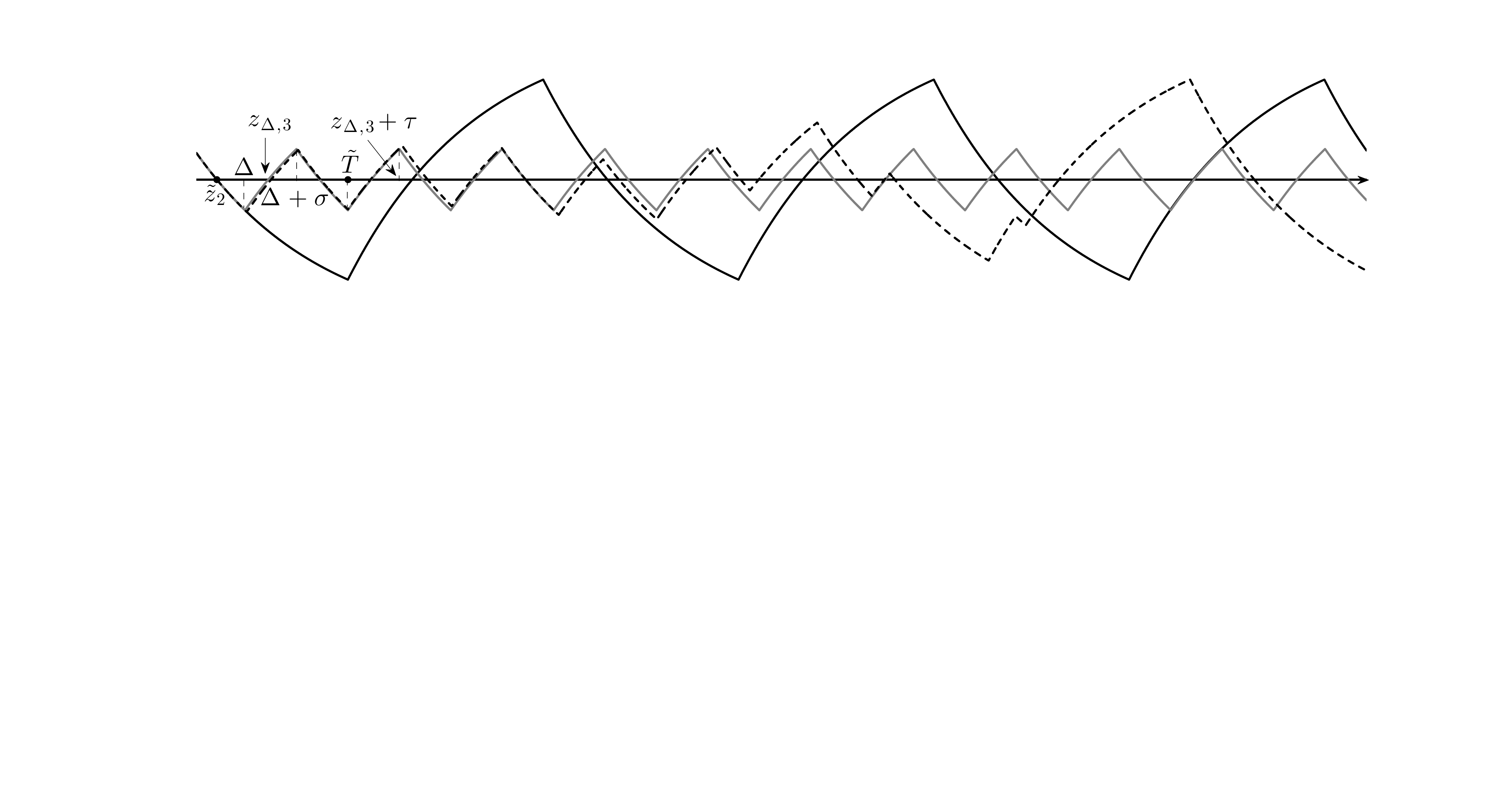}
\caption{A schematic representation of the unperturbed limit cycle~\eqref{dxdt} (solid black line), a rapidly oscillating unstable periodic solution (solid gray line) and a perturbed solution (dashed line) for the case \textbf{FNFP2}. Here the parameters are $\tau=1$, $\sigma=0.39235$, $\beta_{U}=0.4$, $\beta_{L}=0.4$, $a=0.8$, and $\Delta=\{\delta_{\infty},2.208\}$, where $\delta_{\infty}\approx 2.19505$ is given by~\eqref{eq.deltaInf2}.}
\label{fig_FNFP2_23}
\end{figure}

In Remark~\ref{rem.FNFPB} we show that the case $\textbf{FNFP2}$ from Remark~\ref{rem.FNFP} can be distinguished between two subcases by including the extra condition $x^{(\Delta)}(z_{\Delta,3}+\tau)< 0$ to define the case $\textbf{FNFP3}$, and the extra conditions $x^{(\Delta)}(z_{\Delta,3}+\tau)\geq 0$ and $x^{(\Delta)}(z_{\Delta,4}+\tau)<0$ to define $\textbf{FNFP4}$.
The examples of solutions shown in Figure~\ref{fig_FNFP1} correspond to the subcase $\textbf{FNFP2}$ and also $\textbf{FNFP4}$ while the examples of solutions shown in Figure~\ref{fig_FNFP24} refers to case $\textbf{FNFP3}$.
The case $\textbf{FNFP2}$ actually splits into an infinite number of subcases, and examples are shown in Figure~\ref{fig_FNFP2_23}.
The solution represented by the dashed line in Figure~\ref{fig_FNFP2_23} shows how long the transient solution can be before it returns to the periodic orbit.

\begin{remark}\label{rem.FNFPB}
\textit{%
The case \textbf{FNFP2}, described in Remark~\ref{rem.FNFP}, can be further divided into two subcases by including extra conditions as follows:
\begin{itemize}
\item\textbf{FNFP3} If $x^{(\Delta)}(z_{\Delta,3}+\tau)< 0$, then
\begin{equation}\label{FNFP3}
\Delta\in(\tilde{z}_{2},\tilde{z}_{2}\!+\!\tau\!-\!\sigma)\!\cap\!(-\infty,\delta_{2})\!\cap\!(-\infty,\delta_{4})\!\cap\!(-\infty,\hat{\delta}_{4})=\mathit{I_{FNFP3}};
\end{equation}
\item\textbf{FNFP4} If $x^{(\Delta)}(z_{\Delta,3}+\tau)\geq 0$ and $x^{(\Delta)}(z_{\Delta,4}+\tau)<0$, then
\begin{equation}\label{FNFP4}
\Delta\in(\tilde{z}_{2},\tilde{z}_{2}\!+\!\tau\!-\!\sigma)\!\cap\!(-\infty,\delta_{2})\!\cap\!(-\infty,\delta_{4})\!\cap\![\hat{\delta}_{4},\infty)\!\cap\!(-\infty,\delta_{5})=\mathit{I_{FNFP4}};
\end{equation}
\end{itemize}
\noindent where the constants $\delta_{2}$, $\delta_{4}$, $\hat{\delta}_{4}$ and $\delta_{5}$ are respectively given by~\eqref{delta2},~\eqref{delta4},~\eqref{eq.delta4hat} and~\eqref{eq.delta5}.
}%
\end{remark}

Between the case shown in Figure~\ref{fig_FNFP1} and the rapid limit cycle shown in Figure~\ref{fig_FNFP2_23} there is a sequence of subcases. For each new case the solution oscillates one  more time before approaching the limit cycle. For all solutions of the cases \textbf{FNFP} it is expected that $\ubar{x}\leq x^{(\Delta)}(t)\leq \bar{x}$.

For a perturbation with starting time $\Delta$ we define the {\it resetting time} $F(\Delta)$ as the time interval required for the perturbed solution $x^{(\Delta)}(t)$ to return to the limit cycle $\tilde{x}(t)$ with a new phase. We formally define the resetting time in Definition~\ref{def.F} and give it for each case~\eqref{cases} in Remark~\ref{rem.calcF}.

\begin{definition}\label{def.F}
{\bf Resetting time.}  Let us define the function  $F: [0,\tilde{T}) \ni \Delta\mapsto F(\Delta)\in \mathbb{R}$. Assume that the limit cycle $\tilde{x}(t)$ is perturbed at time $\Delta$, as in~\eqref{eq.dde.pert},  and that after a time interval $F(\Delta)$, defined as the resetting time, the perturbed solution returns to the limit cycle with a new phase $\phi$ and stays on it, i.e. $x^{(\Delta)}(t)=\tilde{x}(t+\phi)$ holds for all $t\geq \Delta + F(\Delta)$.
\end{definition}

\begin{remark}\label{rem.calcF}
Assume that the limit cycle $\tilde{x}(t)$ is perturbed at time $\Delta$ as in~\eqref{eq.dde.pert} and $\delta_{4}<\tilde{z}_{2}$, i.e. $\mathit{I_{FNFP2}}=\emptyset$. For each of the cases~\eqref{cases} the resetting time $F(\Delta)$ is given, respectively,  by:
\begin{enumerate}
\item $F_{RNRN}(\Delta)=\sigma$;
\item $F_{RNRP}(\Delta)=t_{max}+(z_{\Delta,2}-\tilde{z}_{2})-\Delta$;
\item $F_{RPRP}(\Delta)=t_{max}+(z_{\Delta,2}-\tilde{z}_{2})-\Delta$;
\item $F_{RPFP}(\Delta)=\sigma$;
\item $F_{RPFN}(\Delta)=z_{\Delta,2}+\tau-\Delta$;
\item $F_{FPFP}(\Delta)=\sigma$;
\item $F_{FPFN}(\Delta)=z_{\Delta,2}+\tau-\Delta$;
\item $F_{FNFP1}(\Delta)=z_{\Delta,3}+\tau-\Delta$;
\item $F_{FNFN}(\Delta)=\tilde{z}_{2}+\tau-\Delta$;
\item $F_{FNRN}(\Delta)=\sigma$;
\item $F_{FNRP}(\Delta)=\tilde{z}_{3}+\tau-(\tilde{z}_{4}-z_{\Delta,4})-\Delta$.
\end{enumerate}
\end{remark}

By definition, the resetting time $F(\Delta)$ is the minimum time required for the perturbed solution to return to the limit cycle, while the cycle length map $T(\Delta)$ is the time measured between two marker events~\citep{Glass_1988}, one located before the solution is perturbed and the other after the perturbed solution has returned to the limit cycle.  The cycle length map is given by $T(\Delta)=\tilde{T}+\tilde{T}\phi(t)$, where $\phi(t)$ is the phase difference between the two marker events. In~\citet{Mackey_2017} the zeros of the limit cycle, $\tilde{z}_{j}$, and the zeros of the perturbed solution, $z_{\Delta,j}$, were used as marker events to calculate the cycle length map $T(\Delta)$. Thus, by definition, the resetting time is always less than the cycle length map, c.f. Remark~\ref{rem.FT}.
\begin{remark}\label{rem.FT}
\textit{For all $\Delta\in[0,\tilde{T})$ with $\delta_{4}<\tilde{z}_{2}$ ($\mathit{I_{FNFP2}}=\emptyset$) we have $F(\Delta)<T(\Delta)$.}
\end{remark}

In the next section we consider~\eqref{eq.dde.pert} with a periodic perturbation instead of a single perturbation. Thus it will be of interest to know whether $T(\Delta)$ and $F(\Delta)$ are finite,  and to know their lower and upper bounds.
We define the maximum and minimum of the resetting time and cycle length map respectively by
\begin{equation}\nonumber
\bar{F}\coloneqq\!\max\limits_{\Delta\in[0,\tilde{T})}\!F(\Delta),\ \  \ubar{F}\coloneqq\!\min\limits_{\Delta\in[0,\tilde{T})}\!F(\Delta),\ \  \bar{T}\coloneqq\!\max\limits_{\Delta\in[0,\tilde{T})}\!T(\Delta), \ \ \mbox{and}\quad \ubar{T}\coloneqq\!\min\limits_{\Delta\in[0,\tilde{T})}\!T(\Delta).
\end{equation}
As expected, for all $\Delta\in[0,\tilde{T})$ and $\sigma\in(0,\tau]$ the minimum resetting time is equal to or greater than the stimulus duration $\sigma$, as detailed in Remark~\ref{rem.min.ress}.  In Remark~\ref{rem.max.ress} the maximum of $T(\Delta)$ for the case $a<\beta_{U}$ is investigated.
\begin{remark}\label{rem.min.ress}
\textit{For all $\Delta\in[0,\tilde{T})$ and $\sigma\in(0,\tau]$ with $\delta_{4}<\tilde{z}_{2}$ ($\mathit{I_{FNFP2}}=\emptyset$) we have $\ubar{F}=\sigma$.}
\end{remark}
\begin{remark}\label{rem.RPFN}
\textit{For $\delta_{2}>t_{max}$ the inequality $a<\beta_{U}$ holds and $\mathit{I_{RPFN}}=\emptyset$.}
\end{remark}
\begin{remark}\label{rem.max.ress}
\textit{For $a<\beta_{U}$ the maximum of the cycle length map occurs at $\Delta=\delta_{2}$ and is given by}
\begin{equation}\label{Tmax1}
\bar{T} = \tilde{T} +  \ln\left(\frac{\beta_{U}}{\beta_{U}-a(1-\mathrm{e}^{-\sigma})}\right) > \tilde{T}.
\end{equation}
\end{remark}

Remark~\ref{rem.RPFN}, and the fact that the cycle length map $T(\Delta)$ is continuous for $a<\beta_{U}$~\citep[Corollary 4.2]{Mackey_2017}, allows us to compute a maximum of $T(\Delta)$ as stated in Remark~\ref{rem.max.ress}.

As consequence of Remarks~\ref{rem.FT} and~\ref{rem.min.ress}, $\sigma$ is a lower bound for the cycle length map, i.e. $T(\Delta)>\sigma$.
From Remarks~\ref{rem.max.ress} and~\ref{rem.FT} it follows that for $a<\beta_{U}$ the maximum $\bar{T}$~\eqref{Tmax1} is an upper bound for the resetting time, i.e. $F(\Delta)<\bar{T}$.

For $a>\beta_{U}$ the resetting time and cycle length map may not be bounded above. For all of the cases  in~\eqref{cases} the case \textbf{FNFP} is the only one in which the cycle length map may not be finite everywhere~\citep{Mackey_2017}. In Figure~\ref{fig_FNFP2_23} we show perturbed solutions that illustrate how a single perturbation to \textbf{FNFP} can lead to an infinite resetting time due an unstable limit cycle. In Remark~\ref{rem.deltainf} we show that for $\Delta$ equal to~\eqref{eq.deltaInf2} the perturbed solution $x^{(\Delta)}(t)$ settles down on a rapidly oscillating unstable periodic solution with $\ubar{x}<x^{(\Delta)}(t)<\bar{x}$ whose period $(\tilde{T}^{(\infty)})$ satisfies $\tau-\sigma<\tilde{T}^{(\infty)}<\tau$.

%
%
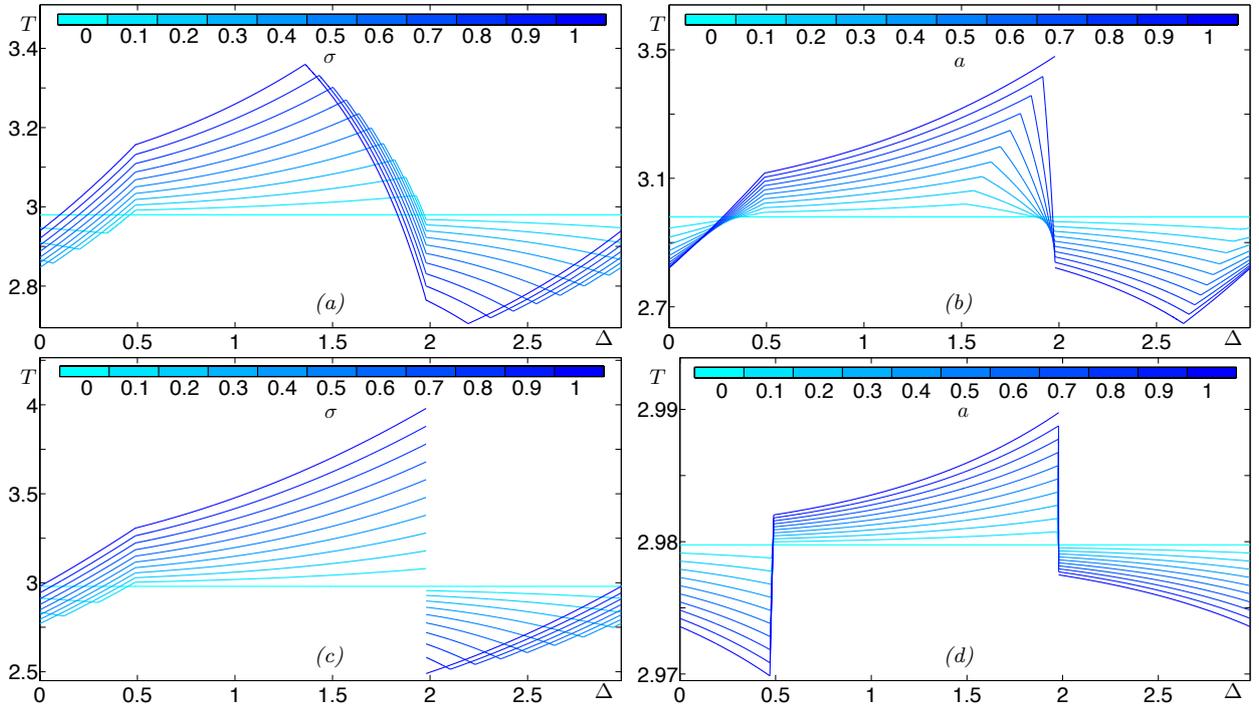
\begin{figure}[!b]
\includegraphics[width=0.49\textwidth,height=0.278\textwidth]{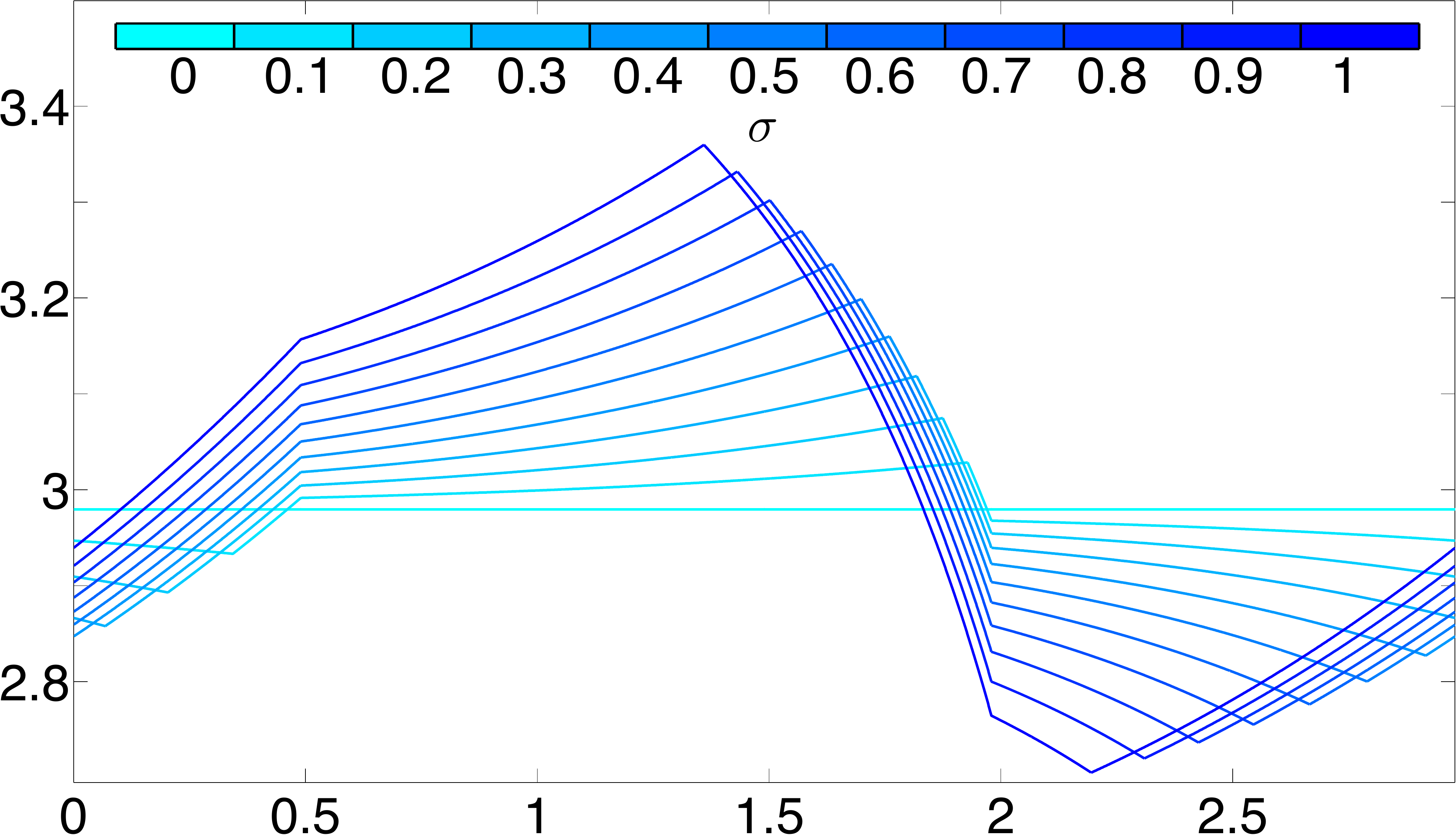}
\put(-115,15){\footnotesize{\textit{(a)}}}
\put(-225,120){\footnotesize{$T$}}
\put(-10,1.5){\footnotesize{$\Delta$}}
\hspace*{1mm}
\includegraphics[width=0.49\textwidth,height=0.278\textwidth]{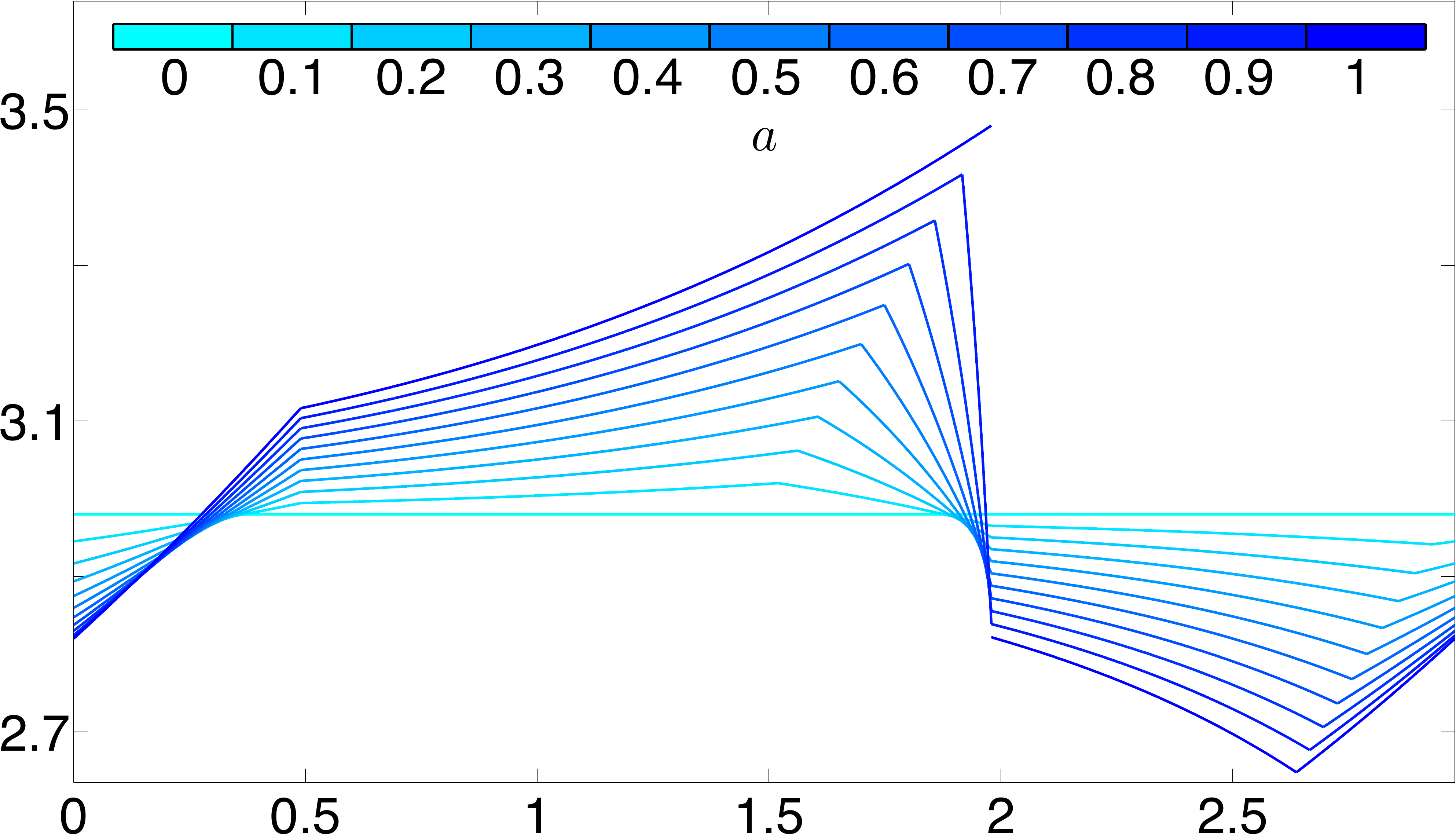}
\put(-115,15){\footnotesize{\textit{(b)}}}
\put(-225,120){\footnotesize{$T$}}
\put(-10,1.5){\footnotesize{$\Delta$}}\\

\vspace*{-3.5mm}
\includegraphics[width=0.49\textwidth,height=0.278\textwidth]{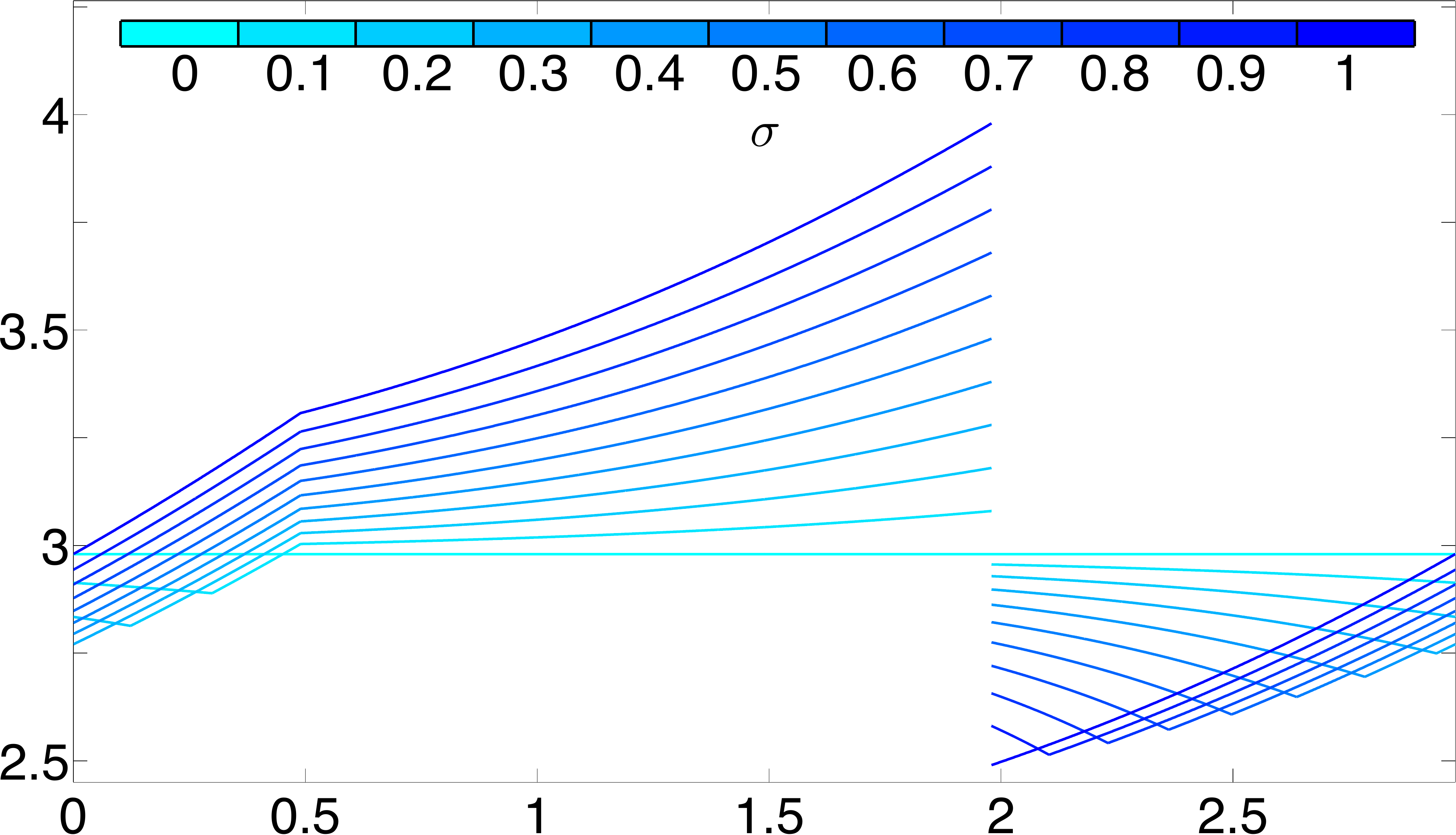}
\put(-115,15){\footnotesize{\textit{(c)}}}
\put(-225,120){\footnotesize{$T$}}
\put(-10,1.5){\footnotesize{$\Delta$}}
\hspace*{1mm}
\includegraphics[width=0.49\textwidth,height=0.278\textwidth]{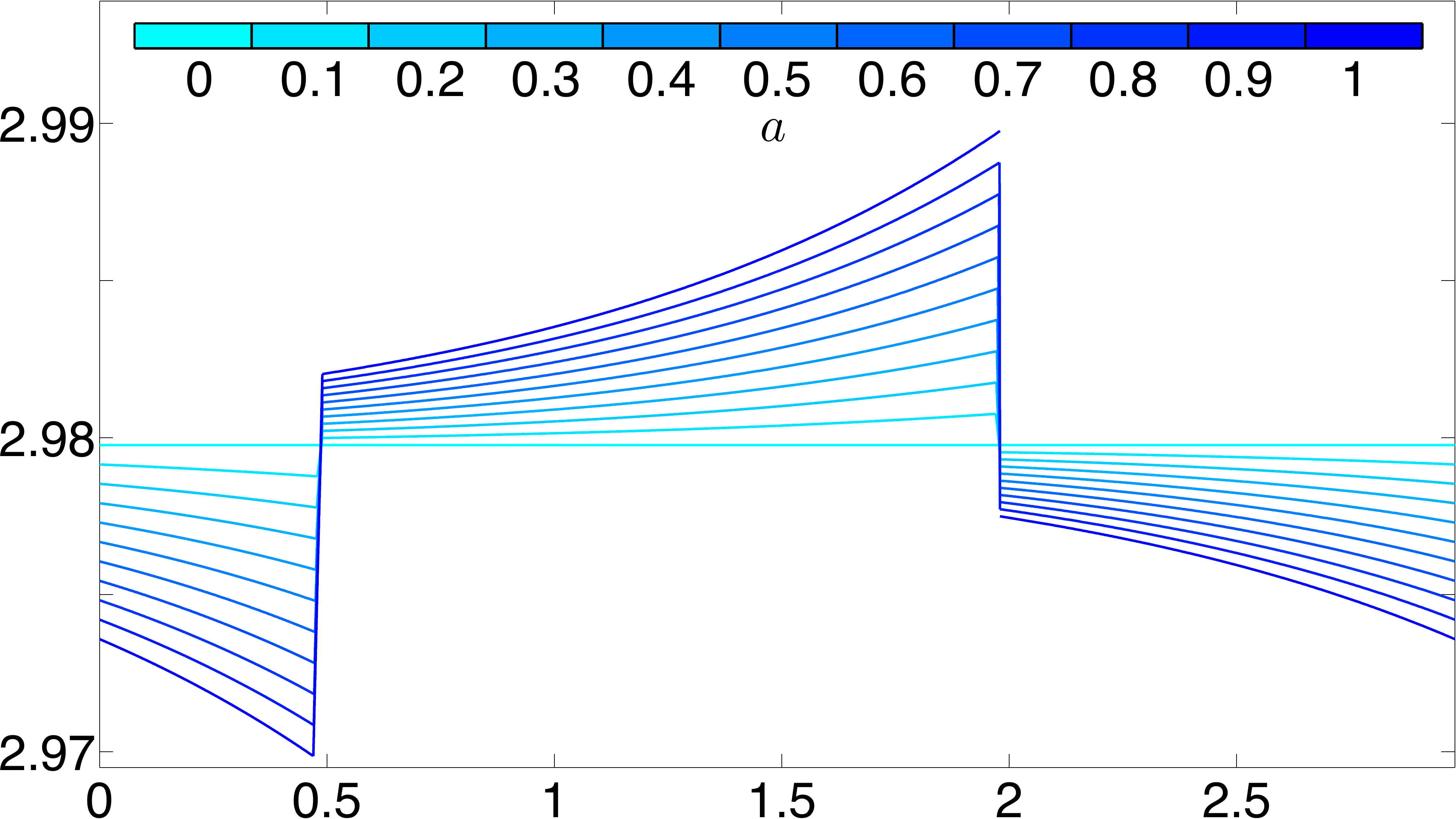}
\put(-115,15){\footnotesize{\textit{(d)}}}
\put(-225,120){\footnotesize{$T$}}
\put(-10,1.5){\footnotesize{$\Delta$}}
\vspace*{-02mm}
\caption{Graphs of the cycle length map $[0,\tilde{T}]\ni\Delta\mapsto T(\Delta)$ for $\beta_{U}=\beta_{L}=\tau=1$ and hence $\tilde{T}=2.97976$, $\tilde{z}_{1}=0.48988$, $\tilde{z}_{2}=1.97976$ with 11 curves transitioning from cyan to dark blue respectively for $\sigma=\{0,0.1,\ldots, 0.9,1\}$ in panels (a) and (c) and $a=\{0,0.1,\ldots,1\}$ in panels (b) and (d) with: $a=0.5$ for panel (a), $\sigma=0.5$ for panel (b), $a=1$ for panel (c) and $\sigma=0.01$ for panel (d).}
\label{fig.clm}
\end{figure}
\begin{remark}\label{rem.deltainf}
\textit{For the case \textbf{FNFP} the cycle length map and the resetting time both tend to infinity when $\Delta$ tends to the constant $\delta_{\infty}$ which is given by
\begin{equation}\label{eq.deltaInf2}
\delta_{\infty} = \tilde{z}_{2} + \sigma + \ln{\frac{a-\beta_{U}(1-\mathrm{e}^{-\sigma})}{a}},
\end{equation}
and satisfies $\tilde{z}_{2}<\delta_{\infty}<\tilde{z}_{2} + \sigma$. For $\Delta=\delta_{\infty}$ the perturbed solution $x^{(\Delta)}(t)$ settles down on a rapidly oscillating unstable periodic solution for which the period $(\tilde{T}^{(\infty)})$ satisfies $\tau-\sigma<\tilde{T}^{(\infty)}<\tau$ and is given by $\tilde{T}^{(\infty)}=\tilde{T}-\delta_{\infty}$. Moreover, the minimum $(\ubar{x}^{(\Delta)})$ and maximum $(\bar{x}^{(\Delta)})$ of the rapid limit cycle are given by
\begin{equation}\nonumber
\left\{\begin{array}{ll}
\ubar{x}^{(\Delta)} = \beta_{U}(\mathrm{e}^{-(\delta_{\infty}-\tilde{z}_{2})} - 1), \\[2mm]
\bar{x}^{(\Delta)} =  a(1-\mathrm{e}^{-\sigma}) + \beta_{U}(\mathrm{e}^{\tilde{z}_{2}-\delta_{\infty}-\sigma} - 1),
\end{array}\right.
\end{equation}
and satisfy $\ubar{x}^{(\Delta)}>\ubar{x}$ and $\bar{x}^{(\Delta)}<\bar{x}$.}
\end{remark}

Next we investigate how the cycle length map $T(\Delta)$ computed in~\citet{Mackey_2017} varies with changes in the  perturbation amplitude $a$ and the pulse duration $\sigma$.

Figures~\ref{fig.clm}(a) and (c) show how $T(\Delta)$ changes as function of $\sigma$, while panels (b) and (d) show how $T(\Delta)$ changes when $a$ is varied. For all examples of the four panels $T(0)=T(\tilde{T})$. All curves $T(\Delta)$ of panel (a) are continuous on $[0,\tilde{T}]$. In panels (b) and (d) we see that for $a=\beta_{U}$ a Type 1 discontinuity~\citep{Winfree_1980} appears in the cycle length map. All curves $T(\Delta)$ show this discontinuity in panel (c).
In the limit $\sigma\to 0$ we have $\delta_{1}\to\tilde{z}_{1}$, $\delta_{2}\to\tilde{z}_{2}$, $\Delta\in\mathit{I_{RNRP}}=[\max\{0,\delta_{1}\},\tilde{z}_{1})\to\emptyset$ and $\Delta\in\mathit{I_{FPFN}}=[t_{max},\tilde{z}_{2})\cap(\delta_{2},\infty)\to\emptyset$. The decreasing length of both intervals $\mathit{I_{RNRP}}$ and $\mathit{I_{FPFN}}$ as $\sigma$ is decreased can be seen by comparing the curves with $a<\beta_{u}$ from panels (b) and (d). For panel (b) $\sigma=0.5$ while for panel (d) $\sigma=0.01$ and the intervals $\mathit{I_{RNRP}}$ and $\mathit{I_{FPFN}}$ approach vertical lines for $\Delta\approx \tilde{z}_{1}$ and $\Delta\approx \tilde{z}_{2}$, respectively.

In Figure~\ref{fig.clm} the intersections points $\Delta^{*}$ where $T(\Delta^{*})=\tilde{T}$ are so-called fixed points of the cycle length map and they are unstable if $T^{\prime}(\Delta^{*})>0$ and stable if $T^{\prime}(\Delta^{*})<0$~\citep{GRANADA20091}, 
where $T^{\prime}(\Delta^{*})$ denotes the derivative of $T(\Delta)$ with respect to $\Delta$ evaluated at $\Delta^{*}$.

\begin{figure}[!t]
\centering
\includegraphics[width=0.49\textwidth,height=0.278\textwidth]{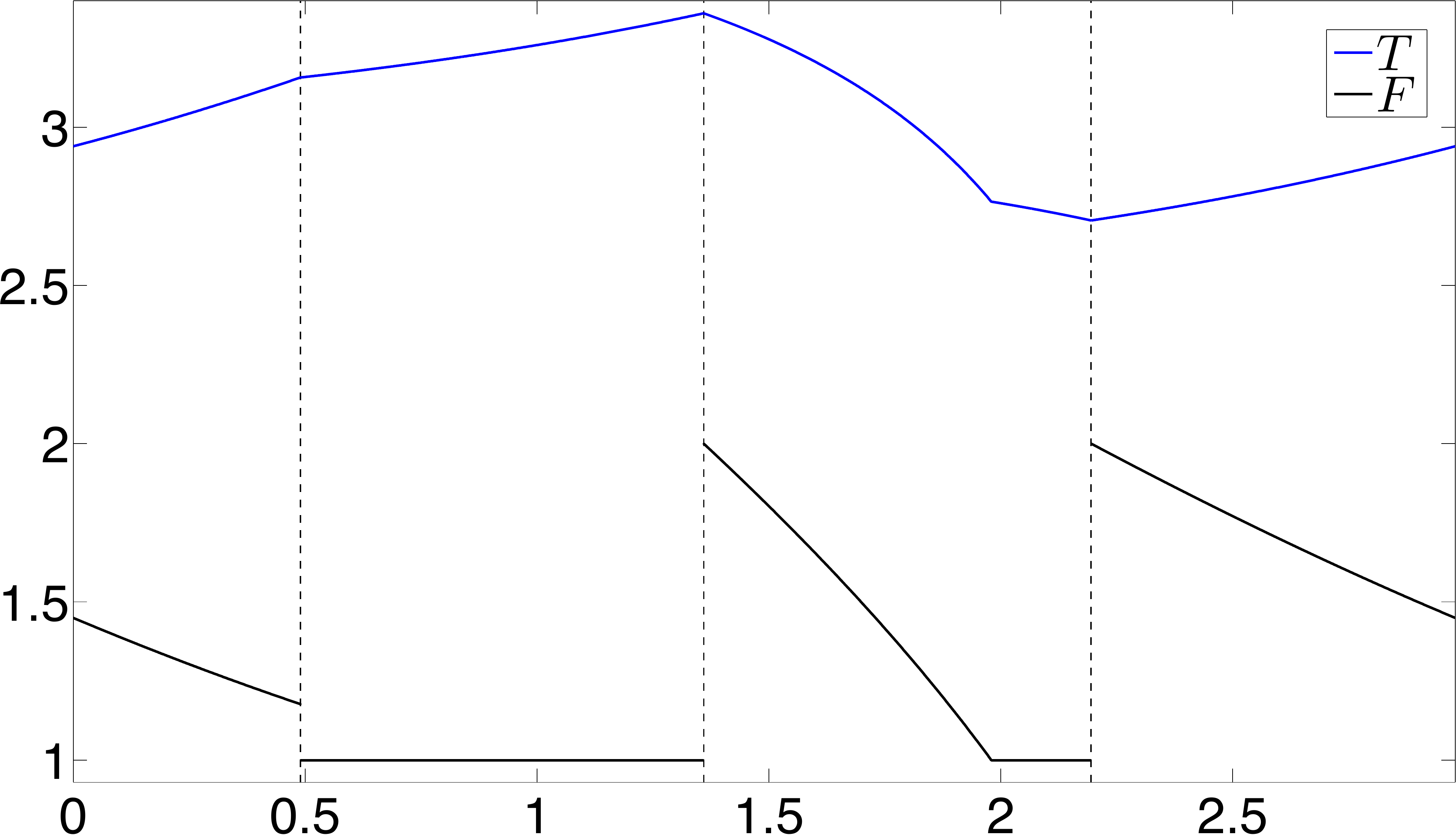}
\put(-214,120){\footnotesize{\textit{(a)}}}
\put(-10,0.5){\footnotesize{$\Delta$}}
\hspace*{0mm}
\includegraphics[width=0.49\textwidth,height=0.278\textwidth]{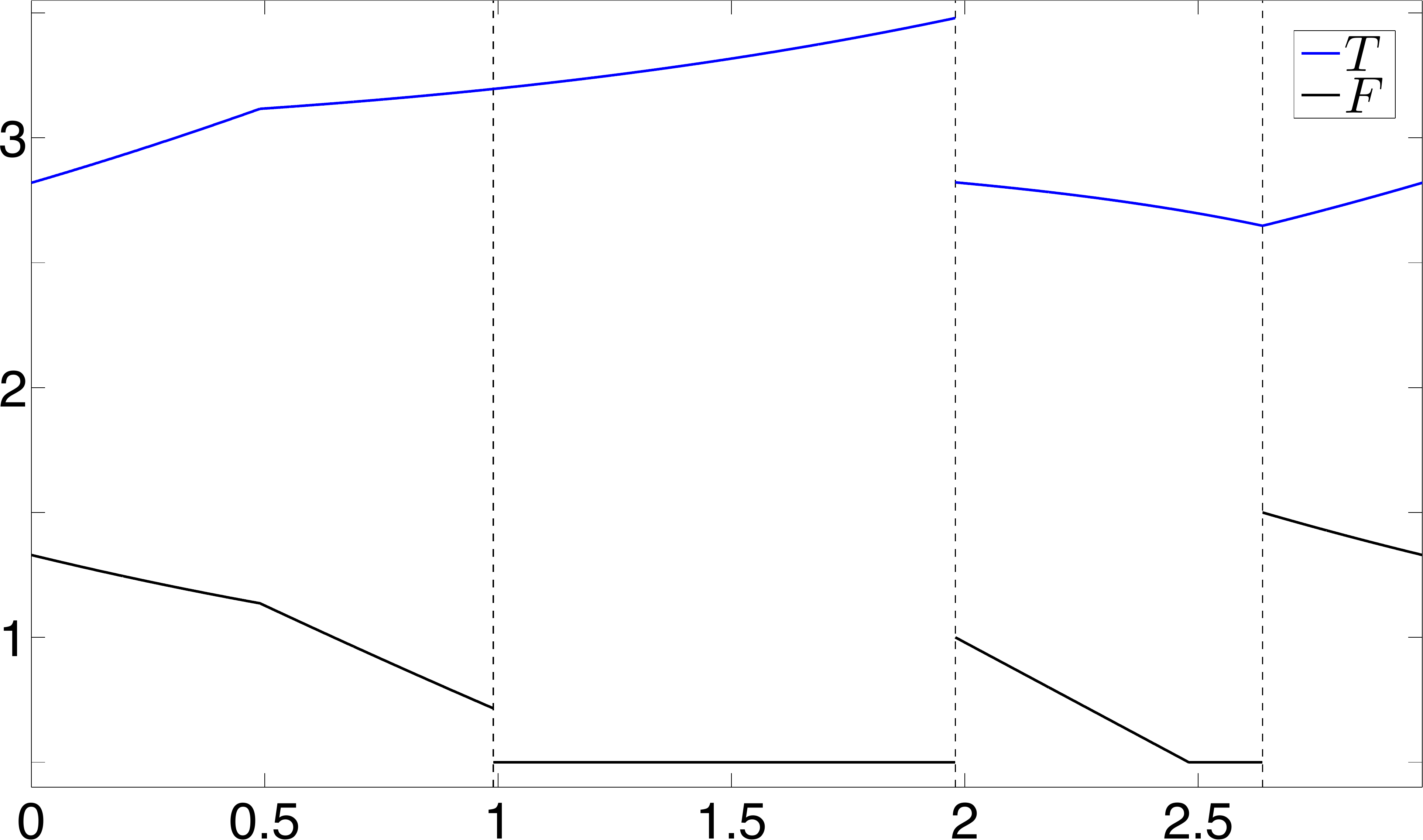}
\put(-220,120){\footnotesize{\textit{(b)}}}
\put(-10,1.5){\footnotesize{$\Delta$}}
\vspace*{-02mm}
\caption{Cycle length map $[0,\tilde{T}]\ni\Delta\mapsto T(\Delta)$ and time resetting $[0,\tilde{T}]\ni\Delta\mapsto F(\Delta)$ for $\beta_{U}=\beta_{L}=\tau=1$ and hence $\tilde{T}=2.97976$.  For panel (a) $a=0.5$, $\sigma=\tau$ and for panel (b) $\sigma=0.5$ and $a=\beta_{U}$. The vertical dashed lines indicate the points where $F(\Delta)$ is discontinuous.}
\label{fig_FT}
\end{figure}

For $a<\beta_{U}$ the cycle length map is continuous~\citep[Corollary 4.2]{Mackey_2017}. Figure~\ref{fig_FT} shows the cycle length map $T(\Delta)$ and time resetting $F(\Delta)$ for the examples with dark blue lines from Figure~\ref{fig.clm}(a) and (b). In panel (a) we have $a<\beta_{U}$, $T(\Delta)$ is continuous and the vertical dashed lines indicate the points where $F(\Delta)$ is discontinuous. From left to right the lines respectively correspond to $\Delta$ equal to: $z_{1}=t_{max}-\sigma=0.48988$, $\delta_{2}=1.35965$ and $\tilde{T}+\delta_{1}=2.19500$. The graphs correspond to the following sequence of cases: \textbf{RNRP}, \textbf{RPRP}, \textbf{RPFP}, \textbf{RPFN}, \textbf{FPFN}, \textbf{FNFN}, \textbf{FNRN}, \textbf{FNRP}. In   panel (b) we have $a=\beta_{U}$, $T(\Delta)$ has a discontinuity at $\delta_{2}=z_{2}$ and the vertical dashed lines indicate the points where $F(\Delta)$ is discontinuous. From left to right the lines are respectively given by $\Delta$ equal to: $t_{max}-\sigma= 0.98988$, $\delta_{2}=z_{2}=1.97976$ and $\tilde{T}+\delta_{1}=2.63784$. For this example we have the following sequence of cases: \textbf{RNRP}, \textbf{RPRP}, \textbf{RPFP}, \textbf{FPFP}, \textbf{FPFN}, \textbf{FNFN}, \textbf{FNRN}, \textbf{FNRP}.

In the examples from both panels of Figure~\ref{fig_FT} we have $T(0)=T(\tilde{T})$, $F(0)=F(\tilde{T})$. Thus the discontinuity of $T(\Delta)$ at $\delta_{2}=z_{2}$ in Figure~\ref{fig_FT}(b) is due to how the cycle length map is defined and not due to the dynamics. This type of discontinuity is called Type 1~\citep{Winfree_1980}. On the other hand, the resetting time presents three discontinuities on each panel of Figure~\ref{fig_FT} and independently of the way that $F(\Delta)$ is defined, it has two discontinuities which are defined as Type 0~\citep{Winfree_1980}.

%
\section{Phase Response to a Periodic Stimulus}\label{sec.per}
In a clinical setting, it is more likely that a patient would receive periodic administrations of a cytokine to ameliorate the effects of periodic chemotherapy.  However, there has been controversy about how to best time these administrations.

Here we examine the response of the DDE~\eqref{dxdt} to a periodic perturbation $p(t)=p(t+T_{p})$. We keep $\tau>0$, $-\beta_{U}<0<\beta_{L}$, $a>0$ and $\sigma\in(0,\tau]$.
The perturbation is ON during a time interval $\sigma$ and OFF during a time interval $\alpha>0$, so the period is $T_{p} = \sigma+\alpha$. The perturbation is defined by
\begin{equation}\label{eq.pert.per.disc}
p(t) = \left\{\begin{array}{ll}
a, \qquad\mbox{if}\qquad t\in[\Delta_{n},\Delta_{n}+\sigma],\\
0,  \qquad\mbox{if}\qquad t\in(\Delta_{n}+\sigma,\Delta_{n+1}),
\end{array}\right.
\end{equation}

\noindent where $\Delta_{n+1}=\Delta_{n}+T_{p}$, $n\in\mathbb{N}$ and $\Delta_{0}\in[0,\tilde{T})$.

We denote the solution of the perturbed DDE by $x^{(p)}:\mathbb{R}\longrightarrow\mathbb{R}$ which, up to $t=\Delta_{0}$, is equal to the periodic solution $\tilde{x}$. For $t\geq\Delta_{0}$ the function $x^{(p)}$ is defined by
\begin{equation}\label{eq.diff.pert}
x^{\prime}(t) = -x(t)+f(x(t-\tau))+p(t),
\end{equation}
\noindent with $f(x(t-\tau))$ and $p(t)$ given by~\eqref{ftau} and~\eqref{eq.pert.per.disc}, respectively.  The solution of~\eqref{eq.diff.pert} is built up piecewise by functions of the form  $A_{k}+B_{k}\mathrm{e}^{-(t-\phi_{k})}$ for each interval $[\phi_{k},\phi_{k+1}]$ with $k\in\mathbb{N}$ and $\phi_{1}=\Delta_{0}$, where $\phi_{k}$ are the points where the derivative is discontinuous. These discontinuity points are known as \textit{breaking points} in the literature~\citep{Bellen_2003}. Along the solutions of~\eqref{eq.diff.pert} the breaking points are located at the points where $f(x(t-\tau))$ switches from a negative to non-negative value or $p(t)$ switches from positive to zero. Thus for $k,n\in\mathbb{N}$ the solution of~\eqref{eq.diff.pert}, $x^{(p)}(t) $, is given by

\begin{equation}\nonumber
\left\{\begin{array}{ll}
\!\!\beta_{L}+a+(x^{(p)}(\phi_{k})-\beta_{L}-a)\mathrm{e}^{-(t-\phi_{k})}\!, \ &\!\mbox{if}\!\ \  t\in[\Delta_{n},\Delta_{n}+\sigma]\ \ \:\!\mbox{and}\!\ \  x(t-\tau)<0,\\
\!\!\beta_{L}+(x^{(p)}(\phi_{k})-\beta_{L})\mathrm{e}^{-(t-\phi_{k})}\!,  \ &\!\mbox{if}\!\ \  t\in(\Delta_{n}+\sigma,\Delta_{n+1})\ \ \!\mbox{and}\!\ \  x(t-\tau)<0,\\
\!\!-\beta_{U}+a+(x^{(p)}(\phi_{k})+\beta_{U}-a)\mathrm{e}^{-(t-\phi_{k})}\!, \ &\!\mbox{if}\!\ \  t\in[\Delta_{n},\Delta_{n}+\sigma]\ \ \:\!\mbox{and}\!\ \  x(t-\tau)\geq 0,\\
\!\!-\beta_{U}+(x^{(p)}(\phi_{k})+\beta_{U})\mathrm{e}^{-(t-\phi_{k})}\!, \ &\!\mbox{if}\!\ \  t\in(\Delta_{n}+\sigma,\Delta_{n+1})\ \ \!\mbox{and}\!\ \  x(t-\tau)\geq 0.\\
\end{array}\right.
\end{equation}

%
%
\noindent In order to examine the effect of the periodic stimulus we will consider the perturbation period $T_{p}$ and resetting time $F(\Delta)$ due to each perturbation. Recall that by definition $\alpha>0$ and from Remark~\ref{rem.min.ress} we have $\ubar{F}=\sigma$ for all $\Delta\in[0,\tilde{T})$, then $T_{p}>\ubar{F}$. The analysis of the periodic perturbation must distinguish between two cases, $T_{p}>\bar{F}$ and $\ubar{F}<T_{p}\leq\bar{F}$.
\begin{figure}[t]
\centering
\includegraphics[width=0.73\textwidth]{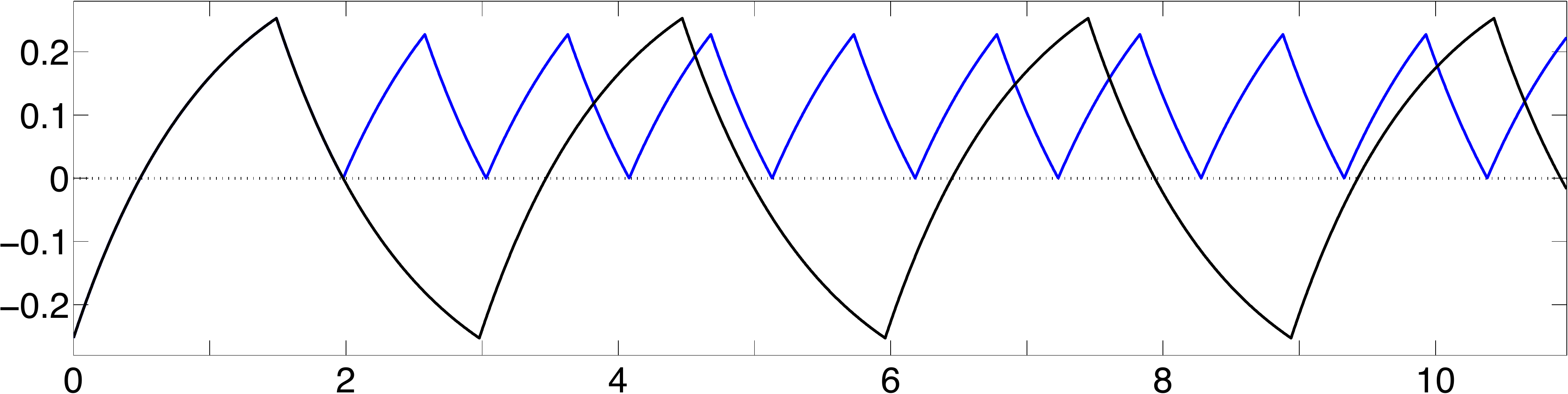}
\put(-334,81){\footnotesize{$x$}}
\put(-6,1.5){\footnotesize{$t$}}
\vspace*{-02mm}
\caption{Example of Proposition~\ref{prop.LC1} for $\beta_{U}=0.4$, $\beta_{L}=0.4$, $\sigma=0.6$, $\tau=1$, $\alpha=0.45$, $\Delta_{0}=z_{2}$ and $a=a_{1}=0.90384$, where $a_{1}$ is given by~\eqref{eq.a1}. The numerical solution with periodic perturbations ($T_{p} = 1.05$) is represented by the blue line and the solution without perturbations by the black line ($\tilde{T}=2.97976$).}
\label{fig_PropLC2}
\end{figure}

\vspace{2mm}
\noindent\textbf{Case (i)} $T_{p}>\bar{F}$\textbf{.} After each perturbation the solution $x^{(p)}$ returns to the periodic solution $\tilde{x}$, with a new phase, before the next perturbation starts. Hence, in this case the solution $x^{(p)}(t)$ can be computed on each interval $\Delta_{n}\leq t\leq\Delta_{n}+\sigma$ and $\Delta_{n}+\sigma\leq t\leq\Delta_{n+1}$, with $n\in\mathbb{N}$, as was done in~\citet[Section 5]{Mackey_2017}. For the special case where $T_{p}=F(\Delta_{0})$ and $x^{(p)}(\Delta_{0})=x^{(p)}(\Delta_{1})$ the solution $x^{(p)}$ is periodic. For this singular case the time interval that the solution takes to return to the limit cycle $\tilde{x}$ is equal to the perturbation period and is such that the solution  returns to $\tilde{x}$ at the same phase that it had when it  left the limit cycle. For the cases where $T_{p}\neq T(\Delta_{0})$ or $T_{p}=F(\Delta_{0})$ and $x^{(p)}(\Delta_{0})\neq x^{(p)}(\Delta_{1})$ the solution may become periodic after a finite number of stimuli or may continue to be non-periodic.

\vspace{2mm}
\noindent\textbf{Case (ii)} $\ubar{F}<T_{p}\leq\bar{F}$\textbf{.} For this case the solution will be periodic if $(T_{p}/\tilde{T})\in\mathbb{N}_{>0}$. Again for the special case where $T_{p}=F(\Delta_{0})$ and $x^{(p)}(\Delta_{0})=x^{(p)}(\Delta_{1})$ the solution $x^{(p)}$ is periodic, but here $T(\Delta_{0})\leq\bar{F}$, as is shown in the example of Figure~\ref{fig_PropLC2}. For the cases where $T_{p}\neq T(\Delta_{0})$ or $T_{p}=F(\Delta_{0})$ and $x^{(p)}(\Delta_{0})\neq x^{(p)}(\Delta_{1})$, after the end of the first perturbation, $t>(\Delta_{0}+\sigma)$, the next perturbation may start before or after the perturbed solution returns to the limit cycle. In both cases the second stimulus will not start at the same phase of the limit cycle where the first stimulus started. Thus, each perturbation starts with a different phase, with respect to the previous, and it may occur that the phases repeat after a finite number of stimuli, resulting in a periodic solution. Otherwise, the phases will not repeat during successive perturbations and the solution may be quasi-periodic or non-periodic.

The numerical solutions of Figure~\ref{fig_PropLC2} and onwards were computed using the MATLAB \texttt{dde23} routine~\citep{Matlab}. All solutions of~\eqref{eq.diff.pert} are composed of piece-wise function segments on intervals $[\phi_{k},\phi_{k+1}]$, where $\phi_{k}$ are the breaking points. These breaking points were detected and included in the solution meshes by using the MATLAB \texttt{events} function~\citep{Matlab,Shampine_2003}.

\par In Proposition~\ref{prop.LC1} below we show that if $a\geq a_{1}$ and $\Delta_{0}=z_{2}$, then the solution converges to the orbit given by~\eqref{eq.xpt1} and~\eqref{eq.xpt2}.
\begin{proposition}\label{prop.LC1}
If $\Delta_{0}=z_{2}$ and $a\geq a_{1}$, where
\begin{equation}\label{eq.a1}
a_{1} = \beta_{U}\frac{(\mathrm{e}^{\alpha}-\mathrm{e}^{-\sigma})}{(1-\mathrm{e}^{-\sigma})},
\end{equation}
\noindent then the solution of~\eqref{eq.diff.pert} is given by
\begin{equation}\label{eq.xpt1}
x^{(p)}(t) = -\beta_{U}+a+(x^{(p)}(\Delta_{n})+\beta_{U}-a)\mathrm{e}^{-(t-\Delta_{n})}, \qquad\mbox{if}\qquad t\in[\Delta_{n},\Delta_{n}+\sigma],
\end{equation}
\begin{equation}\label{eq.xpt2}
x^{(p)}(t) = -\beta_{U}+(x^{(p)}(\Delta_{n}+\sigma)+\beta_{U})\mathrm{e}^{-(t-\Delta_{n}-\sigma)}, \qquad\mbox{if}\qquad t\in[\Delta_{n}+\sigma,\Delta_{n+1}],
\end{equation}
\noindent where
\begin{equation}\label{eq.xpd}
x^{(p)}(\Delta_{n}) = \beta_{U}(\mathrm{e}^{-nT_{p}}-1)+a(1-\mathrm{e}^{-\sigma})\mathrm{e}^{-\alpha}\sum_{k=0}^{n}\mathrm{e}^{-kT_{p}}, \quad\mbox{for}\quad n\in\mathbb{N}_{>0},
\end{equation}
\begin{equation}\label{eq.xpdsigma}
x^{(p)}(\Delta_{n}+\sigma) = \beta_{U}(\mathrm{e}^{-nT_{p}-\sigma}-1)+a(1-\mathrm{e}^{-\sigma})\sum_{k=0}^{n}\mathrm{e}^{-kT_{p}}, \quad\mbox{for}\quad n\in\mathbb{N}.
\end{equation}
\end{proposition}

\par In Proposition~\ref{prop.LC2} we show that the solution of Proposition~\ref{prop.LC1} converges to a limit cycle. We compute this periodic solution and show that if the phase of the first perturbation $\Delta_{0}$ overlaps with the minimum point of the limit cycle, then the perturbed solution settles down on this limit cycle.
\begin{proposition}\label{prop.LC2}
If $a\geq a_{1}$, then: (i) the solution of~\eqref{eq.diff.pert} converges to the limit cycle given by~\eqref{eq.xp}, (ii) for $\Delta_{0}=0$ and a non-negative history function $\varphi(t)\geq 0$, with $t\in[-\tau,0)$, and such that $\varphi(0)=\ubar{x}^{(p)}$, the solution of~\eqref{eq.diff.pert} settles down on the limit cycle given by
\begin{equation}\label{eq.xp}
x^{(p)}(t) = \left\{\begin{array}{ll}
-\beta_{U}+a+(\ubar{x}^{(p)}+\beta_{U}-a)\mathrm{e}^{-(t-\Delta_{n})}, &\quad\mbox{if}\quad t\in[\Delta_{n},\Delta_{n}+\sigma],\\[1mm]
-\beta_{U}+(\bar{x}^{(p)}+\beta_{U})\mathrm{e}^{-(t-\Delta_{n}-\sigma)}, &\quad\mbox{if}\quad t\in[\Delta_{n}+\sigma,\Delta_{n+1}],
\end{array}\right.
\end{equation}

\noindent for $n\in\mathbb{N}$, where $\ubar{x}^{(p)}$ and $\bar{x}^{(p)}$ are defined by
\begin{equation}\label{eq.xpn1}
\ubar{x}^{(p)}  = -\beta_{U}+a\frac{(1-\mathrm{e}^{-\sigma})}{(\mathrm{e}^{\alpha}-\mathrm{e}^{-\sigma})},
\end{equation}
\begin{equation}\label{eq.xpn2}
\bar{x}^{(p)} = -\beta_{U}+a\frac{(1-\mathrm{e}^{-\sigma})}{(\mathrm{e}^{\alpha}-\mathrm{e}^{-\sigma})}\mathrm{e}^{\alpha}.
\end{equation}
\end{proposition}

\par In Propositions~\ref{prop.LC1} and~\ref{prop.LC2} it is shown that the perturbed solution converges to a limit cycle if $a\geq a_{1}$ and $\Delta_{0}=z_{2}$. For $a=a_{1}$ we have $\ubar{x}^{(p)}=0$ and $\bar{x}^{(p)}=\beta_{U}(\mathrm{e}^{\alpha}-1)$.
\begin{figure}[!htbp]
\centering
\includegraphics[width=0.49\textwidth,height=0.278\textwidth]{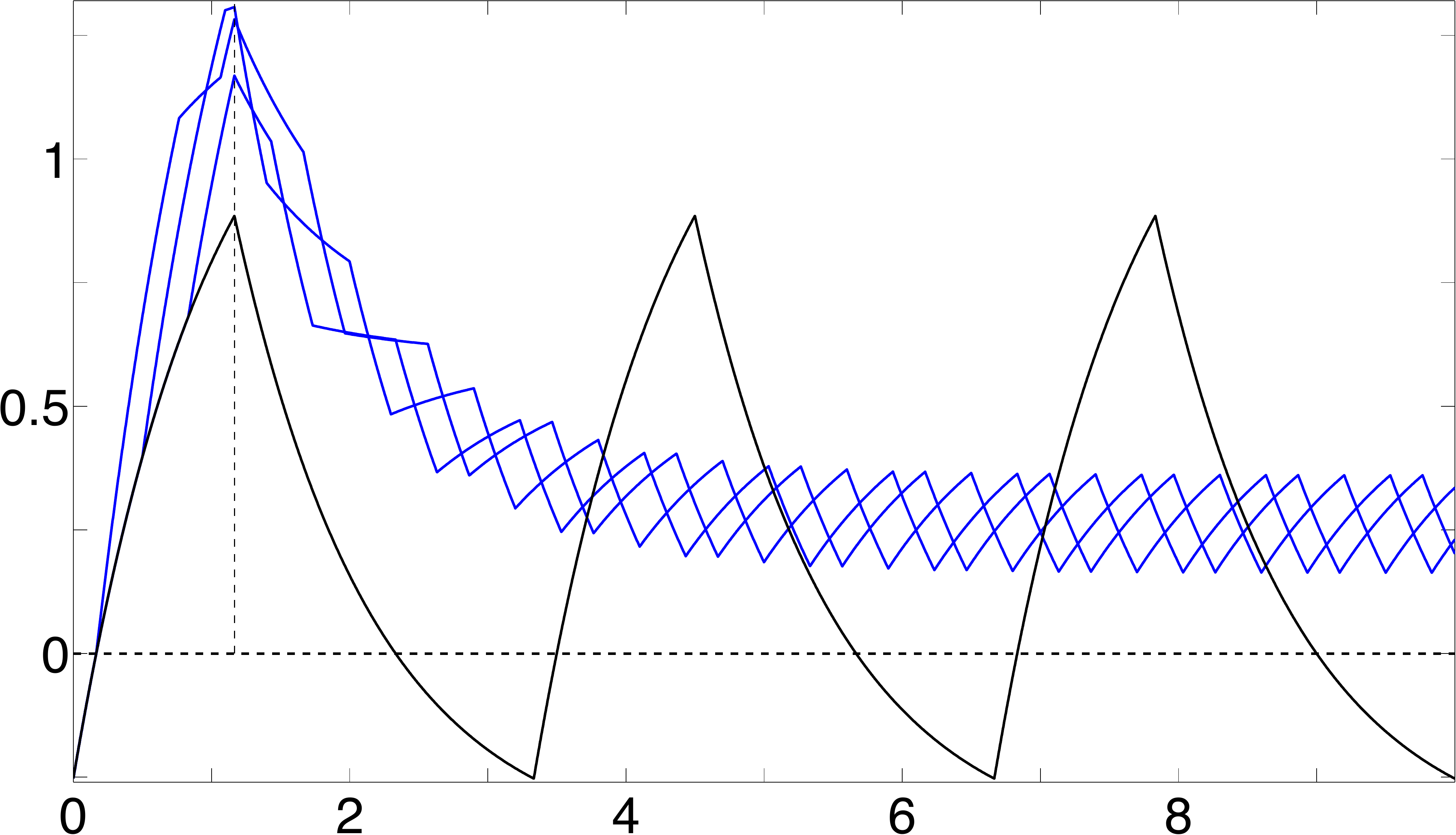}
\put(-15,120){\footnotesize{\textit{(a)}}}
\put(-15,0.5){\footnotesize{\textit{t}}}
\put(-225,120){\footnotesize{\textit{x}}}
\put(-201,21){\footnotesize{$t_{max}$}}
\hspace*{0mm}
\includegraphics[width=0.49\textwidth,height=0.278\textwidth]{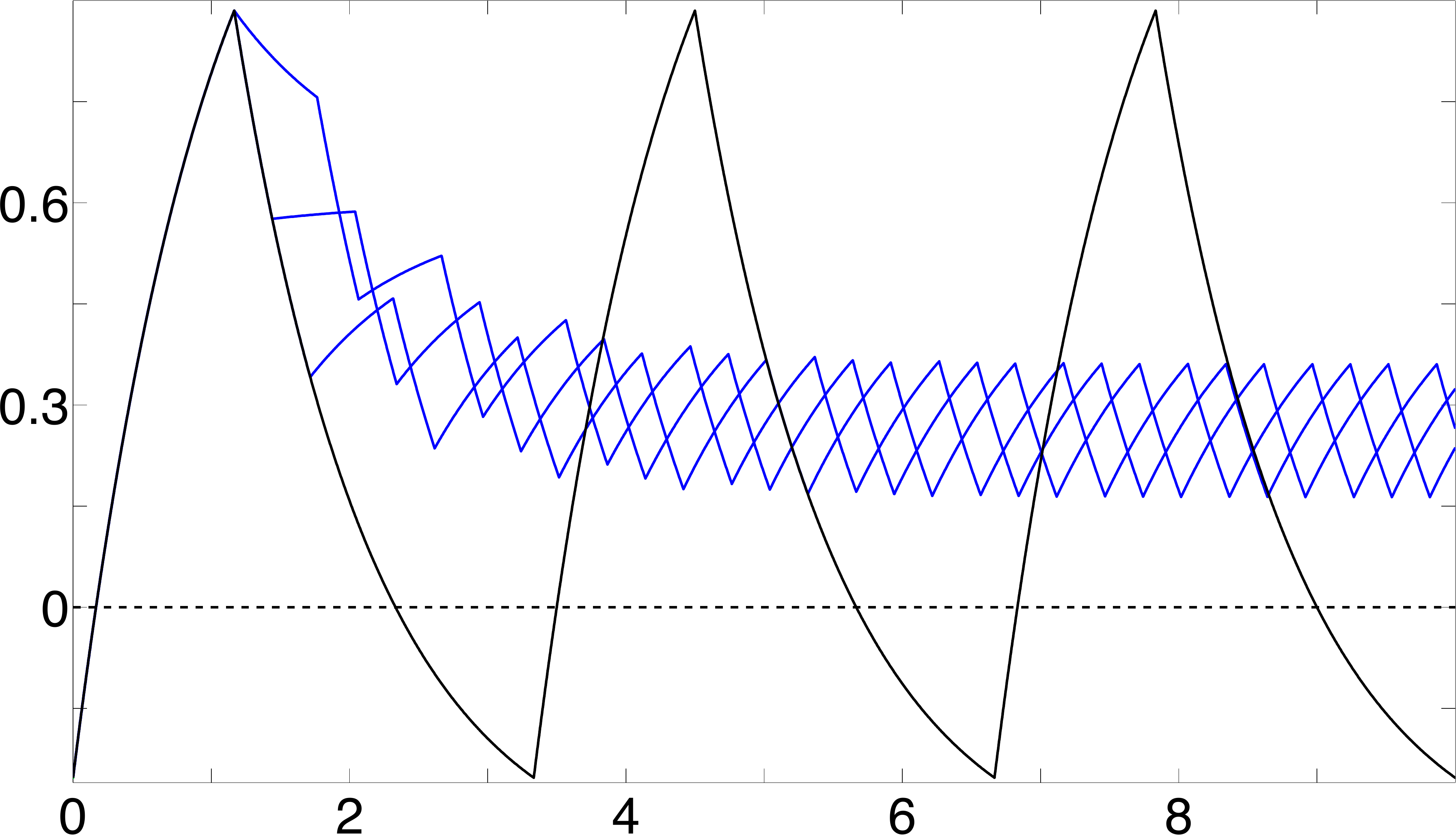}
\put(-15,120){\footnotesize{\textit{(b)}}}
\put(-15,0.5){\footnotesize{\textit{t}}}
\put(-225,120){\footnotesize{\textit{x}}}\\
\vspace*{1mm}
\includegraphics[width=0.49\textwidth,height=0.278\textwidth]{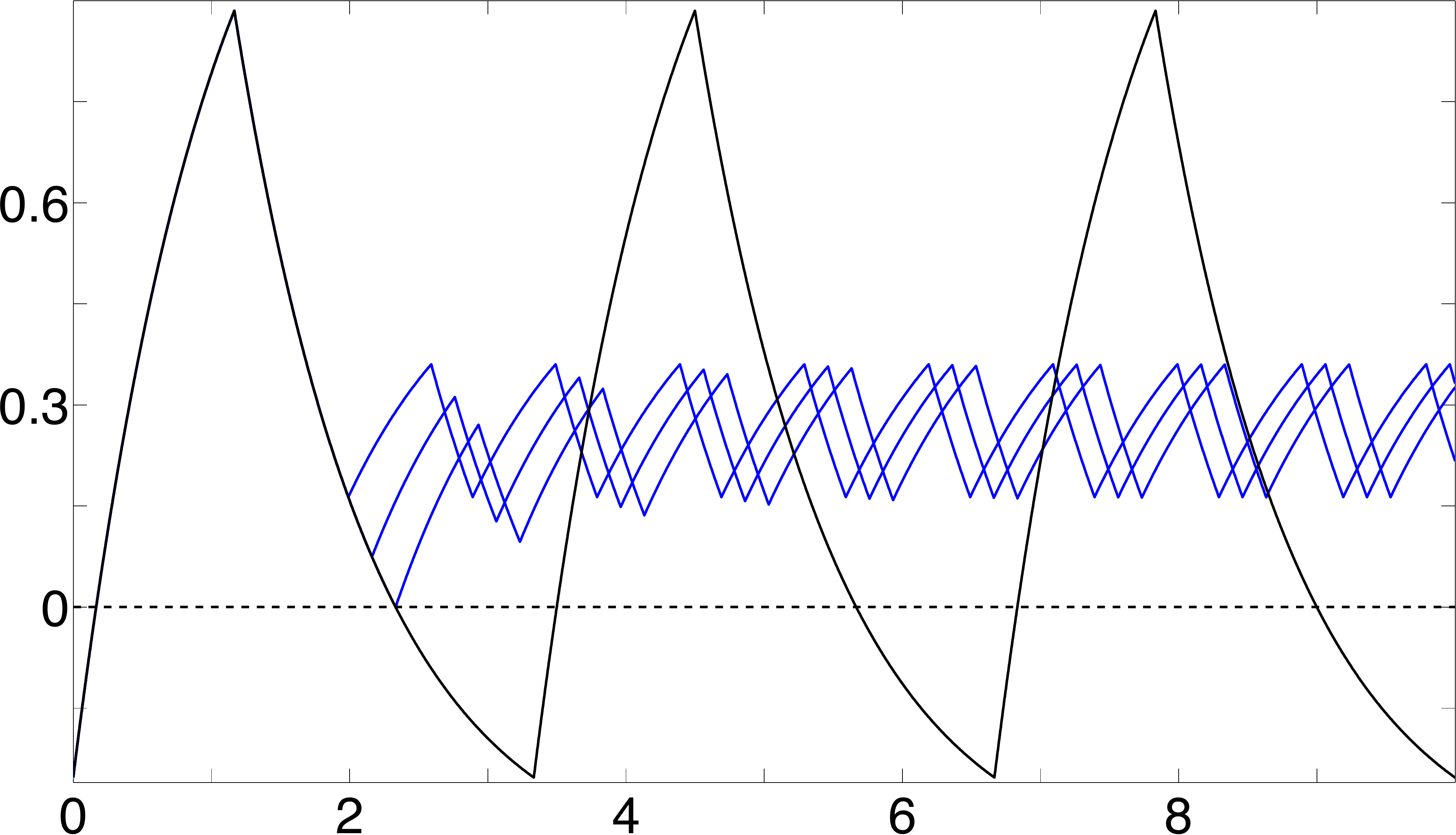}
\put(-15,120){\footnotesize{\textit{(c)}}}
\put(-15,0.5){\footnotesize{\textit{t}}}
\put(-225,120){\footnotesize{\textit{x}}}
\hspace*{0mm}
\includegraphics[width=0.49\textwidth,height=0.278\textwidth]{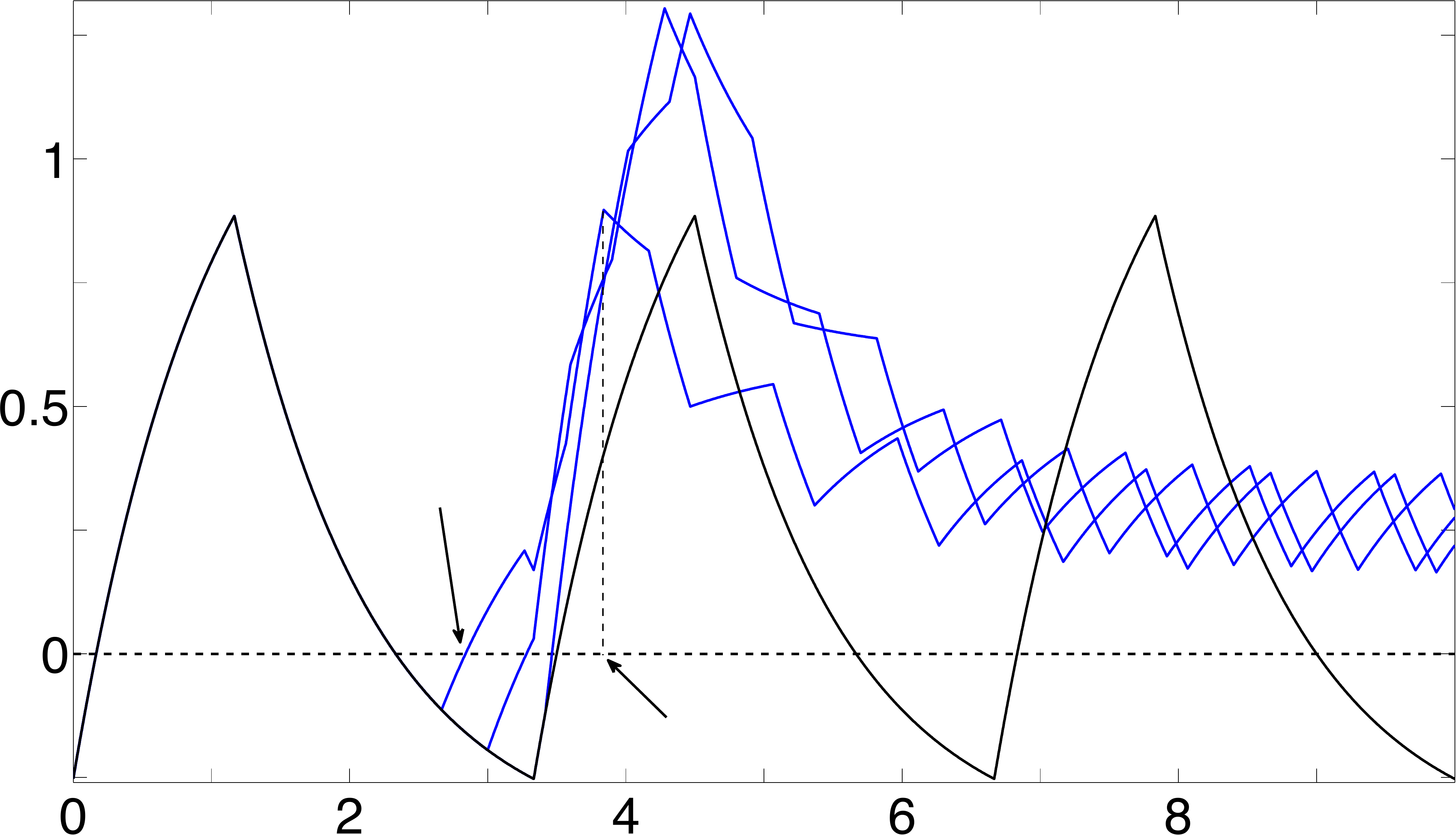}
\put(-15,120){\footnotesize{\textit{(d)}}}
\put(-15,0.5){\footnotesize{\textit{t}}}
\put(-225,120){\footnotesize{\textit{x}}}
\put(-123,17){\footnotesize{$z_{p,1}+\tau$}}
\put(-167,55){\footnotesize{$z_{p,1}$}}
\vspace*{-02mm}
\caption{Proposition~\ref{prop.LC} cases with $\beta_{U}=0.4$, $\beta_{L}=1.4$, $\sigma=0.6$, $a=1$, $\tau=1$, $\alpha=0.3$, and with $\Delta_{0}\in[0,\tilde{T})$. The numerical solutions with periodic perturbations are represented by the blue line and the solution without perturbation by the black line. For each panel the $\Delta_{0}$ interval is: (a) $(\tilde{z}_{1},t_{max})$, (b) $[t_{max},\Delta_{l})$, (c) $[\Delta_{l},\tilde{z}_{2})$ with $\Delta_{l}$ given by~\eqref{DeltaL}, (d) $[\tilde{z}_{2},\tilde{T}+\tilde{z}_{1}]$.}
\label{fig_PropLC4x4}
\end{figure}
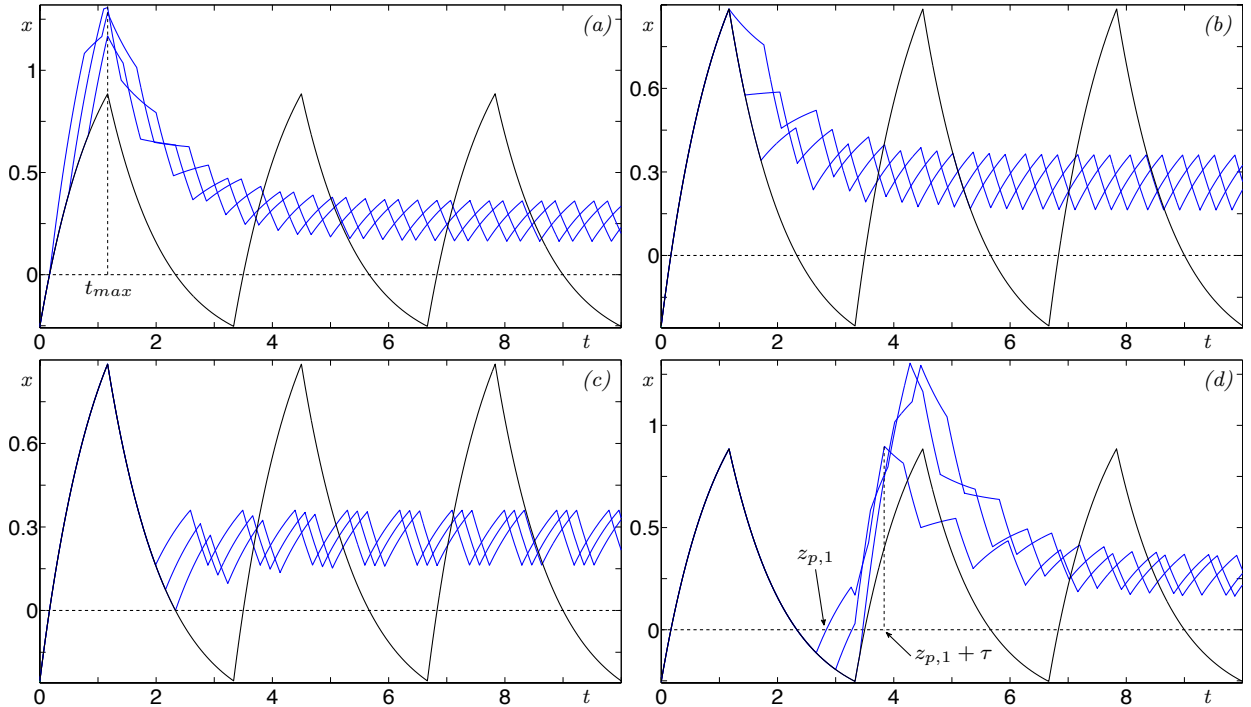

\par For $T_{p}\neq T(\Delta_{0})$ or $T_{p}=F(\Delta_{0})$ and $x^{(p)}(\Delta_{0})\neq x^{(p)}(\Delta_{1})$ the long-time behavior of the solutions for both \textbf{Cases (i)-(ii)} described earlier does not depend on the value of $\Delta_{0}\in[0,\tilde{T})$. Indeed, the results of Propositions~\ref{prop.LC1} can be extended for all $\Delta_{0}\in[0,\tilde{T})$, which was done in Proposition~\ref{prop.LC}.
Figure~\ref{fig_PropLC4x4}  shows examples of solutions, distinguishing four cases, which satisfy the conditions of Proposition~\ref{prop.LC}.  All solutions of Figure~\ref{fig_PropLC4x4} converge to a limit cycle given by~\eqref{eq.xp}.  In panel (b) and (c) the solution points $x^{(p)}(\Delta_{n})$ and $x^{(p)}(\Delta_{n}+\sigma)$ exponentially converge to~\eqref{eq.xpn1} and~\eqref{eq.xpn2}, respectively.  In panel (a) this convergence occurs for $t\geq t_{max}$ while for panel (d) it occurs for $t>z_{p,1}+\tau$, where $z_{p,1}$ is the first zero of $x^{(p)}(t)$ with $t>\Delta_{0}$.

\begin{proposition}\label{prop.LC}
For any initial phase $\Delta_{0}\in[0,\tilde{T})$ the Proposition~\ref{prop.LC1} holds and the solution converges to the limit cycle given by Proposition~\ref{prop.LC2}.
\end{proposition}

\section{Treatment Implications}\label{sec:treatment}

At this point it is interesting to consider our results with this extremely simple model  in the context of a hypothetical patient with cyclic circulating blood cell numbers (e.g. cyclic neutropenia, or cycling induced by chemotherapy) being treated with periodic G-CSF administration. We denote the normal level of neutrophils by $x_{norm}$. For cyclic neutropenia the circulating neutrophil numbers typically oscillate from normal levels to very low levels with a period of about 19 to 21 days~\citep{Colijn_2005b}.  The period for cycling induced by periodic chemotherapy is approximately equal to the period of the chemotherapy~\citep{Craig_2016}.

One issue of interest is whether or not it is possible to abrogate the severe neutropenic phases of the oscillation in the model as is done in practice by keeping the circulating neutrophil levels equal to or greater than normal. We can answer this in the affirmative, since the condition $\ubar{x}^{(p)}\geq x_{norm}$ is sufficient to end the neutropenia.
Using~\eqref{eq.a1} and~\eqref{eq.xpn1}, the condition to end the neutropenia can be written as
\begin{equation}\label{eq.xp3}
\ubar{x}^{(p)} =\beta_{U}\left(\frac{a}{a_{1}}-1\right)\geq x_{norm}.
\end{equation}
From~\eqref{eq.xp3} we obtain the condition $a\geq a_{1}(1+x_{norm}/\beta_{U})$, which satisfies the condition $a\geq a_{1}$ from Proposition~\ref{prop.LC1} and increases the nadir of the oscillations to, or above, the normal level $x_{norm}$.

During the periodic G-CSF administration the hypothetical patient will have a neutrophil oscillation described by a limit cycle $x^{(p)}$ given by Proposition~\ref{prop.LC2} with an oscillation amplitude given by
\begin{equation}\nonumber
(\bar{x}^{(p)}-\ubar{x}^{(p)}) = a\dfrac{(1-\mathrm{e}^{-\sigma})}{(\mathrm{e}^{\alpha}-\mathrm{e}^{-\sigma})}(\mathrm{e}^{\alpha}-1)>0.
\end{equation}
Notice that $\bar{x}^{(p)}\geq\ubar{x}^{(p)}\geq 0$ and we only have $\bar{x}^{(p)}=\ubar{x}^{(p)}$ for $\alpha=0$, or $\sigma=0$, or $a=0$.

\begin{figure}[!tb]
\includegraphics[width=0.49\textwidth,height=0.276\textwidth]{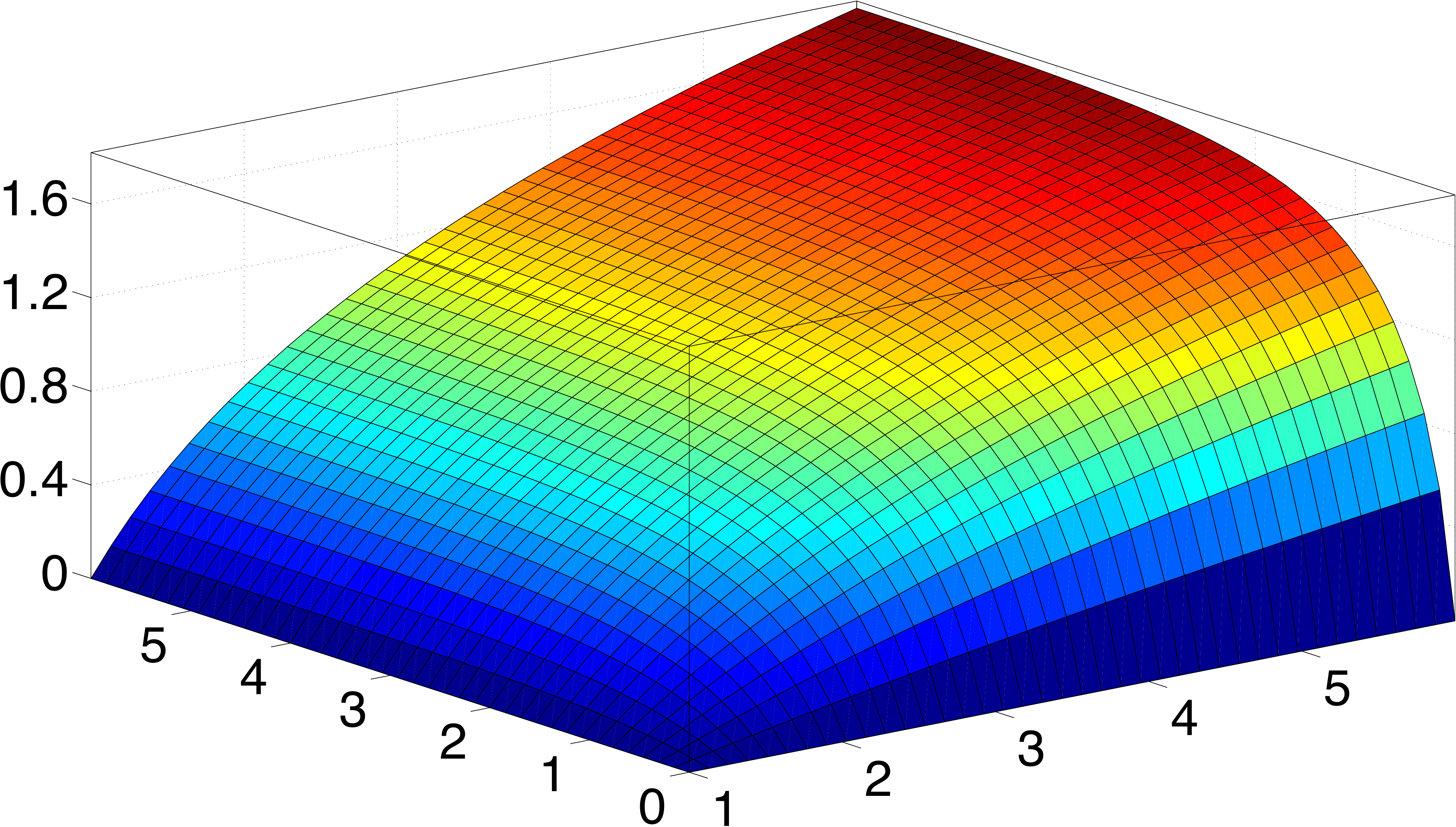}
\put(-229,2.5){\footnotesize{\textit{(a)}}}
\put(-54,5){\footnotesize{$a$}}
\put(-189,10){\footnotesize{$\sigma$}}
\put(-226,106){\footnotesize{$\alpha$}}
\hspace*{1mm}
\includegraphics[width=0.49\textwidth,height=0.276\textwidth]{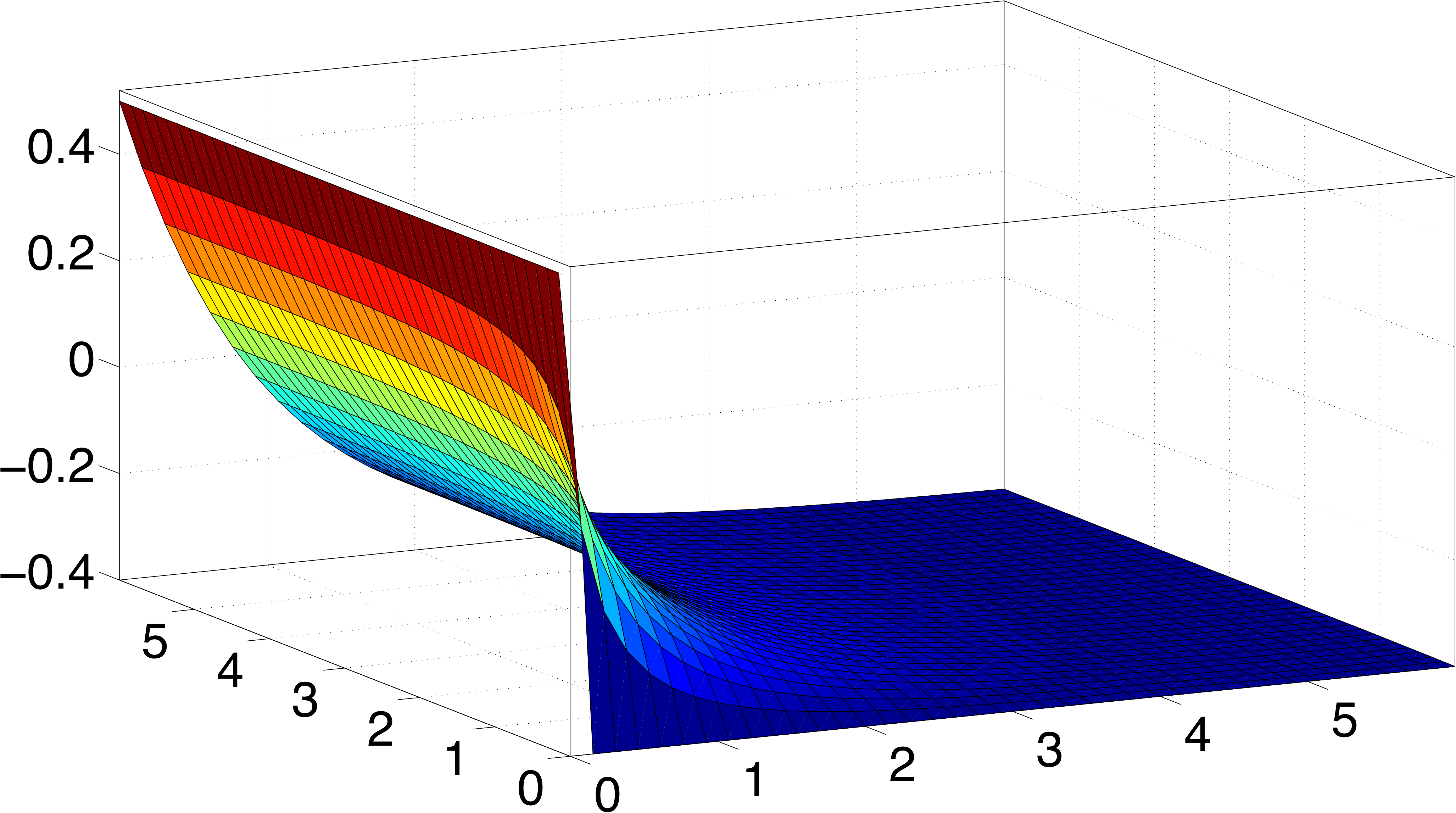}
\put(-229,2.5){\footnotesize{\textit{(b)}}}
\put(-53,2.5){\footnotesize{$\alpha$}}
\put(-198,15){\footnotesize{$\sigma$}}
\put(-227,116){\footnotesize{$\protect\ubar{x}^{(p)}$}}
\vspace*{-02mm}
\caption{For $\beta_{L}=2$, $\beta_{U}=0.4$ and $x_{norm}=0.6$: (a) the time interval between cytokine administrations $\alpha$, given by~\eqref{eq.alpha}, increases as the duration $\sigma$ and/or amplitude $a$ of the cytokine dose increase. (b) with $a=0.9$ the minimum $\protect\ubar{x}^{(p)}$, given by~\eqref{eq.xpn1}, is decreasing with respect to $\alpha$ and increasing with respect to $\sigma$.}
\label{fig_neutrop}
\end{figure}
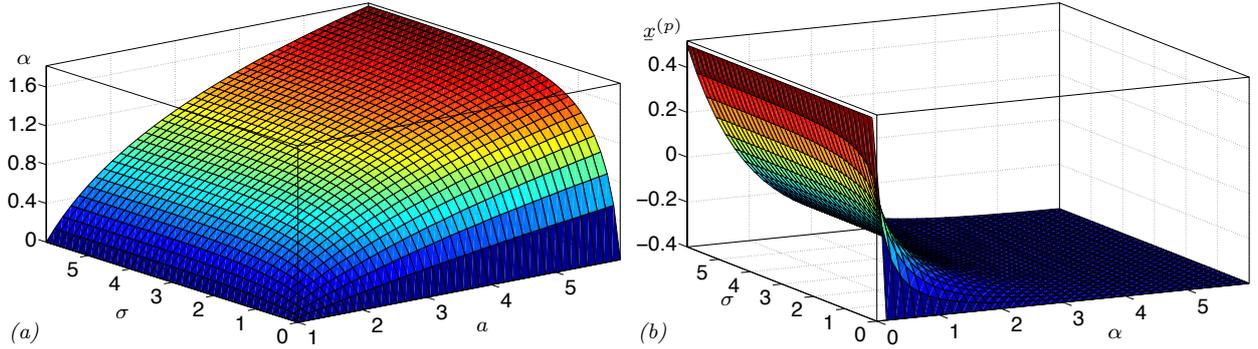

To avoid neutropenia, in the limiting case we can force $\ubar{x}^{(p)}=x_{norm}$ in~\eqref{eq.xp3}, which gives $a=a_{1}(1+x_{norm}/\beta_{U})$. This satisfies the condition $a\geq a_{1}$ from Proposition~\ref{prop.LC1} and in combination with~\eqref{eq.a1} yield the minimal interval between administrations $\alpha$ to avoid neutropenia as function of the duration $\sigma$ and the amplitude $a$ of cytokine administration:
\begin{equation}\label{eq.alpha}
\alpha = \ln\left(\frac{a+(x_{norm}+\beta_{U}-a)\mathrm{e}^{-\sigma}}{x_{norm}+\beta_{U}}\right).
\end{equation}
In Figure~\ref{fig_neutrop}(a) we show that for $a>x_{norm}+\beta_{U}$ the minimal interval between cytokine administrations $\alpha$ to avoid neutropenia, given by~\eqref{eq.alpha},  increases as the duration $\sigma$ and/or the amplitude $a$ of the cytokine dose increase. In panel (b) we show that the minimum $\ubar{x}^{(p)}$, given by~\eqref{eq.xpn1}, is decreasing with respect to the time interval between  cytokine administration $\alpha$ and increasing with respect to the duration $\sigma$ of the administration.  Both effects are what one would intuitively expect.

\begin{figure}[!tb]
\centering
\includegraphics[width=0.73\textwidth]{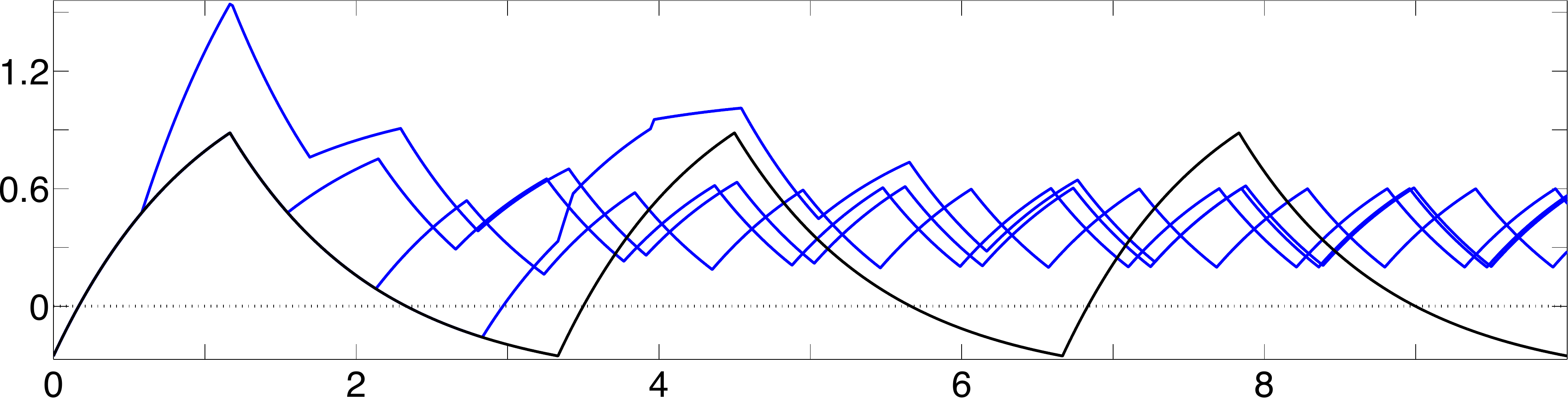}
\put(-337,81){\footnotesize{$x$}}
\put(-6,1.5){\footnotesize{$t$}}
\vspace*{-02mm}
\caption{Examples of orbits with $a\geq a_{1}$ and $\beta_{U}=0.4$, $\beta_{L}=1.4$, $\sigma=0.6$, $\tau=1$, $a=1.48655$, $\alpha=0.510826$ for four different values of $\Delta_{0}$:  $(\tilde{z}_{1}+t_{max})/2$, $(t_{max}+\Delta_{l})/2$ with $\Delta_{l}$ given by~\eqref{DeltaL}, $(\Delta_{l}+\tilde{z}_{2})/2$, $(\tilde{z}_{2}+\tilde{T}+\tilde{z}_{1})/2$. The numerical solutions with periodic perturbations are represented by the blue line and the solution without perturbations by the black line.}
\label{fig_PropLC9}
\end{figure}
We also investigate the hypothetical situation where the periodic G-CSF administration still results in neutropenia and the neutrophils oscillate between $f_{min}x_{norm}$ and $f_{max}x_{norm}$ where $f_{min} \leq 1 \leq f_{max}$. We mimic this situation by imposing the conditions
\begin{equation}\label{fminfmax}
\left\{\begin{array}{ll}
\ubar{x}^{(p)}-f_{min}x_{norm}=0,\\[1mm]
\bar{x}^{(p)}-f_{max}x_{norm}=0,
\end{array}\right.
\end{equation}
where $\ubar{x}^{(p)}$ and $\bar{x}^{(p)}$ are respectively given by~\eqref{eq.xpn1} and~\eqref{eq.xpn2}.

For a healthy adult human the normal circulating neutrophil level fluctuates around $0.22-0.85\times{10}^{9}\,\mathrm{cells/kg}$ of body mass~\citep{Craig_2016}. Assuming that we want to administer G-CSF to maintain the oscillation within these normal bounds and taking  $x_{norm}\approx 0.4$, then we have $f_{min}\approx 0.5$ and $f_{max}\approx 1.5$. In order to solve~\eqref{fminfmax} we need to fix one of the triplet $(a,\sigma, \alpha)$ and then solve for the remaining two parameters.

Both limit cycle extrema $\protect\ubar{x}^{(p)}$ and $\bar{x}^{(p)}$ are nonlinear increasing functions with respect to $\sigma$, nonlinear decreasing functions with respect to $\alpha$, and linearly increasing with respect to $a$. Thus we take $a$ and $\alpha$ as unknowns and solve~\eqref{fminfmax} with the set of parameters: $\beta_{U}=0.4$, $\beta_{L}=1.4$, $\sigma=0.6$ and $\tau=1$,  $x_{norm}=0.4$, $f_{min}=0.5$ and $f_{max}=1.5$. Using the MATLAB \texttt{fsolve} routine~\citep{Matlab} to solve~\eqref{fminfmax} gives $a=1.48655$ and $\alpha=0.510826$.

\begin{figure}[!t]
\includegraphics[width=0.49\textwidth,height=0.276\textwidth]{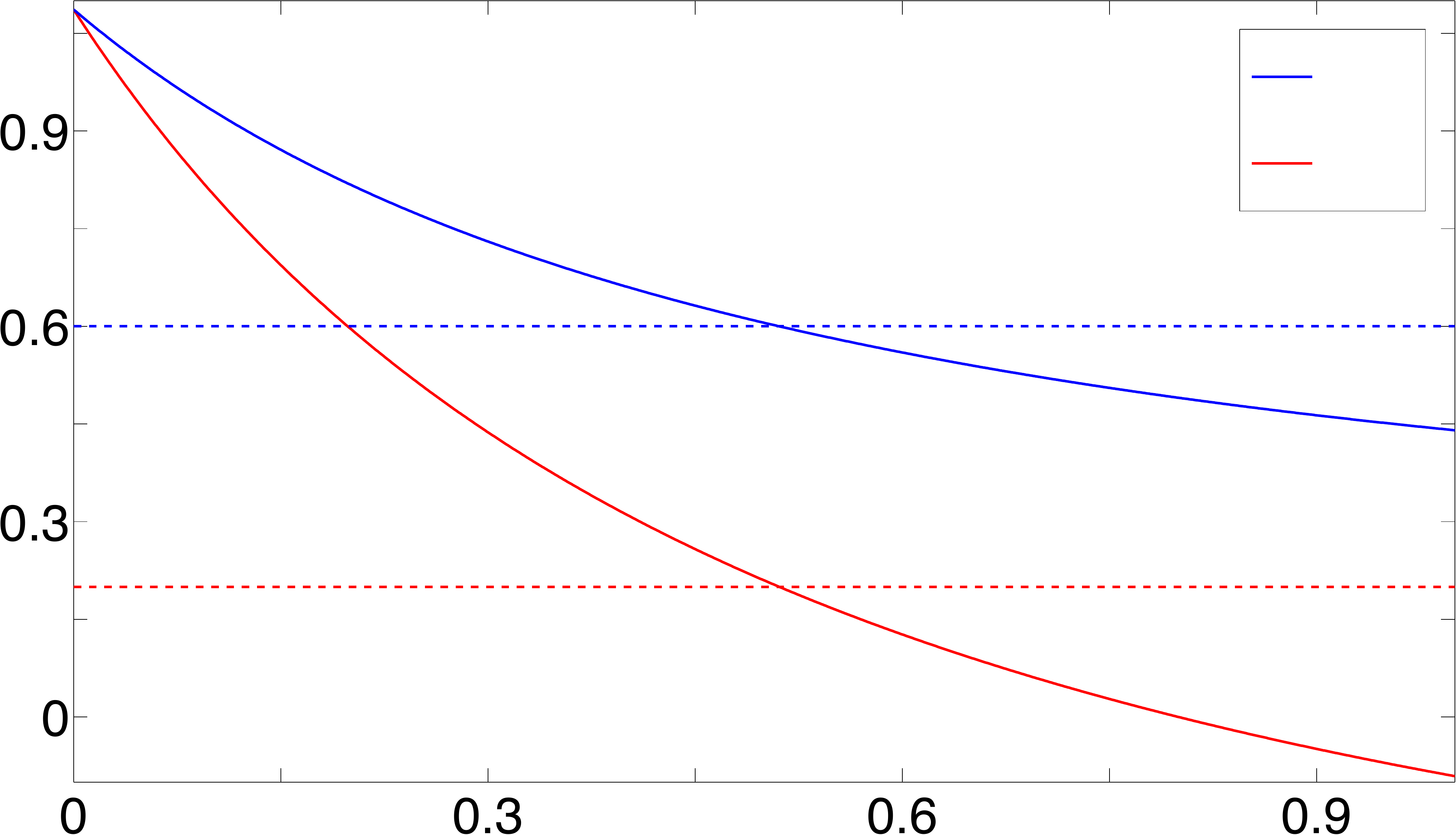}
\put(-115,120){\footnotesize{\textit{(a)}}}
\put(-9,1){\footnotesize{$\alpha$}}
\put(-20.5,115.6){\footnotesize{$\bar{x}^{(p)}$}}
\put(-20.5,102.2){\footnotesize{$\ubar{x}^{(p)}$}}
\hspace*{1mm}
\includegraphics[width=0.49\textwidth,height=0.276\textwidth]{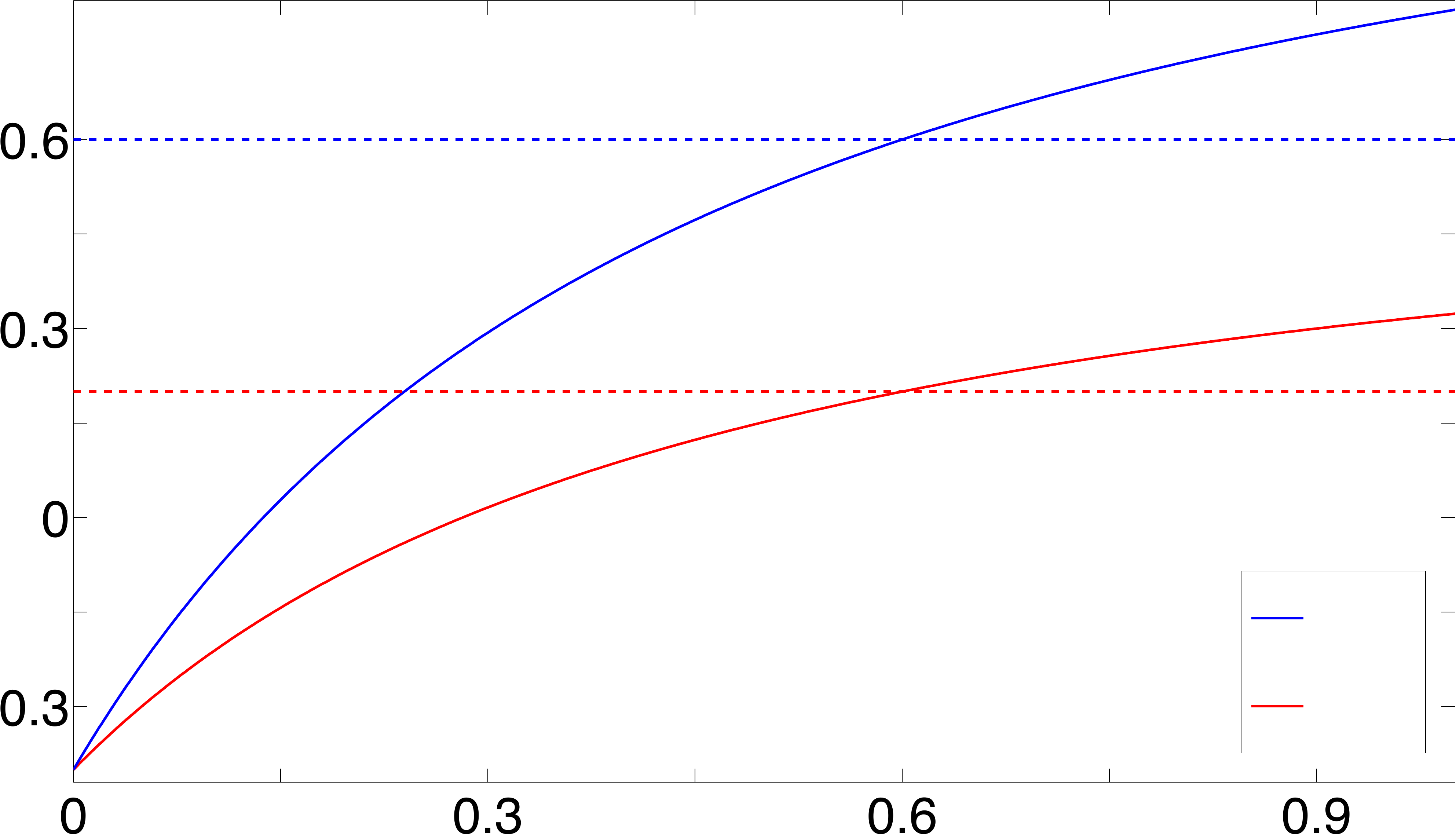}
\put(-115,120){\footnotesize{\textit{(b)}}}
\put(-9,1){\footnotesize{$\sigma$}}
\put(-20.5,31.8){\footnotesize{$\bar{x}^{(p)}$}}
\put(-20.5,18.4){\footnotesize{$\ubar{x}^{(p)}$}}\\

\vspace*{-3.5mm}
\includegraphics[width=0.49\textwidth,height=0.276\textwidth]{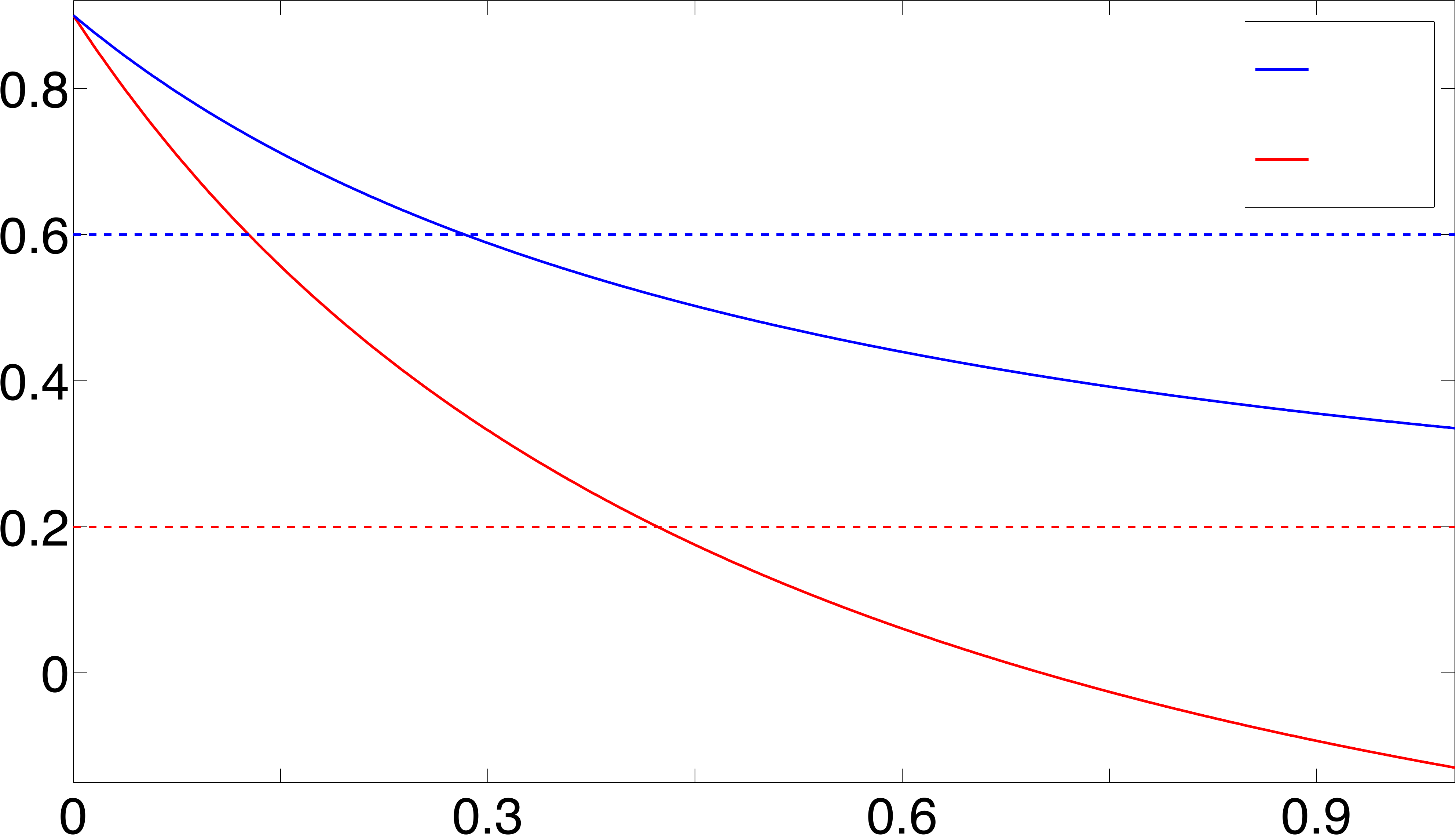}
\put(-115,120){\footnotesize{\textit{(c)}}}
\put(-9,1){\footnotesize{$\alpha$}}
\put(-20.5,115.6){\footnotesize{$\bar{x}^{(p)}$}}
\put(-20.5,102.2){\footnotesize{$\ubar{x}^{(p)}$}}
\hspace*{1mm}
\includegraphics[width=0.49\textwidth,height=0.276\textwidth]{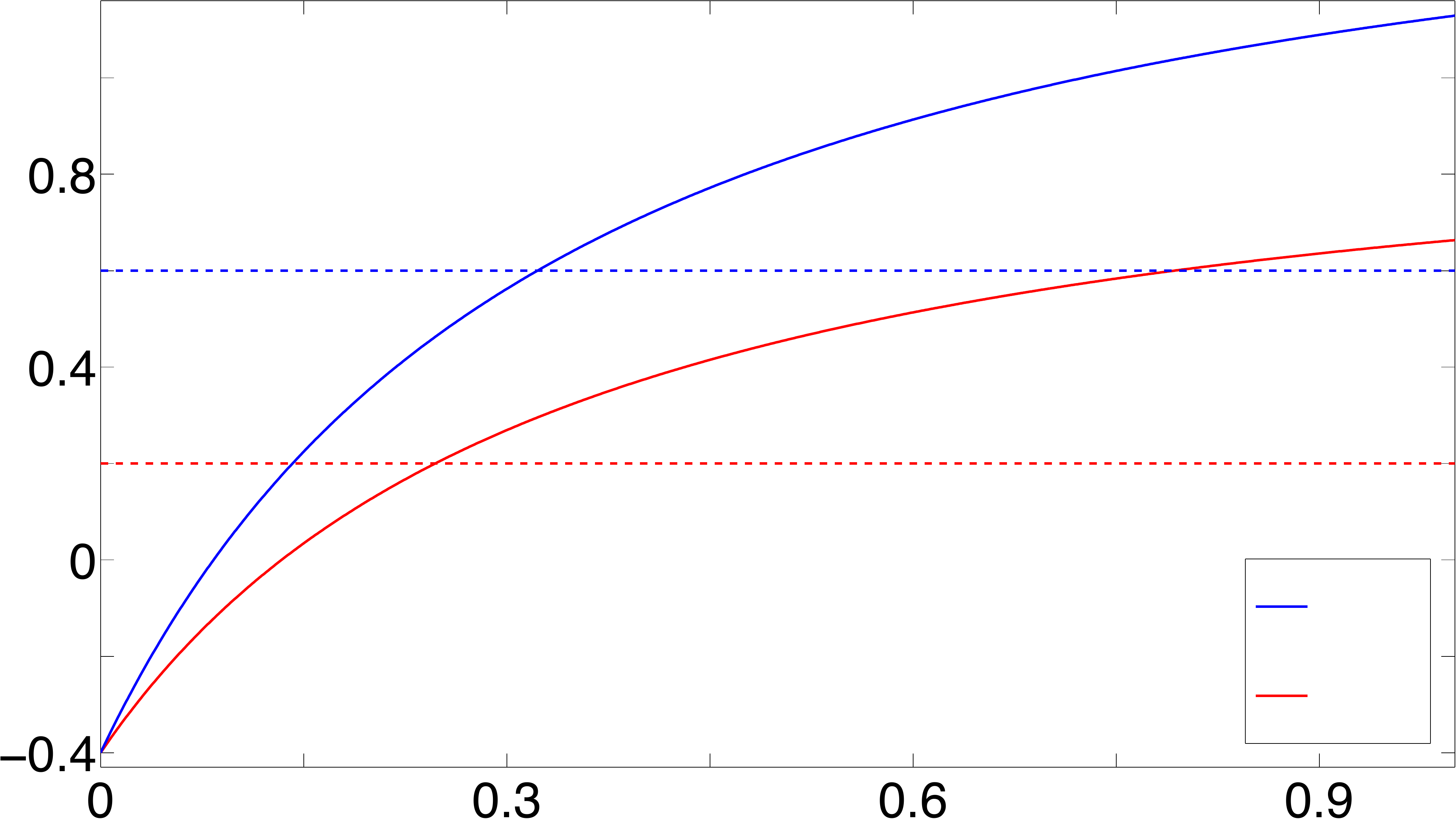}
\put(-115,120){\footnotesize{\textit{(d)}}}
\put(-9,1){\footnotesize{$\sigma$}}
\put(-20.5,31.8){\footnotesize{$\bar{x}^{(p)}$}}
\put(-20.5,18.4){\footnotesize{$\ubar{x}^{(p)}$}}
\vspace*{-02mm}
\caption{The extrema $\protect\ubar{x}^{(p)}$ and $\bar{x}^{(p)}$ are plotted as function of $\alpha$ in panel (a) and (c) and as function of $\sigma$ in panel (b) and (d) using the same parameters of Figure~\ref{fig_PropLC9} in panel (a) and (b), but with $a=1.3$ in panel (c) and $a=1.8$ in panel (d). The blue dashed lines correspond to  $x=f_{max}x_{norm}$ and the red dashed lines correspond to $x=f_{min}x_{norm}$.
}%
\label{fig_osc}
\end{figure}
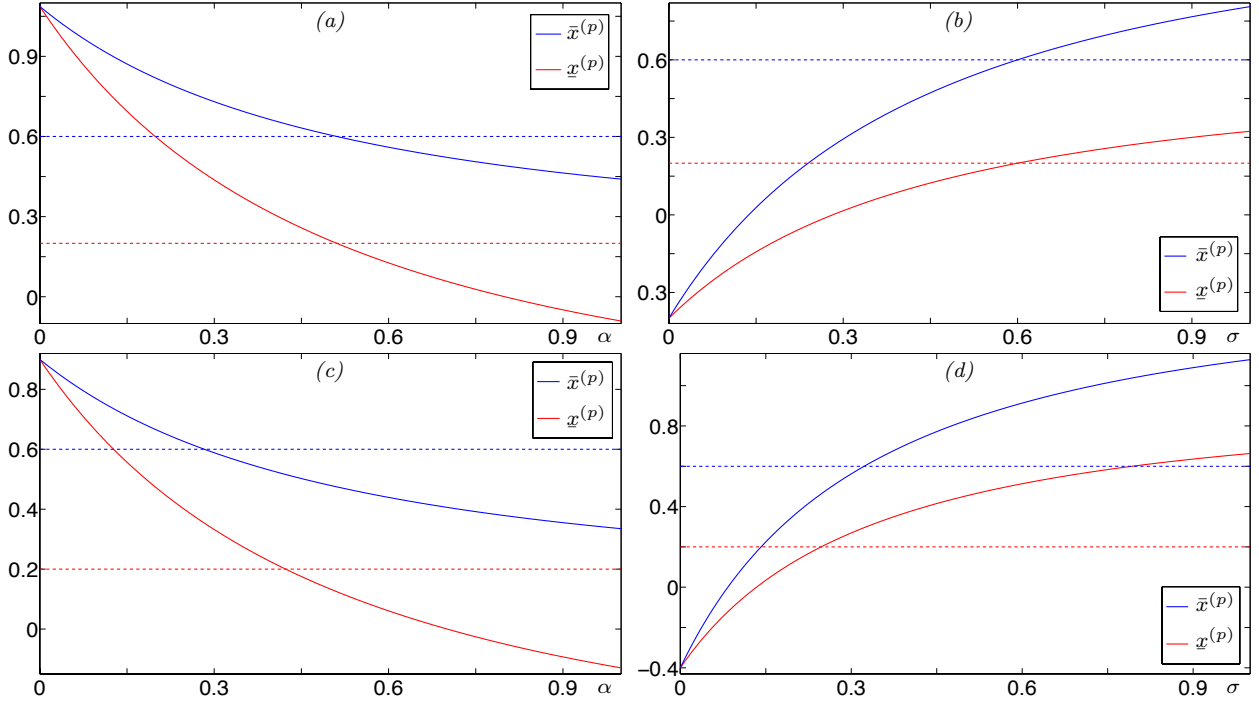

In Figure~\ref{fig_PropLC9} all the four solutions $x^{(p)}$ with different initial perturbation phases $\Delta_{0}$ converge to the same limit cycle, which is given by Proposition~\ref{prop.LC2}, but with different phases.
We also investigated how $\protect\ubar{x}^{(p)}$ and $\bar{x}^{(p)}$ change as function of $\alpha$ in Figure~\ref{fig_osc}(a) and (c) and as function of $\sigma$ in Figure~\ref{fig_osc}(b) and (d) considering the same parameters of Figure~\ref{fig_PropLC9} in panel (a) and (b), but with $a=1.3$ in panel (c) and $a=1.8$ in panel (d). We have  $\protect\ubar{x}^{(p)}=f_{min}x_{norm}$ and $\bar{x}^{(p)}=f_{max}x_{norm}$ for $\alpha=0.510826$ in panel (a) and for $\sigma=0.6$ in panel (b).  In both panels (c) and (d) there is a parameter interval such that $\protect\ubar{x}^{(p)}$ and $\bar{x}^{(p)}$ are bounded above by $x=f_{max}x_{norm}$ and bounded below by $x=f_{min}x_{norm}$.

Several mathematical models have used a DDE similar to~\eqref{dxdt}-\eqref{ftau} to study blood cell dynamics, but with different feedback functions~\eqref{ftau}, and it would be of interest to know how our results compare with those obtained in more complicated (and physiologically more realistic) treatments.  One could think of, for example, the \citet{HEARN1998167} model for canine cyclical neutropenia, or the neutrophil models of \citet{zhuge2012neutrophil} and~\citet{Brooks_2012} which consider GCS-F administration during chemotherapy in humans.

Here we consider the neutrophil model from~\citet{zhuge2012neutrophil}.  Equation (2) from~\citet{zhuge2012neutrophil} describes the dynamics of bone marrow hematopoietic stem cells $(Q)$ and circulating neutrophils $(N)$. We decouple the neutrophil dynamics from the stem cell dynamics by assuming that stem cells are at their normal steady state concentration $Q_{*}=1.1\times10^{6}\,\mathrm{cells/kg}$ of body mass. Thus the second component of equation (2) from~\citet{zhuge2012neutrophil} becomes
\begin{equation}\label{n1}
N^{\prime}(t) = -\gamma_{N}N(t) + A_{N}(t)\kappa_{N}(N(t-\tau_{N}))Q_{*},
\end{equation}
where $\gamma_{N}=2.4\,\mathrm{days}^{-1}$ is the neutrophil apoptosis rate, $A_{N}(t)$ is the amplification factor, $\kappa_{N}(N(t-\tau_{N}))$ is the rate that stem cells  commit to differentiate to neutrophil precursors, and $\tau_{N}$ is the total time it takes to a neutrophil be produced which is defined as the sum of the neutrophil proliferation time $(\tau_{N\!P})$ and the neutrophil maturation time $(\tau_{N\!M})$, i.e. $\tau_{N}= \tau_{N\!P}+\tau_{N\!M}$~\citep{zhuge2012neutrophil}. The neutrophil proliferation phase duration is constant and equal to $5\,\mathrm{days}$, while the neutrophil maturation phase duration depends of the G-CSF serum level. We consider that for a G-CSF dose of $30\mathrm{\mu g/kg/day}$ the neutrophil maturation time is equal to $4.3\,\mathrm{days}$, as estimated by~\citet{zhuge2012neutrophil}. The normal level of circulating neutrophils is taken to be $N_{*}=6.3\times10^{8}\,\mathrm{cells/kg}$ of body mass~\citep{zhuge2012neutrophil}.

Assuming periodic administration of G-CSF the amplification can be expressed as $A_{N}(t)=A_{N}+\xi(t)$, where $A_{N}=6.55\times10^{4}$ and $\xi(t)$ is the periodic perturbation due to the G-CSF ~\citep{zhuge2012neutrophil}.
Without G-CSF the amplification factor is constant and equal to $A_{N}$~\citep{zhuge2012neutrophil}. We can rewrite~\eqref{n1} as
\begin{equation}\label{n2}
N^{\prime}(t) = -\gamma_{N}N(t) + F(N(t-\tau_{N})) + P(t),
\end{equation}
where $F(N(t-\tau_{N}))=A_{N}\kappa_{N}(N(t-\tau_{N}))Q_{*}$ is the delayed negative feedback,  $P(t)$ is  a periodic perturbation assumed to be of the type~\eqref{eq.pert.per.disc}, but with $t$, $\sigma$ and $\alpha$ in units of days and $a$ in units of $\mathrm{cells/kg/day}$. As in~\citet{zhuge2012neutrophil}, we consider that the effects of G-CSF are maintained for one day $(\sigma=1\,\mathrm{day})$ and that the interval between consecutive administrations is also one day $(\alpha=1\,\mathrm{day})$.
We approximate the delayed negative feedback $F(N(t-\tau_{N}))$ by the piecewise constant function
\begin{equation}\label{n3}
F(N(t-\tau_{N})) = \left\{\begin{array}{ll}
b_{L}&\qquad\mbox{for}\qquad N(t-\tau_{N})<N_{*},\\
b_{U}&\qquad\mbox{for}\qquad N(t-\tau_{N})\geq N_{*},
\end{array}\right.
\end{equation}
and take $b_{L}$ equal to the feedback maximum value $A_{N}\kappa_{N}(0)Q_{*}$, where $\kappa_{N}(0)=f_{0}=0.4\,\mathrm{days}^{-1}$~\citep{zhuge2012neutrophil}, and $b_{U}= \varepsilon b_{L}$ with $0<\varepsilon \ll 1$.
%
The change of variables $N(t)\to x(t)$, $\gamma_{N}\to\gamma$, $\tau_{N}\to\tau$, $N_{*}\to\theta$, followed by
the change of variables~\eqref{transform} together with $p(t)=P(t)/\gamma_{N}$ transform Eq.~\eqref{n2} with $F(N(t-\tau_{N}))$ given by~\eqref{n3} into the form~\eqref{eq.diff.pert}. Taking $\varepsilon=10^{-4}$ and the parameter values earlier described we compute the parameters from~\eqref{eq.diff.pert} as being   $\tau=\tau_{N}$, $\gamma_{N}=22.32$, $\sigma=\alpha=\gamma_{N}=2.4\,\mathrm{days}^{-1}$, $\beta_{L}= -N_{*}+b_{L}/\gamma_{N} =11.3783$ and $\beta_{U}=N_{*}-b_{U}/\gamma_{N}=0.6288$, where $\beta_{U}$ and $\beta_{L}$ are in units of $10^{9}\,\mathrm{cells/kg}$ of body mass.

\begin{figure}[!htbp]
\centering
\includegraphics[width=0.73\textwidth]{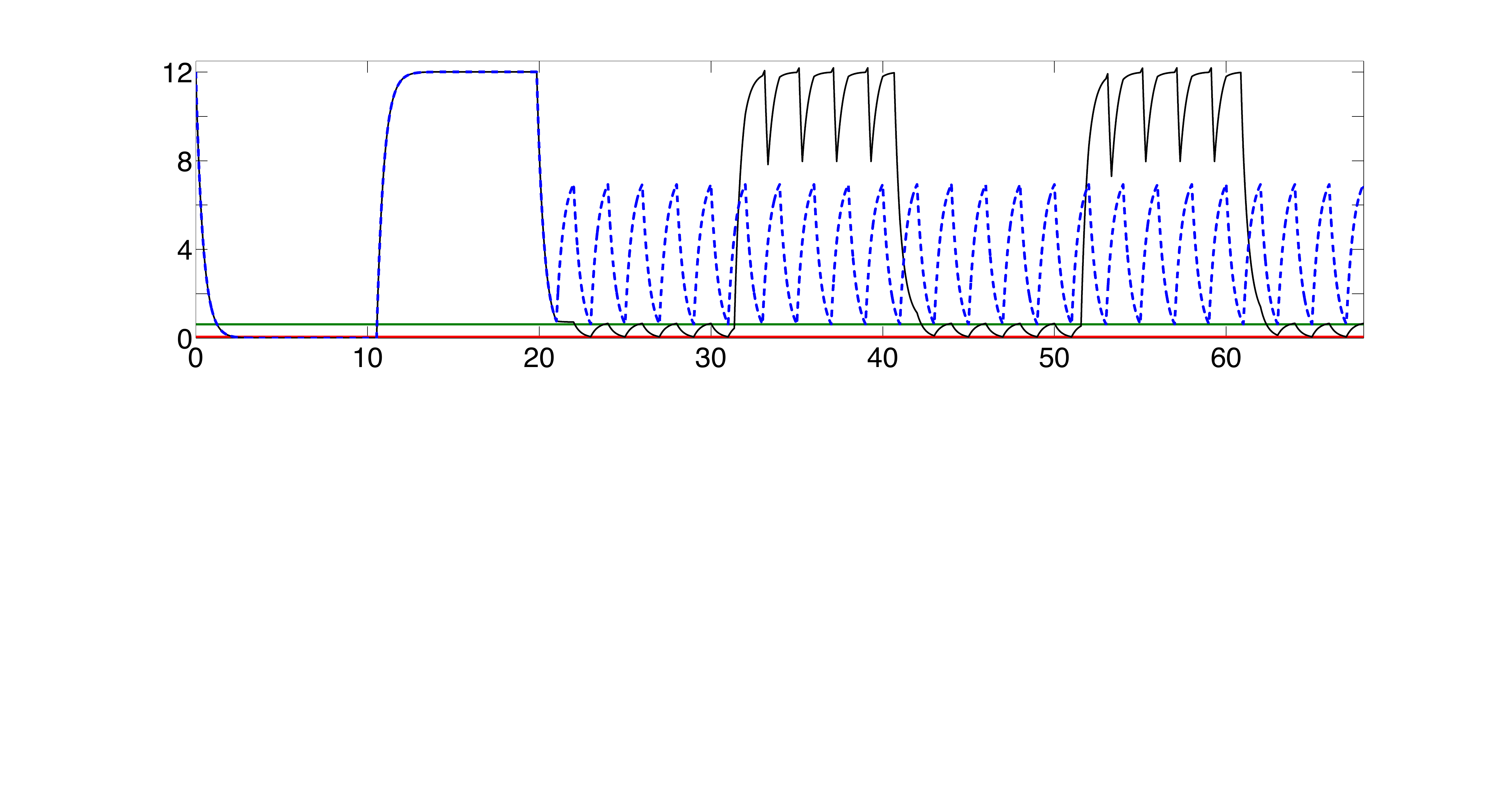}
\put(-341.5,69){\footnotesize{$N$}}
\put(-8,0.8){\footnotesize{$\hat{t}$}}
\caption{For both numerical simulations the neutrophil concentration is given by $N(\hat{t}\hspace{.10em})=x(t)+N_{*}$, where $\hat{t}=t/\gamma_{N}$, $x(t)$ is the solution of~\eqref{eq.diff.pert} for $\beta_{L}=11.3783$, $\beta_{U}=0.6288$, $\tau=22.32$, $\sigma=\alpha=2.4\,\mathrm{days}^{-1}$, $N_{*}=0.63$, and the perturbation begins at day 21, i.e. $\Delta_{0}=21\gamma_{N}$. The perturbation amplitude $a=0.719$ for the black line and $a_{1}=7.5602$ for the blue line, where $a_{1}$ is given by~\eqref{eq.a1}. The green straight line $(N=0.63)$ corresponds to the healthy neutrophil concentration while the red line $(N=0.061)$ is the reference level for severe neutropenia. The time $\hat{t}$ is in units of days while the neutrophil concentration $N(t)$ and the parameters $\beta_{U}$, $\beta_{L}$ and $a$ are in units of $10^{9}\,\mathrm{cells/kg}$ of body mass.}
\label{fig3_Sim_2x2}
\end{figure}

For this set of parameters, $N(t)$ oscillates between $0.0012$ and  $12\times10^{9}\,\mathrm{cells/kg}$ with a period of $19.85$ days  as shown in the first $21$ day portion of
Figure~\ref{fig3_Sim_2x2} where we plot $N(\hat{t}\hspace{.10em})=x(t)+N_{*}$, where $\hat{t}=t/\gamma_{N}$ and $x(t)$ is the solution of~\eqref{eq.diff.pert}.  The remainder of Figure~\ref{fig3_Sim_2x2} for days $21$ to $70$ shows the simulated effects of G-CSF treatments for two different perturbation amplitudes.
The black line corresponds to a G-CSF dose of $a=0.719\times10^{9}\,\mathrm{cells/kg}$ while for the blue line $a=a_{1}=7.5602\times10^{9}\,\mathrm{cells/kg}$.
Both perturbations start at day 21, i.e. $\Delta_{0}=21\gamma_{N}$.
The green line corresponds to the normal neutrophil concentration $(N=N_{*})$ while the red line is the reference level for severe neutropenia, namely $N=0.061\times10^{9}\,\mathrm{cells/kg}$\footnote{Human neutropenia is classified as severe if the neutrophil concentration is below $0.061\times10^{9}\,\mathrm{cells/kg}$ (of body mass), which corresponds to an absolute neutrophil count (ANC) of $500\,\mathrm{cells/\!\mu l}$~\citep{Craig_2015}}.
For the black line the perturbation amplitude $a$ is such that the neutrophil concentrations are greater than or equal to the reference level for severe neutropenia, while for the blue line the amplitude $a$ is such that the neutrophil concentrations are equal to or greater than the normal level $N_{*}$.
For $\hat{t}\geq 21\,\mathrm{days}$ the black line has a maximum value of about $12.2\times10^{9}\,\mathrm{cells/kg}$ (of body mass) which is slightly larger than the maximum value reached without perturbation.


In order to compare the effects of the same G-CSF perturbation in both models, we estimate the amplitude of perturbation $P(t)$ from~\eqref{n2} through the negative feedback function $A_{N}(t)Q_{*}\kappa_{N}(N_{*})\equiv A_{N}(t)Q_{*}f_{0}/2$ from~\eqref{n1} by computing its variation for the value of $A_{N}(t)$ under the G-CSF effects $(A_{N}^{gcs\:\!\!f})$ and without the effects of G-CSF $(A_{N})$. Under the effects of a G-CSF dose of $30\mathrm{\mu g/kg/day}$ the amplification rate $A_{N}(t)$ can be approximated by $A_{N}^{gcs\:\!\!f}=\exp(\eta_{N\!P}^{\max}\tau_{N\!P}-\gamma_{0}^{\min}\tau_{N\!M}^{gcs\:\!\!f})$, where $\eta_{N\!P}^{\max}=3.0552\,\mathrm{days}^{-1}$ is the neutrophil maximal proliferation rate, $\gamma_{0}^{\min}=0.12\,\mathrm{days}^{-1}$ is the neutrophil minimal death rate in maturation, and $\tau_{N\!M}^{gcs\:\!\!f}=4.3\,\mathrm{days}$ is the maturation time, yielding the figure $A_{N}^{gcs\:\!\!f}=2.5715\times10^{6}$~\citep{zhuge2012neutrophil}. This figure together with the change of variables $p(t)=P(t)/\gamma_{N}$ gives the approximation for perturbation amplitude $a\approx (A_{N}^{gcs\:\!\!f}-A_{N})Q_{*}f_{0}/(2\gamma_{N}) \approx 229.72\times10^{9}\,\mathrm{cells/kg}$ of body mass. For this amplitude and considering the other parameters as are in Figure~\ref{fig3_Sim_2x2}, the dynamics of neutrophil concentration $N(\hat{t}\hspace{.10em})$ (not shown) is similar to the behaviour observed for the blue line from Figure~\ref{fig3_Sim_2x2}, it cycles with the same period of the G-CSF administration, but with $N(\hat{t}\hspace{.10em})$ varying from $19.1$ to $210.6\times10^{9}\,\mathrm{cells/kg}$ of body mass. This variation is much larger than the oscillation due to the G-CSF administration obtained in~\citet{zhuge2012neutrophil}, see the first few days of the simulation with red line presented by the authors in Figure 3(b). This perturbation amplitude value also is about 30-fold larger than the amplitude necessary to end neutropenia $(a_{1})$ considered for the simulation with blue line shown in Figure~\ref{fig3_Sim_2x2}.

The discrepancy in the response of both models to the same perturbation amplitude is related with the difference between their negative feedback functions, but it also might be related with the fact and that the response to G-CSF is highly variable for the model from~\citet{zhuge2012neutrophil}, as indicated by the authors.

Consider the limit cycle $N(\hat{t}\hspace{.10em})=x(t)+N_{*}$ shown in Figure~\ref{fig3_Sim_2x2} for $\hat{t}\in[0,21]$, where $\hat{t}=t/\gamma_{N}$ and $x(t)$ is the solution of~\eqref{eq.diff.pert}, the perturbation amplitude is zero $a=0$ and the other parameters are as in Figure~\ref{fig3_Sim_2x2}. Define the minimum, maximum and period of oscillation for $N(\hat{t}\hspace{.10em})$ respectively by $\ubar{N}\coloneqq \ubar{x}+N_{*}$, $\bar{N}\coloneqq \bar{x}+N_{*}$ and $\hat{T}\coloneqq \tilde{T}/\gamma_{N}$, where $\ubar{x}$, $\bar{x}$ and $\tilde{T}$ are given by~\eqref{max.min} and~\eqref{z1z2T}. This together with $\mathrm{e}^{-\tau}\approx 0$ and $\beta_{L}=g\beta_{U}$, where $g\approx 18.09$, yields $\ubar{N}\approx -\beta_{U}+N_{*} = 1.2\times10^{6}\,\mathrm{cells/kg}$ (of body mass), $\bar{N}\approx \beta_{L}+N_{*}=12.0083\times10^{9}\,\mathrm{cells/kg}$ (of body mass) and $\hat{T}=2\tau_{N}+\ln({(g+1)}^{2}/g)/\gamma_{N}=19.85\,\mathrm{days}$. The choice of values for $b_{L}$ and $b_{U}$ from the feedback function~\eqref{n3} defines the parameters $\beta_{L} = -N_{*}+b_{L}/\gamma_{N}$ and $\beta_{U} = N_{*}-b_{U}/\gamma_{N}$ and therefore it is important to define the minimum and maximum of the oscillation, since $\ubar{N}\approx b_{U}/\gamma_{N}$ and $\bar{N}\approx b_{L}/\gamma_{N}$,  but it does not play an important role for defining the oscillation period. Indeed, $\hat{T}$ increases linearly with the delay $\tau_{N}$ and for $g=\beta_{L}/\beta_{U}$ varying from $10^{-3}$ to $10^{3}$ the period stays between $19.1176$ and $21.4791\,\mathrm{days}$, which are close to the periods of $19-21\,\mathrm{days}$ observed for CN~\citep{Colijn_2005b}.

For a wide range of values of the parameters $\beta_{L}$, $\beta_{U}$ (and so $b_{L}$, $b_{U}$) the solution of Eq.~\eqref{n2} with $F(N(t-\tau_{N}))$ given by~\eqref{n3} yields a neutrophil dynamics characteristic of CN, with periods close to $19-21\,\mathrm{days}$ and minimum and maximum of the oscillation given approximately by $\ubar{N}\approx b_{U}/\gamma_{N}$ and $\bar{N}\approx b_{L}/\gamma_{N}$, respectively. So the piecewise linear feedback function~\eqref{n3} may be used as reference to construct nonlinear continuous feedback functions such as those from the mathematical models~\citet{zhuge2012neutrophil} and~\citet{Brooks_2012} in order to model CN dynamics, namely the period, maximum and minimum of the oscillations.

Daily administrations of G-CSF in both humans and grey collies with CN has the effect of reducing the oscillation period and increasing both the neutrophil nadir and the oscillation amplitude~\citep{HEARN1998167,Colijn_2005b}. For the model~\eqref{n2} with $F(N(t-\tau_{N}))$ given by~\eqref{n3} the G-CSF perturbation does increase the neutrophil nadir, as shown in Figure~\ref{fig3_Sim_2x2} simulations. However, for the simulation with black line shown in Figure~\ref{fig3_Sim_2x2} the oscillation amplitude is only slightly increased and the period stays close to the period of oscillation, while for the simulation with blue line both the amplitude and period of oscillation decreases. These results indicates that the perturbation $P(t)$ from~\eqref{n2} does not capture the G-CSF effects on either the period and amplitude from CN oscillations.  
Thus we must conclude that the simple model considered here deviates considerably from the supposedly more physiologically realistic model of~\cite{zhuge2012neutrophil}. While disappointing it is hardly surprising considering the differences in the two models.
%

We can also compare the model effects of periodic administration of chemotherapy for a healthy human by reducing Eq.~\eqref{n2} with $F(N(t-\tau_{N}))$ given by~\eqref{n3} into the form~\eqref{eq.diff.pert}. So the neutrophil dynamics  is given by $N(\hat{t}\hspace{.10em})=x(t)+N_{*}$, where $\hat{t}=t/\gamma_{N}$ and $x(t)$ is the solution of~\eqref{eq.diff.pert}. For this onset of chemotherapy the amplitude of the periodic perturbation $p(t)$ must be negative (only for this onset we assume that $a<0$). We simulate the situation where a normal human receives chemotherapy doses with period of administration varying from 1 to 40 days, as is considered in~\citet{zhuge2012neutrophil}.

As mentioned earlier, for a healthy adult human the normal circulating neutrophil level fluctuates around $0.22-0.85\times{10}^{9}\,\mathrm{cells/kg}$ of body mass. These levels are similar to the neutrophil concentrations shown in Figure 2(b) of~\citet{zhuge2012neutrophil} for the first few days, before the chemotherapy begins. Thus, before chemotherapy we have $p(t)=0$ and we obtain that the neutrophil concentration $N(\hat{t}\hspace{.10em})$ is a limit cycle which oscillates from $\ubar{N}=0.22$ to $\bar{N}=0.85$ by taking $\beta_{U}=0.41$, $\beta_{L}=0.225$, $a=0$ and the other parameters as in Figure~\ref{fig3_Sim_2x2}, with $\ubar{N}$, $\bar{N}$, $\beta_{U}$ and $\beta_{L}$ in units of ${10}^{9}\,\mathrm{cells/kg}$ of body mass. For this set of parameters the period of oscillation is $\hat{T}=19.21\,\mathrm{days}$.
%
\begin{figure}[!t]
\hspace*{-0.9mm}
\includegraphics[width=0.497\textwidth,height=0.278\textwidth]{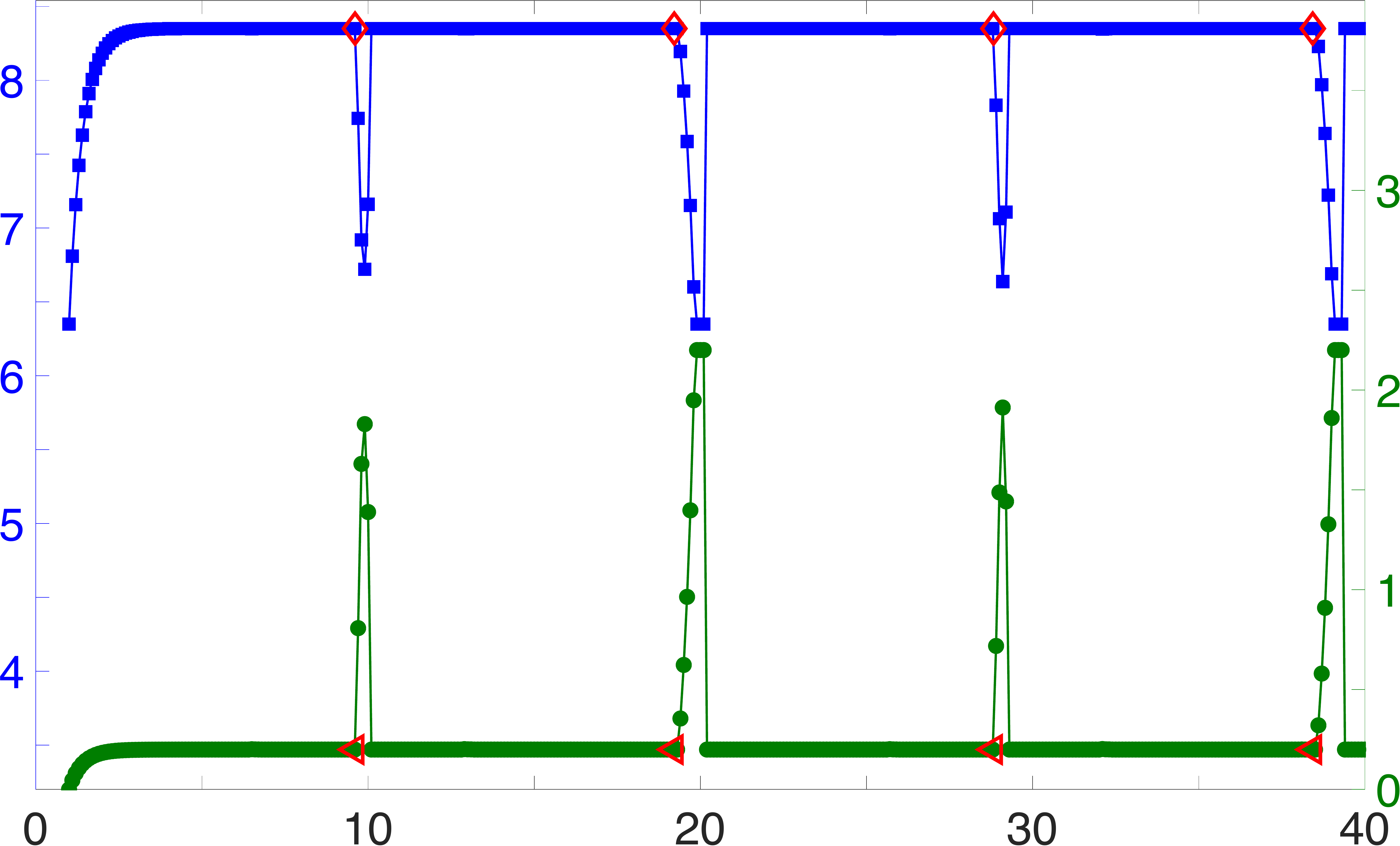}
\put(-151,75){\footnotesize{\textit{(a)}}}
\put(-29,1.5){\scriptsize{$\hat{T}_{p}$}}
\put(-220,36){\rotatebox{90}{\scriptsize{\textcolor{blue}{Amplitude}}}}
\put(-22,56){\rotatebox{90}{\scriptsize{\textcolor{green2}{Nadir}}}}
\hspace*{-0.5mm}
\includegraphics[width=0.497\textwidth,height=0.278\textwidth]{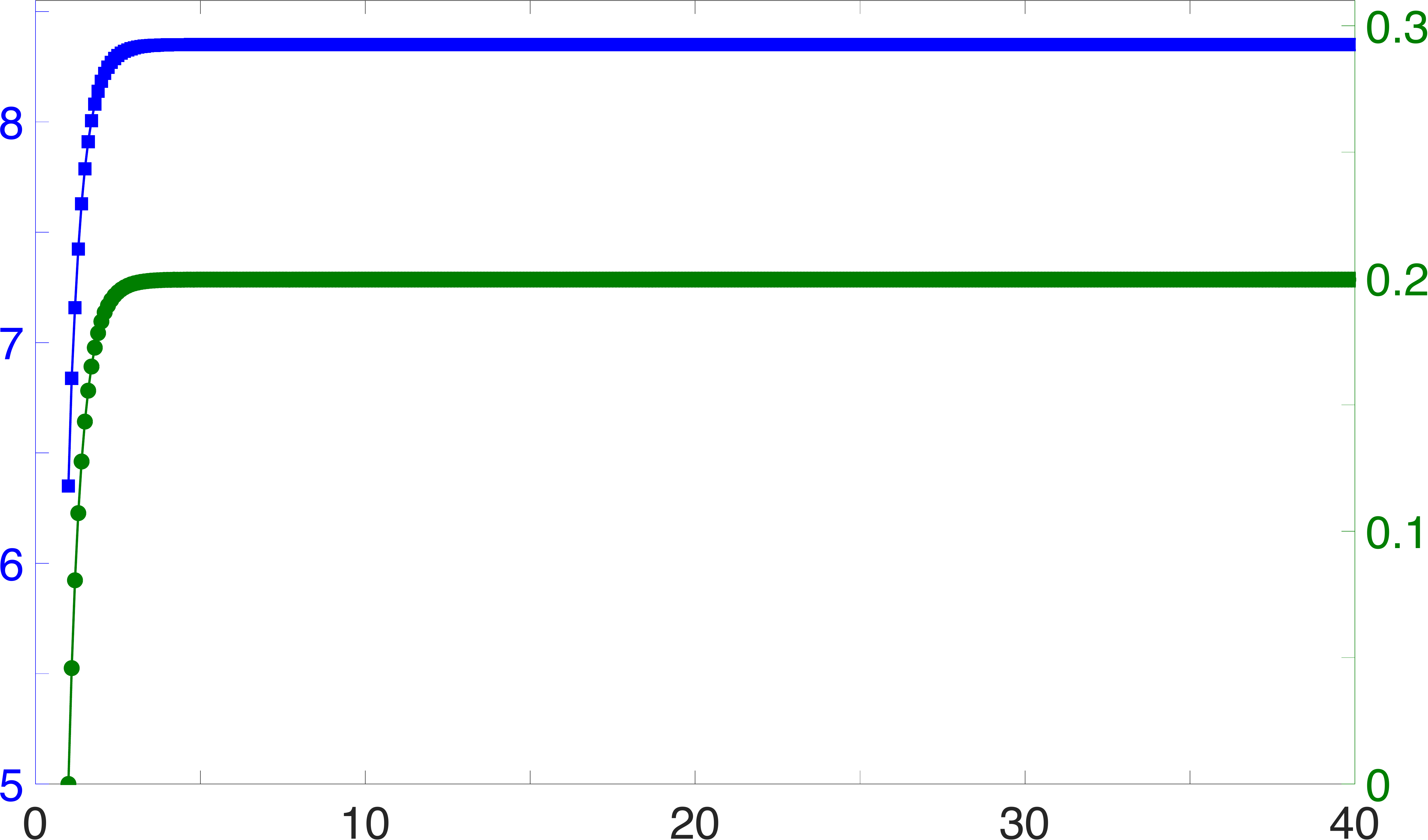}
\put(-151,75){\scriptsize{\textit{(b)}}}
\put(-36,1.5){\scriptsize{$\hat{T}_{p}$}}
\put(-213,37){\rotatebox{90}{\scriptsize{\textcolor{blue}{Amplitude}}}}
\put(-23,48){\rotatebox{90}{\scriptsize{\textcolor{green2}{Nadir}}}}
\caption{Both panels show the Amplitude (left axis and blue squares connected with a solid blue line) and Nadir (right axis and green circles connected with a solid green line) of the neutrophil concentration $N(\hat{t}\hspace{.10em})=x(t)+N_{*}$, where $\hat{t}=t/\gamma_{N}$ and $x(t)$ is the solution of~\eqref{eq.diff.pert}, as function of the period $\hat{T}_{p}$. For both panels $\beta_{L}=0.225$, $\beta_{U}=0.41$, $\tau=22.32$, $\sigma=2.4\,\mathrm{days}^{-1}$, $\Delta_{0}=0$, $N_{*}=0.63$ and the $\hat{T}_{p}$ meshes have 391 equally spaced points from 1 to 40 days. For each mesh point $N(\hat{t}\hspace{.10em})$ was numerically computed for $\hat{t}\in[0,2000]$. The amplitude is given by $(\max(N(\hat{t}\hspace{.10em}))-\min(N(\hat{t}\hspace{.10em})))$ and the nadir is given by $\min(N(\hat{t}\hspace{.10em}))$ with $\hat{t}\in[600,2000]$ for panel (a) and $\hat{t}\in[0,50]$ for panel (b). In panel (a) the triangle and diamond red dots show, respectively, the nadir and amplitude with $\hat{t}\in[600,2000]$ for $\hat{T}_{p}=n\hat{T}/2$ and $n=1,2,3,4$, where $\hat{T}\approx 19.21\,\mathrm{days}$ is the period of the unperturbed limit cycle. The amplitude and nadir are in units of $10^{8}\,\mathrm{cells/kg}$ of (body mass), the period $\hat{T}_{p}$ is in units of days and the parameters $\beta_{U}$, $\beta_{L}$ and $a$ are in units of $10^{9}\,\mathrm{cells/kg}$ of body mass.}
\label{fig_BifAlpha4D}
\end{figure}
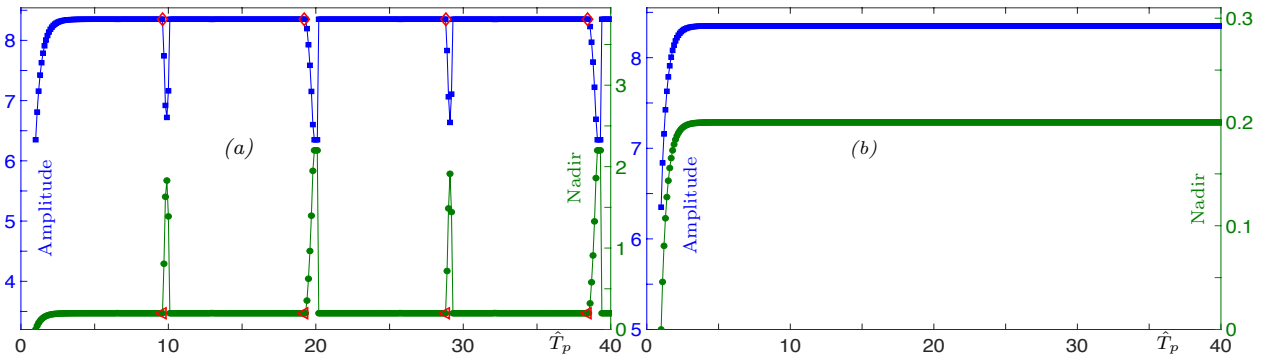

As in~\citet{zhuge2012neutrophil}, we consider that the effects of chemotherapy is maintained for one day, so $\sigma = (1\,\mathrm{day})\gamma_{N}$. We compute the amplitude of perturbation $(a)$ by considering that for a daily administration of chemotherapy the neutrophil nadir $(\ubar{N}^{p})$ is zero, i.e.  $\ubar{N}^{(p)}\coloneqq\ubar{x}^{(p)}+N_{*}\equiv-(\beta_{U}-a)(1-\mathrm{e}^{-\tau})+N_{*}=0$, from where it follows that
\begin{equation}\label{a.min}
a=\beta_{U}-N_{*}/(1-\mathrm{e}^{-\tau})\approx \beta_{U}-N_{*} \approx -0.22\times{10}^{9}\,\mathrm{cells/kg}.
\end{equation}
For this set of parameters the neutrophil concentration $N(\hat{t}\hspace{.10em})$ was numerically computed for $\hat{t}\in[0,2000]$ for each period of administration $\hat{T}_{p}=1,1.1,1.2,\ldots,40\,\mathrm{days}$, where $\hat{T}_{p}\equiv T_{p}/\gamma_{N}$ with $T_{p}=\alpha+\sigma$ and $\sigma$ fixed at $\sigma=(1\,\mathrm{day})\gamma_{N}$. For each period $T_{p}$ Figure~\ref{fig_BifAlpha4D} shows the amplitude $(\max(N(\hat{t}\hspace{.10em}))-\min(N(\hat{t}\hspace{.10em})))$ and nadir $(\min(N(\hat{t}\hspace{.10em}))$ with $\hat{t}\in[600,2000]$ for panel (a) and $\hat{t}\in[0,50]$ for panel (b). For all simulations the first perturbation begins at zero $(\Delta_{0}=0)$ and the history function of~\eqref{eq.diff.pert} is given by
\begin{equation}\label{hist.min}
x(t) = -\beta_{U}+\beta_{U}\mathrm{e}^{-(t-\sigma+\tau)}\qquad\mbox{for}\qquad t\in[-\tau,0],
\end{equation}
which is equivalent to the segment of the limit cycle~\eqref{tildex} for $t\in[\tilde{T}-\sigma-\tau,\tilde{T}-\sigma]$. For $t\in[0,\sigma]$ the solution of~\eqref{eq.diff.pert} is given by
\begin{equation}\label{sol.min}
x(t) = -\beta_{U}+a+(\beta_{U}\mathrm{e}^{\sigma-\tau}-a)\mathrm{e}^{-t}.
\end{equation}
The first perturbation begins at $\Delta_{0}=0$. So at the end of the first perturbation we have $t=\sigma=(1\,\mathrm{day})\gamma_{N}$ and~\eqref{sol.min} yields
\begin{equation}\label{x.min}
x(\sigma) = -\beta_{U}(1-\mathrm{e}^{-\tau})+a(1-\mathrm{e}^{-\sigma})\approx -\beta_{U}+a(1-\mathrm{e}^{-\sigma}).
\end{equation}
With $N(\hat{t})=x(t)+N_{*}$,~\eqref{x.min},~\eqref{a.min} and $\hat{t}=t/\gamma_{N}=\sigma/\gamma_{N}=1\,\mathrm{day}$ we have
\begin{equation}\label{Nmin}
N(1)=x(\sigma)+N_{*}\approx (N_{*}-\beta_{U})\mathrm{e}^{-\sigma} \approx 0.2\times{10}^{8}\,\mathrm{cells/kg}.
\end{equation}

The history function~\eqref{hist.min} considered here gives the shortest nadir possible for $x(t)$ at the end of the first perturbation. Although we have not proved this results mathematically, it is intuitive that for the scenario of a single perturbation and considering the history function~\eqref{hist.min}, the minimal value of the solution~\eqref{eq.diff.pert} as function of the phase $\Delta\in[0,\tilde{T}]$ occurs for $\Delta=0$ and at the end of the perturbation ($t=\sigma$). This minimal value is given by~\eqref{x.min} and the correspondent neutrophil minimal level is given by~\eqref{Nmin}.

In Figure~\ref{fig_BifAlpha4D}(b) the neutrophil nadir converges to the minimal level $N(1)$~\eqref{Nmin}.  For $1\leq\hat{T}_{p}\lesssim 5$ the neutrophil minimum levels after the transient dynamics is smaller than $N(1)$ and the nadir increases with $\hat{T}_{p}$. For $\hat{T}_{p}\gtrsim 5$ the time interval for which the perturbation is turned off $(\hat{T}_{p}-\sigma/\gamma_{N})$ is large enough for the perturbed solution~\eqref{eq.diff.pert} returns to the unperturbed limit cycle~\eqref{tildex} before the perturbation be turned on again. Thus for $\hat{T}_{p}\gtrsim 5$ the minimum levels of neutrophil only depends on which phases of the unperturbed limit cycle~\eqref{tildex} the perturbation is turned on and the resulting nadir is approximately equal to~\eqref{Nmin}.

The amplitude and nadir curves shown in Figure~\ref{fig_BifAlpha4D}(b) are essentially the same when they are computed for $\hat{t}\in[0,2000]$ instead of $\hat{t}\in[0,50]$. Comparing the nadir values and amplitude values for both time intervals and the same mesh $\hat{T}_{p}$, their maximum absolute difference is smaller than $0.0001\times{10}^{8}\,\mathrm{cells/kg}$ of body mass.

In Figure~\ref{fig_BifAlpha4D}(a) the triangular symbols indicate the nadir for the points $\hat{T}_{p}=n\hat{T}/2$ with $n=1,2,3,4$ and show that they are close to the nadir resonance peaks. For the first resonance peak the maximum occurs for $\hat{T}_{p}\approx 9.9$ days and for the third peak it occurs for $\hat{T}_{p}\approx 29.1$ days. The second and fourth resonance peaks are limited above by the nadir of the unperturbed limit cycle $\ubar{N}=2.2\times{10}^{9}\,\mathrm{cells/kg}$ of body mass. Their maximum values stay constant at $\ubar{N}$ along the three points $19.9$, $20$ and $20.1$ days for the third resonance peak, and along the four points $29.1$, $39.1$, $39.2$ and $39.3$ days for the fourth peak. 

The results shown in Figure~\ref{fig_BifAlpha4D}(a) are quantitatively and qualitatively different from the results shown in Figure 2(a) from~\citet{zhuge2012neutrophil}.
In Figure 2(a) from~\citet{zhuge2012neutrophil} the nadir and amplitude vary considerably along the interval $\hat{T}_{p}\in[1,40]$, and at the resonance points the amplitude increases and the nadir decreases.
While in Figure~\ref{fig_BifAlpha4D}(a) the nadir and amplitude increase rapidly for $1\leq\hat{T}_{p}\lesssim 5$ and are essentially constant for $5\gtrsim \hat{T}_{p} \leq 40$, except at the 4 narrow resonance peaks, where the amplitude decreases and the nadir increases substantially.

Administration protocols of common chemotherapeutic drugs (such as docetaxel, cyclophosphamide, cisplatin, paclitaxel, etc.) often prescribe a chemotherapy cycle of three weeks $(\hat{T}_{p}=21\,\mathrm{days})$~\citep{zhuge2012neutrophil}.
Our simulations in the current model would  suggest that chemotherapy treatments with cycles close to 21 days do not affect the nadir levels, except that for cycles inside of the narrow interval from about 19.9 to 20.1 days would increase the neutrophil nadir along the treatment.

So, again, the results of the simple model considered in this paper seem to be at odds with the results of the model of~\cite{zhuge2012neutrophil}.
%
%
%
\section{Bifurcations in the Face of Periodic Perturbation}\label{sec.num.exp.per.pert}
Periodic perturbation can give rise to periodic solutions through a variety of bifurcations, and in this section we explore these. We examined these mechanisms and how the local maxima and minima of the solution change by performing one parameter continuation on the extremal points.
\begin{figure}[!t]
\includegraphics[width=0.49\textwidth,height=0.278\textwidth]{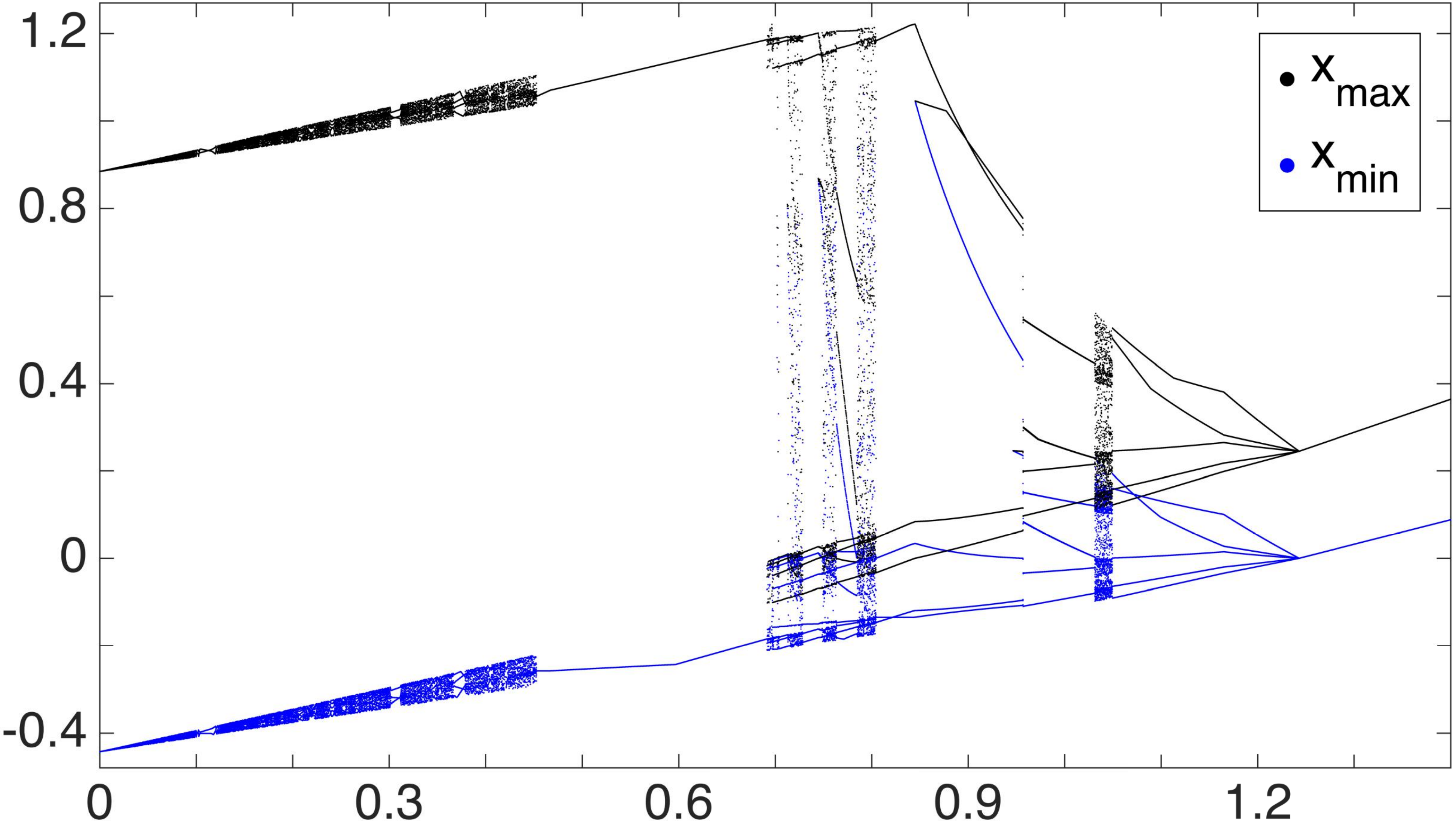}
\put(-209,68){\footnotesize{\textit{(a)}}}
\put(-9,1.5){\footnotesize{$a$}}
\hspace*{1mm}
\includegraphics[width=0.49\textwidth,height=0.278\textwidth]{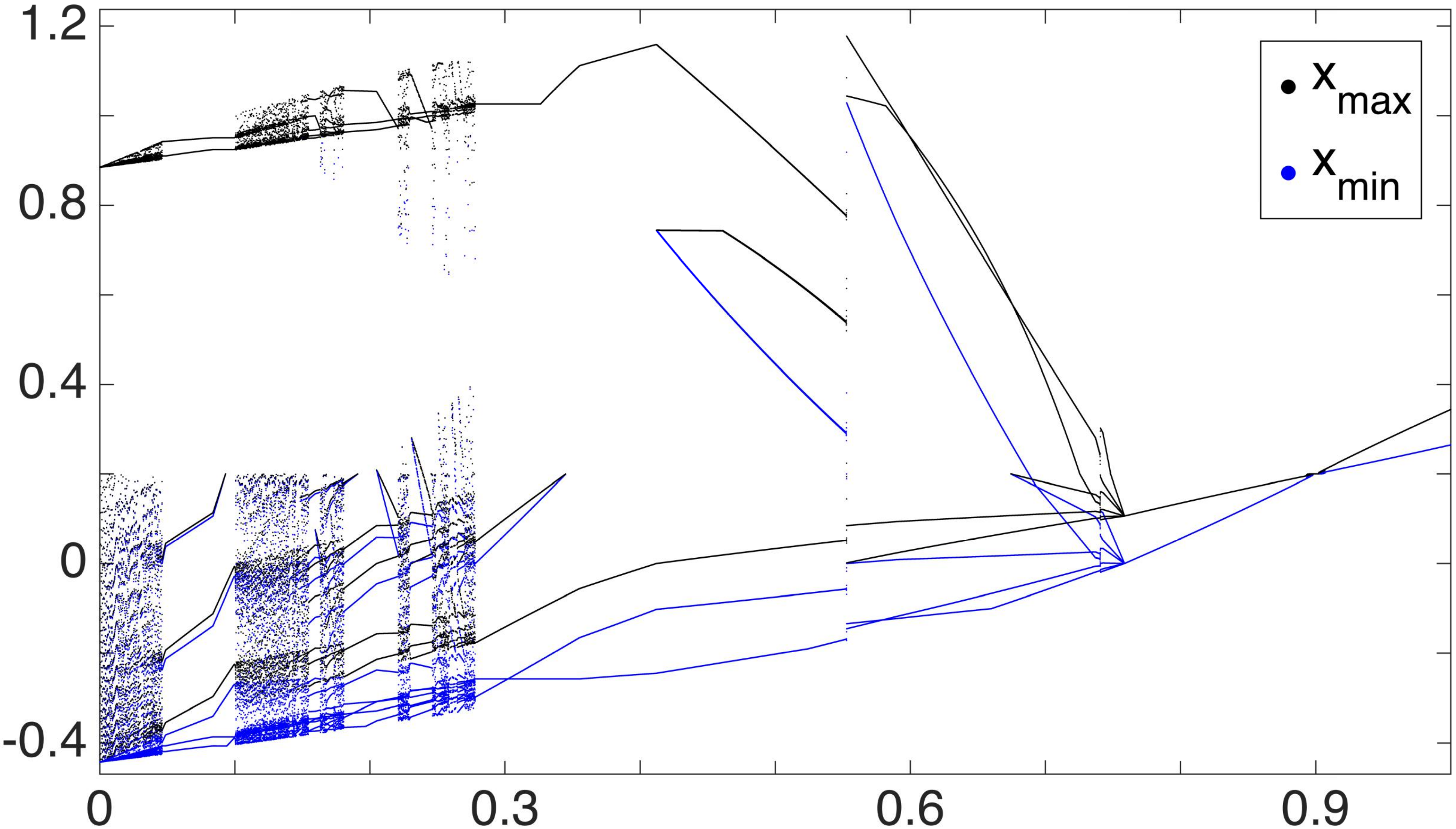}
\put(-209,68){\footnotesize{\textit{(b)}}}
\put(-9,1.5){\footnotesize{$\sigma$}}\\
\includegraphics[width=0.49\textwidth,height=0.278\textwidth]{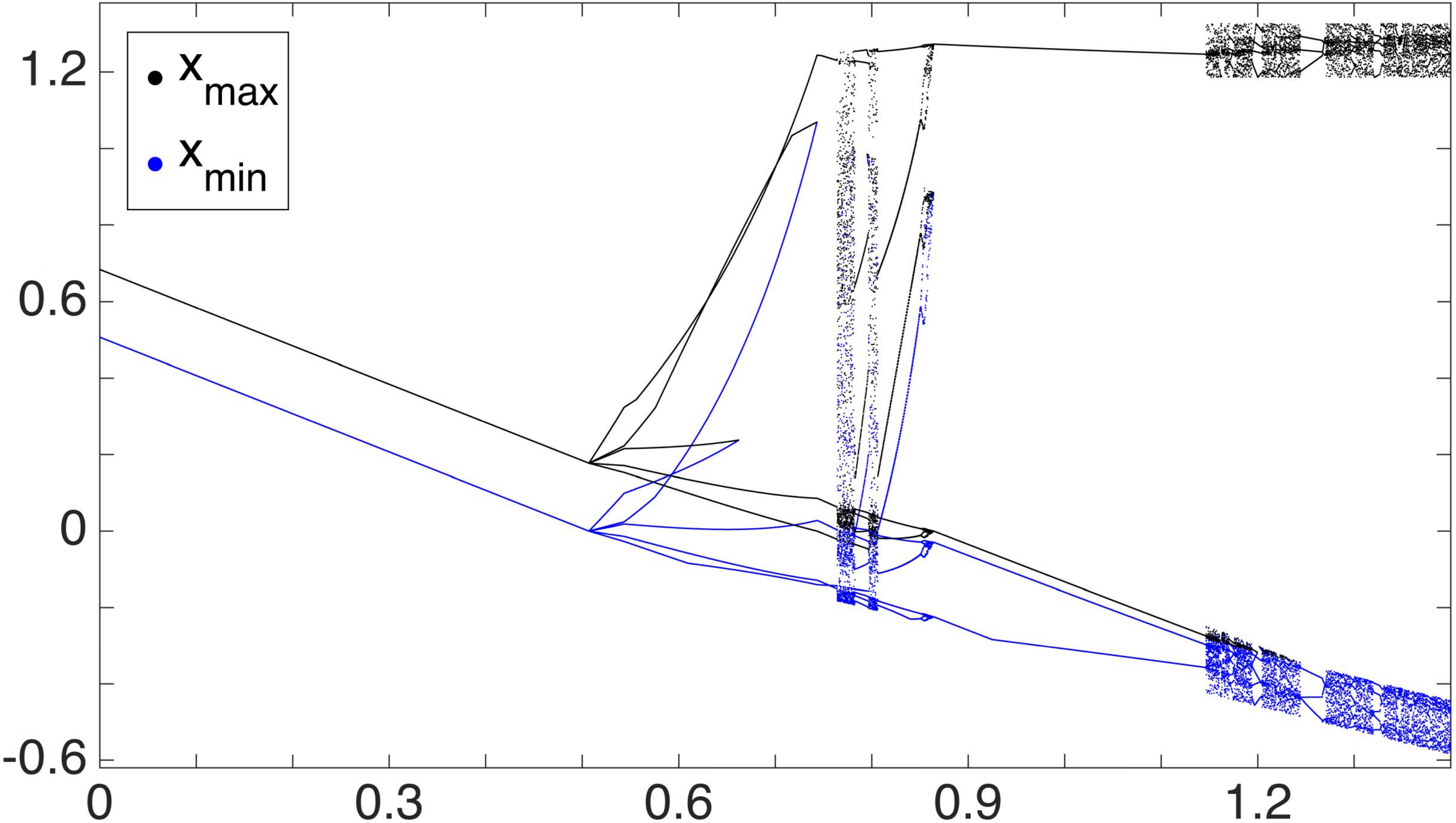}
\put(-209,56){\footnotesize{\textit{(c)}}}
\put(-12,2.1){\footnotesize{$\beta_{U}$}}
\hspace*{1mm}
\includegraphics[width=0.49\textwidth,height=0.278\textwidth]{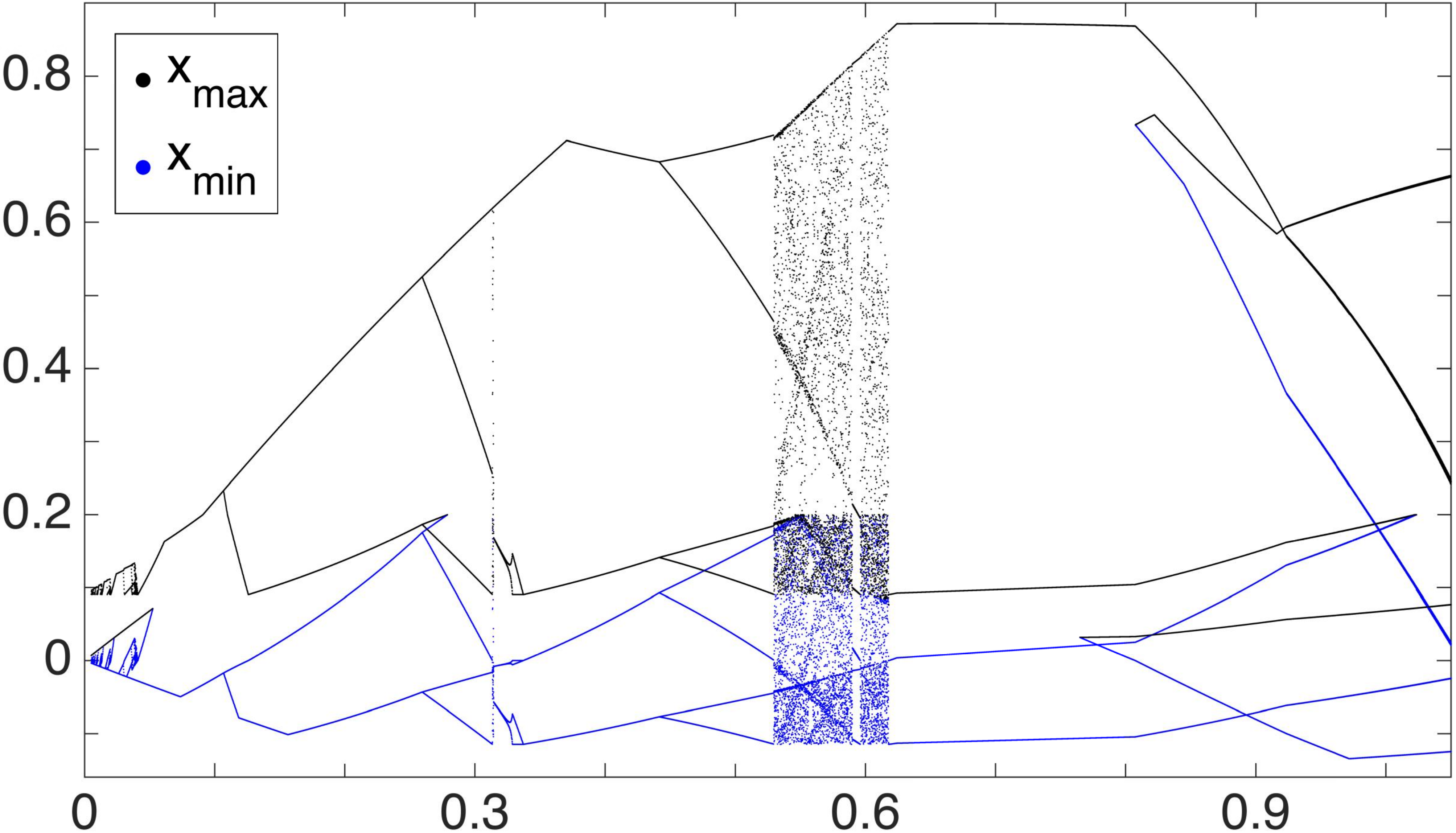}
\put(-212,59){\footnotesize{\textit{(d)}}}
\put(-9,1.5){\footnotesize{$\tau$}}
\vspace*{-02mm}
\caption{Bifurcation diagrams showing local maxima and minima of numerical solution segments of~\eqref{eq.diff.pert} with the parameters $\beta_{L}=1.4$, $\beta_{U}=0.7$, $\sigma=0.6$, $a=0.9$, $\tau=1$, $\alpha=0.3$. Black and blue dots represent respectively the local maxima and minima computed along a mesh with $10^{4}$ points for increasing parameter values for panels (a) and (c) and for decreasing parameter values for panels (b) and (d). We varied $a$ from 0 to 1.4 in (a), $\sigma$ from 0 to 1 ($T_{p}=\sigma+\alpha$ varies from 0.3 to 1.3) in (b), $\beta_{U}$  from 0.5 to 0.9 in (c), and $\tau$ from 0.005 to 1.3 in (d).}
\label{fig_Bif}
\end{figure}
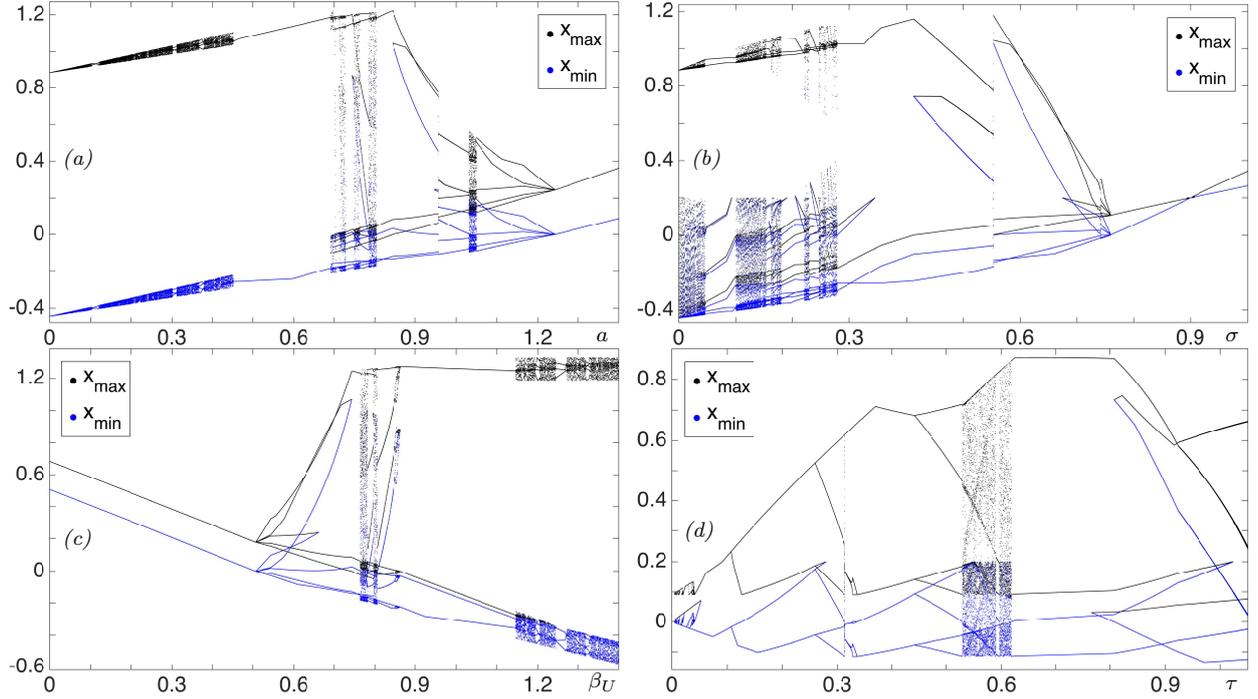

In Figure~\ref{fig_Bif} we present orbit diagrams for~\eqref{eq.diff.pert} as one of each of the parameters $\{a,\sigma,\beta_{U},\tau\}$ is varied. Orbit diagrams are normally produced for maps, but we can reduce the solution of~\eqref{eq.diff.pert} to a map by considering crossings of a Poincar\'{e} section which contains the local maxima and minima of $x(t)$ along the solution~\citep{Daniel_Tony}.

All graphs of Figure~\ref{fig_Bif} were constructed with the following technique. A mesh with $10^{4}$ points was used to compute the solutions for increasing parameter values (noted on each abscissa)   for panels (a) and (c) and for decreasing parameter values for panels (b) and (d). For all mesh points the orbits were computed using the MATLAB \texttt{dde23} routine~\citep{Matlab}, with an absolute error of $10^{-9}$ and relative error of $10^{-9}$. The points where the trajectory is discontinuous were detected and included in the solution mesh by using the MATLAB \texttt{events} function.  For each mesh point we integrated through a transient of $5.5T_{p}$, and then plotted all the maxima and minima that occur over the next $5.5T_{p}$.  The  last $\tau$ time units of the solution is used as the history function to compute the solution at the next mesh point. We also used the function \texttt{Jumps} of the MATLAB \texttt{dde23} routine to include the discontinuity points  in each history function mesh to compute the solution of the next mesh point. For the first mesh point we integrated through a transient of length $220T_{p}$.

In all panels of Figure~\ref{fig_Bif}, the one parameter continuation reveals numerous points of period-doubling bifurcation of periodic orbits and several parameter intervals of periodic dynamics and windows of irregular motion. Panels (a) and (b) respectively show that the smallest minima of the solutions increases, except for small variations, as the perturbation amplitude $a$ or the time duration of perturbation $\sigma$ are increased. Conversely, panels (c) and (d) respectively show that the smallest minima of the solutions decreases, except for small variations, as $\beta_{U}$ or the delay $\tau$ are increased.

Figure~\ref{fig_Bif}(a) shows discontinuities in the extrema for $a\approx 0.957$, while panel (b) shows discontinuities for $\sigma\approx 0.553$ and $\sigma\approx 0.741$. These discontinuities are due to numerical issues and must disappear by increasing the integration time and the mesh size. However, it was verified that doubling the the integration time and/or the mesh size was not enough to remove these discontinuities. While it would be interesting to obtain these graphs with continuous extrema, it might take a long time to find a suitable integration time and mesh size because every attempt requires to solve the DDE~\eqref{eq.diff.pert} numerically with a reasonable precision and also include the breaking points in the solution meshes. We have not pursued this issue in this paper.

\begin{figure}[!t]
\centering
\includegraphics[width=0.49\textwidth,height=0.278\textwidth]{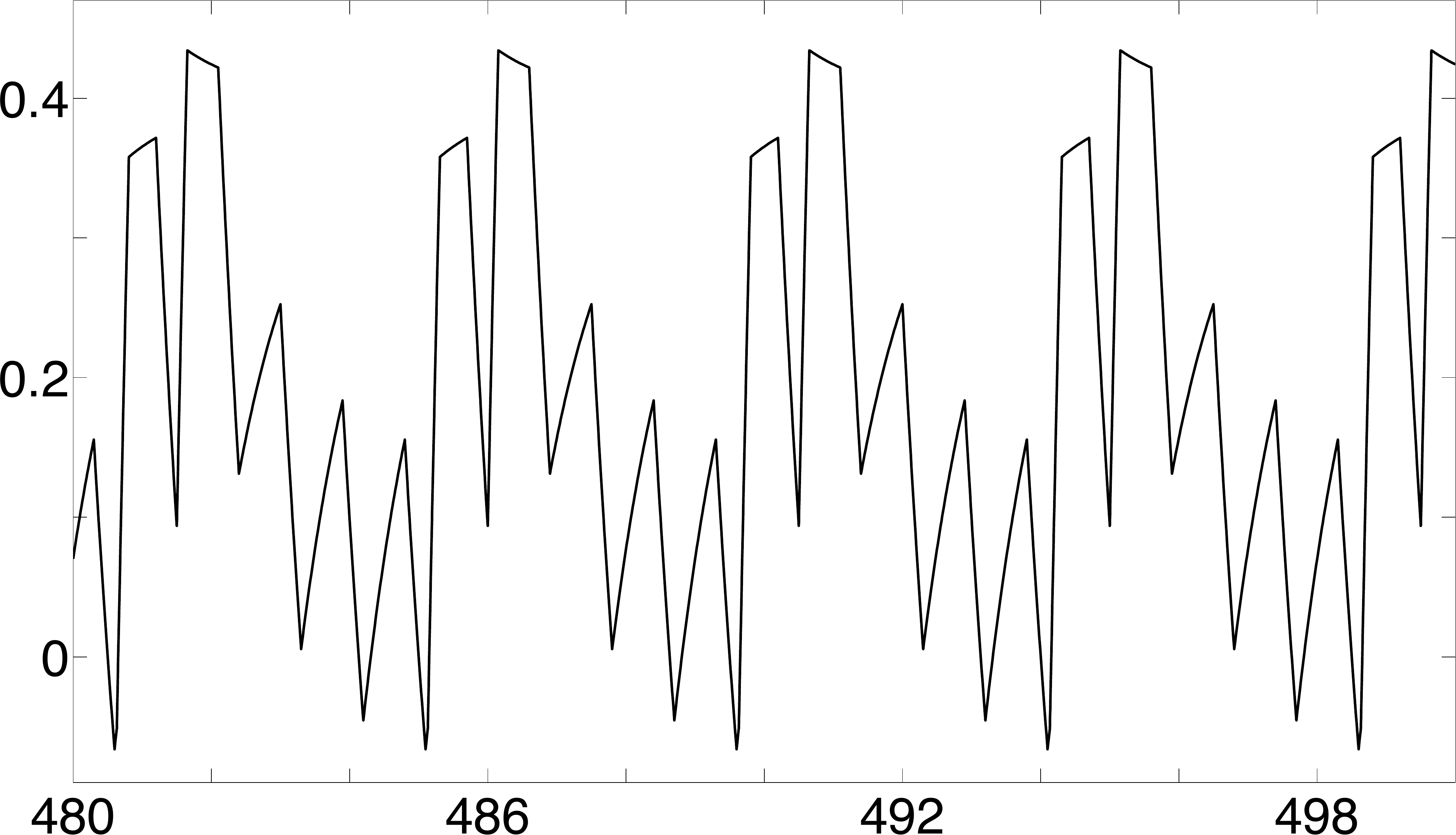}
\put(-18.3,123){\footnotesize{\textit{(a)}}}
\put(-5,0.5){\footnotesize{$t$}}
\put(-225,124){\footnotesize{$x$}}
\hspace*{0mm}
\includegraphics[width=0.49\textwidth,height=0.278\textwidth]{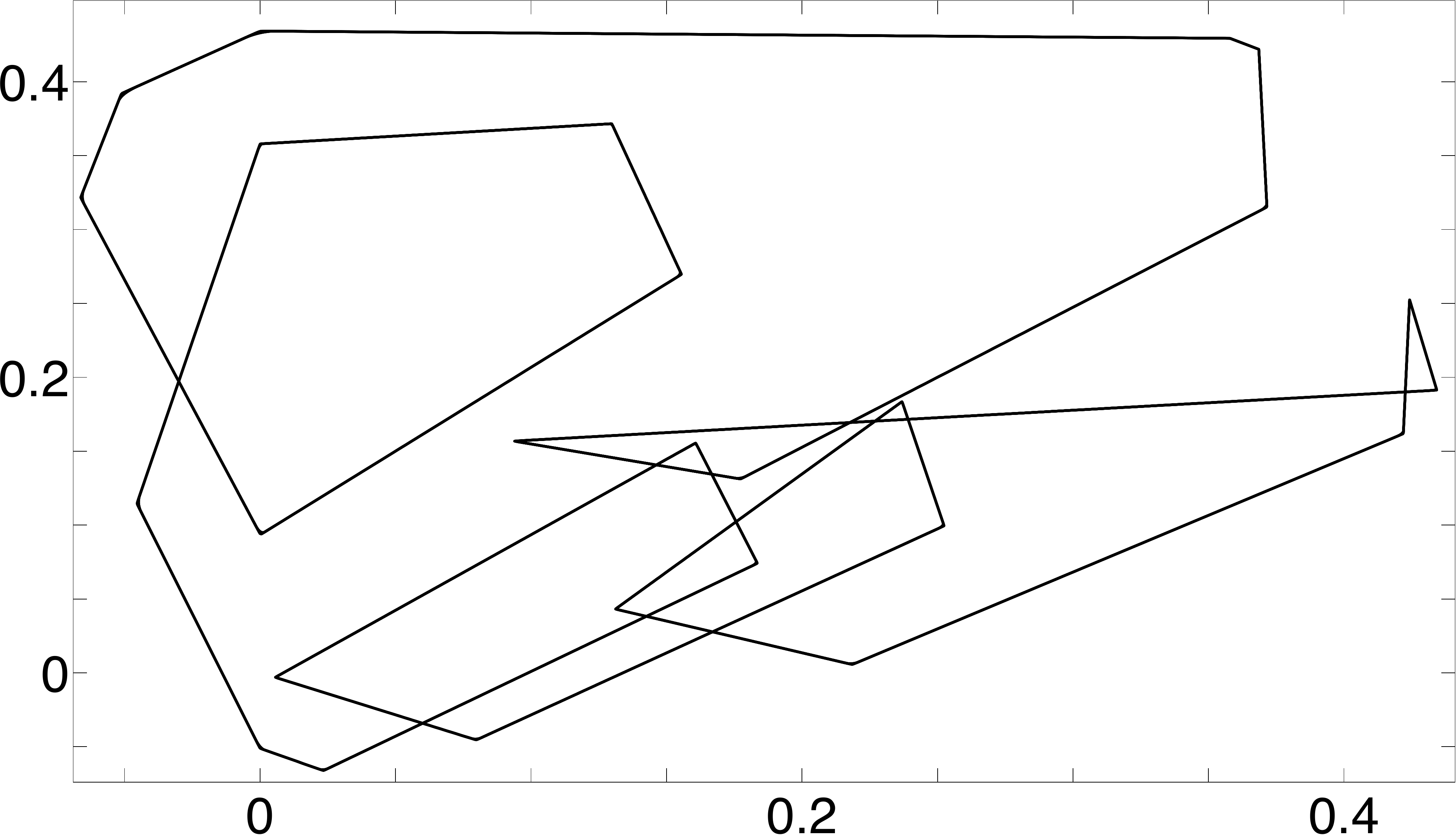}
\put(-15,120){\footnotesize{\textit{(b)}}}
\put(-9.8,1.5){\footnotesize{$x_{\tau}$}}
\put(-225,124){\footnotesize{$x$}}\\
\vspace*{1mm}
\includegraphics[width=0.49\textwidth,height=0.278\textwidth]{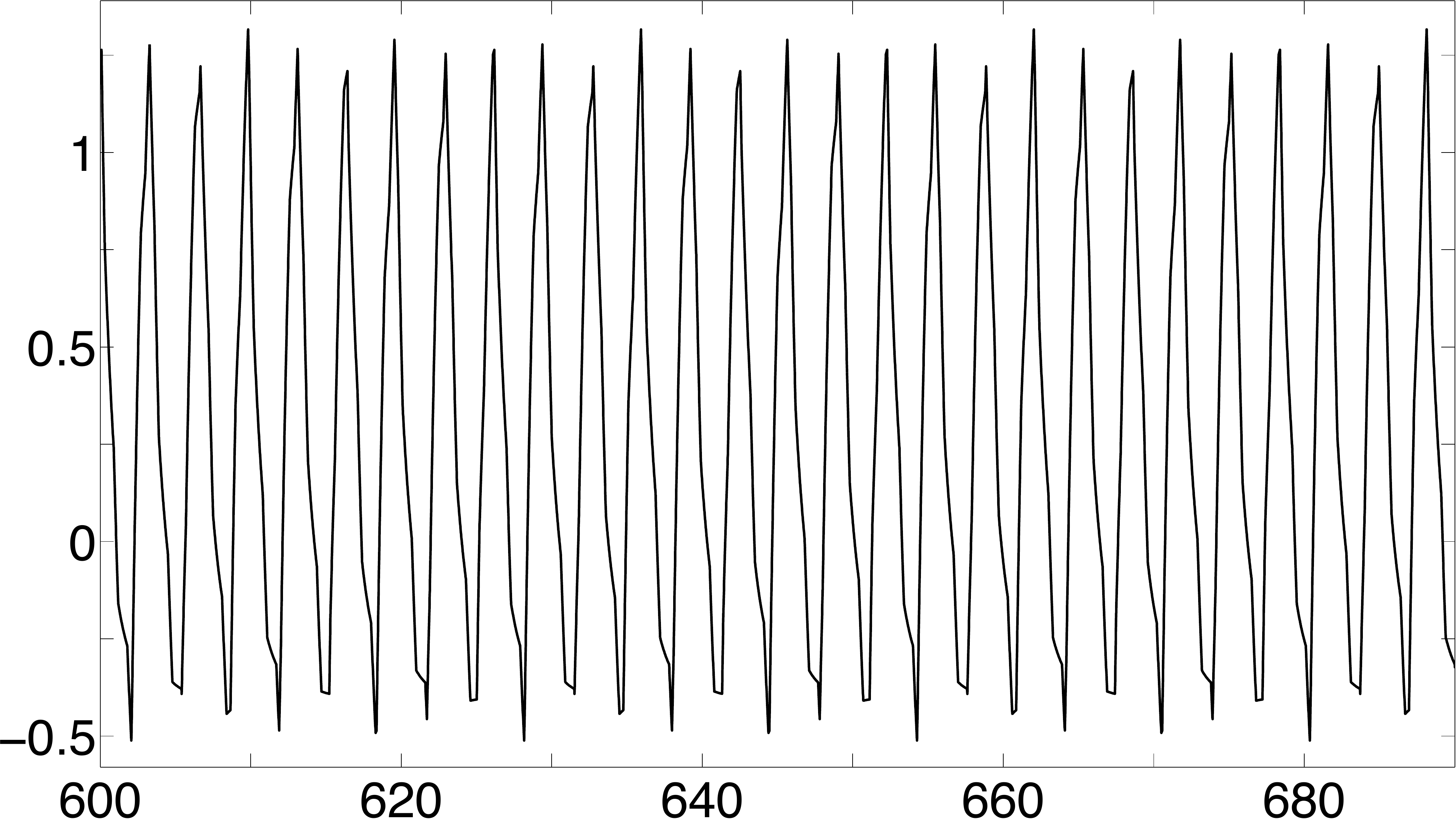}
\put(-18,121){\footnotesize{\textit{(c)}}}
\put(-5,0.5){\footnotesize{$t$}}
\put(-221,124){\footnotesize{$x$}}
\hspace*{0mm}
\includegraphics[width=0.49\textwidth,height=0.278\textwidth]{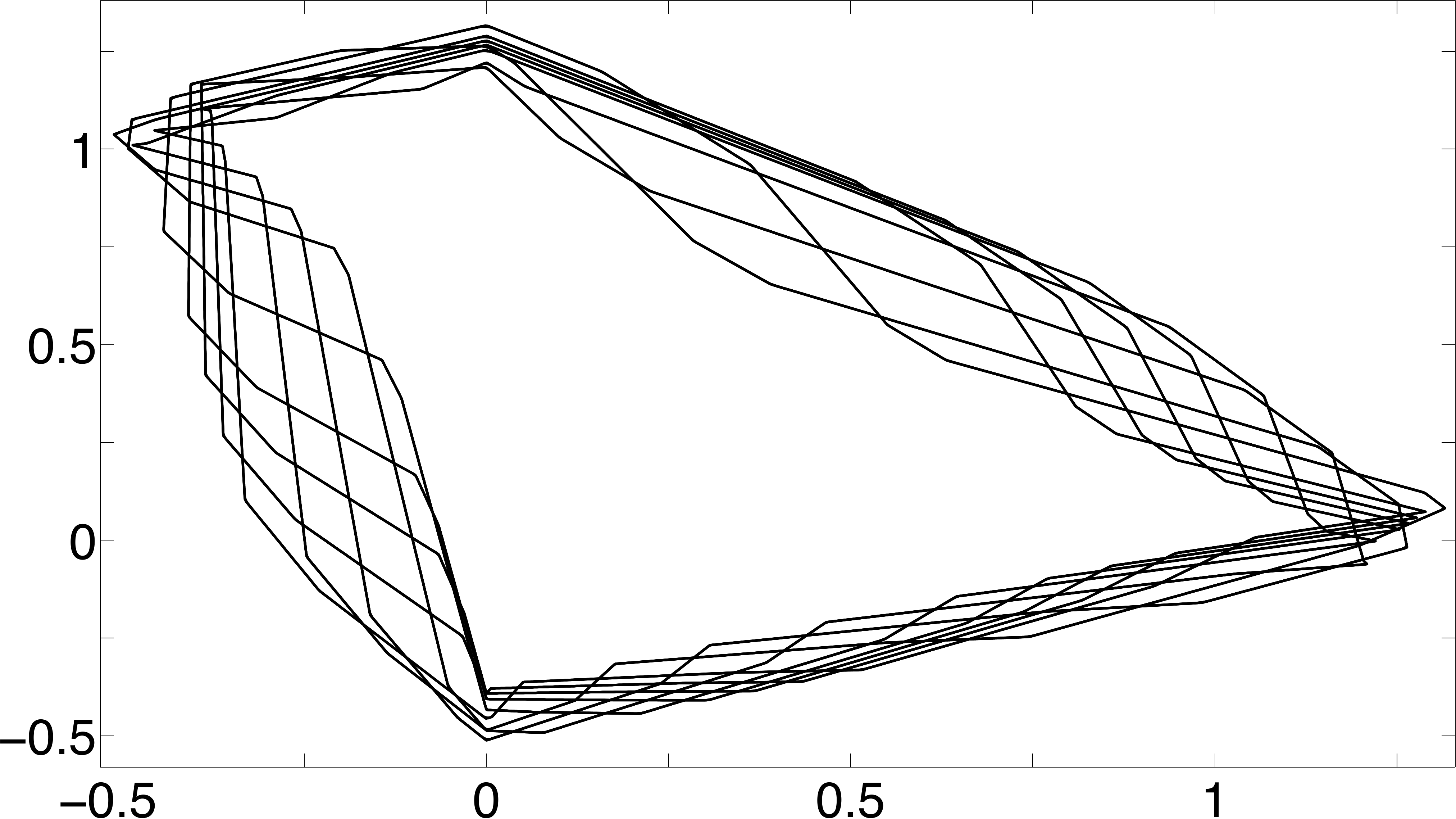}
\put(-15,120){\footnotesize{\textit{(d)}}}
\put(-10,1.5){\footnotesize{$x_{\tau}$}}
\put(-221,124){\footnotesize{$x$}}\\
\vspace*{1mm}
\includegraphics[width=0.49\textwidth,height=0.278\textwidth]{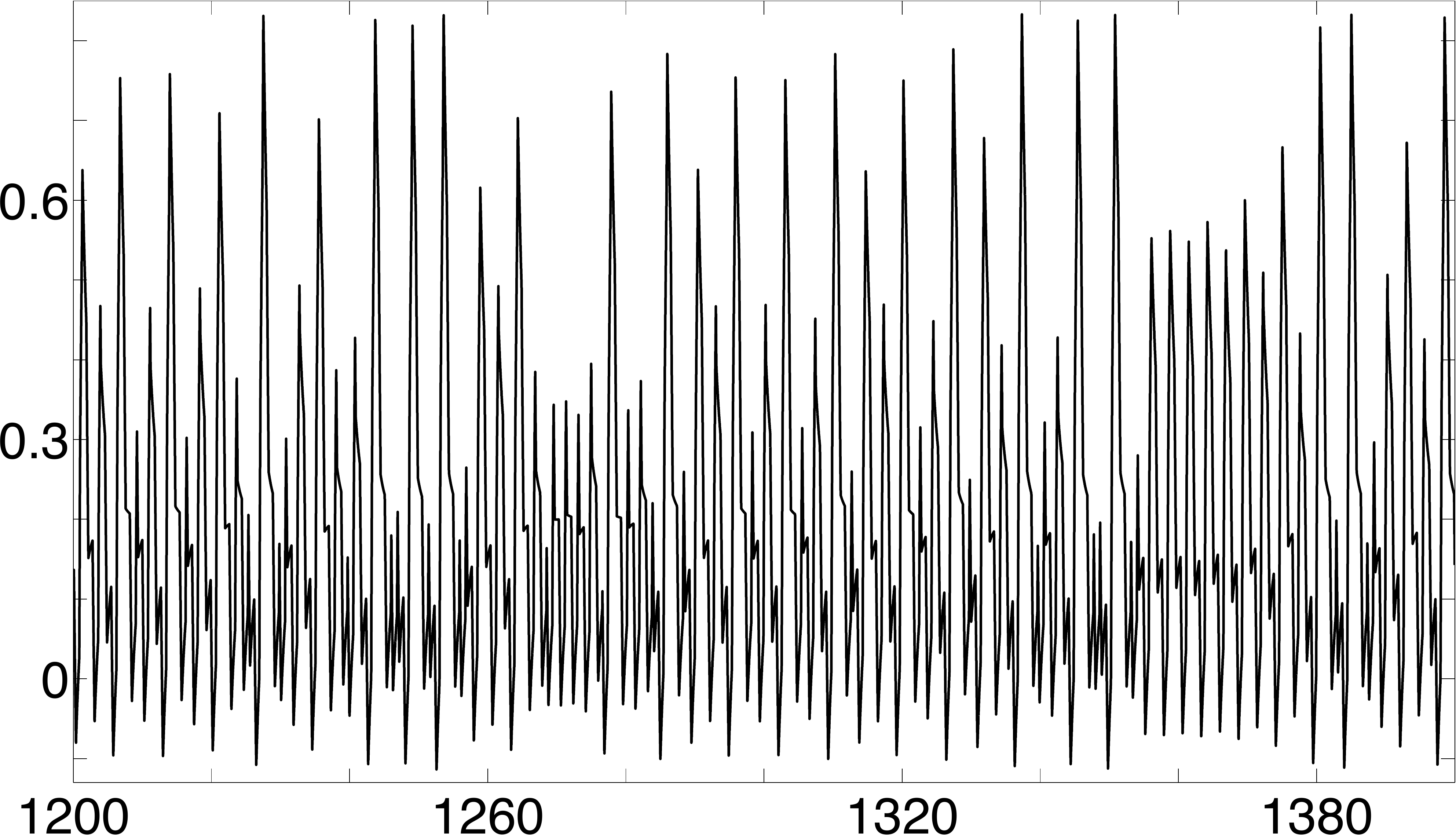}
\put(-17,121){\footnotesize{\textit{(e)}}}
\put(-5,0.5){\footnotesize{$t$}}
\put(-225,125){\footnotesize{$x$}}
\hspace*{0mm}
\includegraphics[width=0.49\textwidth,height=0.278\textwidth]{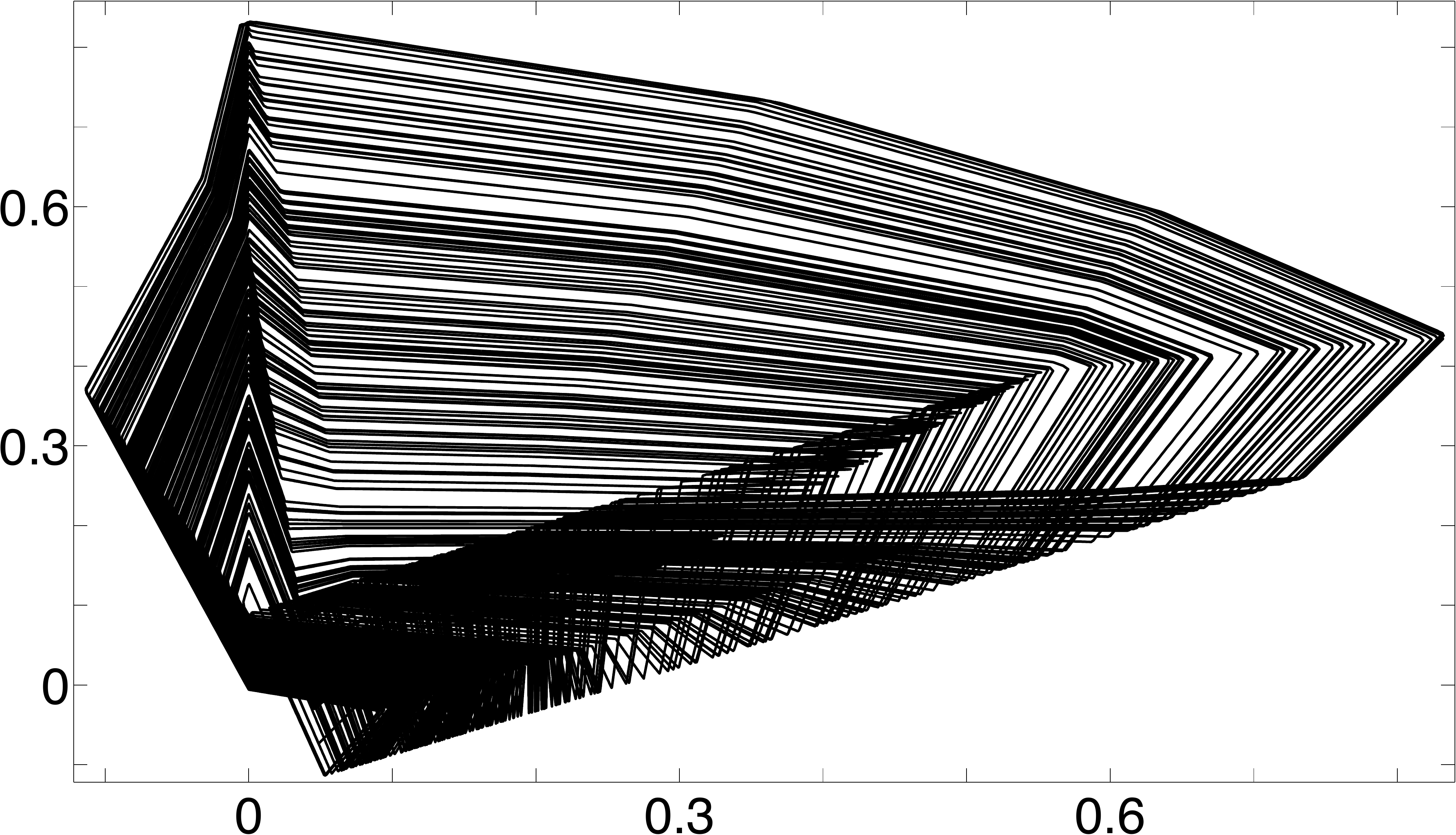}
\put(-15,121){\footnotesize{\textit{(f)}}}
\put(-10,1.5){\footnotesize{$x_{\tau}$}}
\put(-225,124){\footnotesize{$x$}}
\vspace*{-02mm}
\caption{Here we show three numerical solutions of~\eqref{eq.diff.pert} in the left panels and their respective time-delay embeddings $x^{(p)}(t-\tau)\times x^{(p)}(t)$ in the right panels. To construct the embeddings we integrated the solution through a transient for $t\in[0,200]$ and computed $x(t)$ and $x(t-\tau)$ for $t\in[200,400]$.   For panels (a)-(b) we took  $a=1.1$ and used the interval $t\in[450,500]$ for the embedding. For panel for (c)-(d) we took $\beta_{U}=1.3$ and used the interval $t\in[600,700]$ for the embedding. For panel for (e)-(f) we have $\tau=0.6$ and used the interval $t\in[1100,1400]$ for the embedding. All other parameters were taken as in Figure~\ref{fig_Bif}.}
\label{fig_PS}
\end{figure}
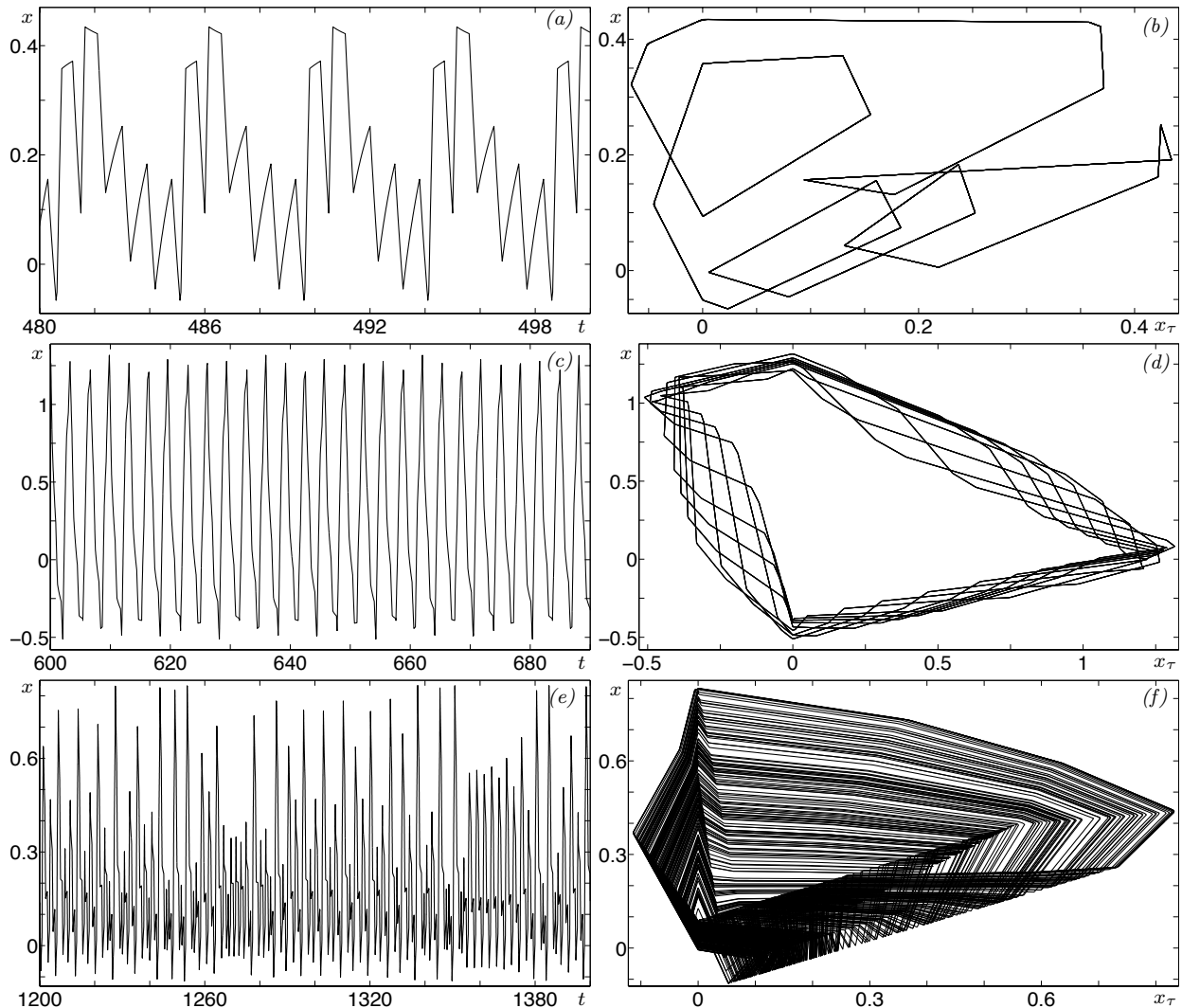

In the left side of Figure~\ref{fig_Bif}(c) we see a parameter interval with a simple limit cycle followed by intervals with period-5, -4 and -3 limit cycles, where the period-3 region ends with an abrupt transition to irregular motion. A small parameter change can thus cause  periodic motion to become irregular, and vice versa.

Figure~\ref{fig_PS} shows a period-5 limit cycle in panels (a)-(b), a period-8 limit cycle in panels (c)-(d) and a complex solution in panels (e)-(f). For all three orbits we used $\Delta_{0}=z_{2}$ and we have $T_{p}<\bar{F}$ and $a<a_{1}$. Orbits from Figure~\ref{fig_PS}(a)-(b) correspond to solutions of Figure~\ref{fig_Bif}(a) for $a=1.1$ and orbits from Figure~\ref{fig_PS}(c)-(d) are related with solutions of Figure~\ref{fig_Bif}(c) for $\beta_{U}=1.3$. The solutions shown in Figure~\ref{fig_PS}(e)-(f) correspond to solutions with $\tau=0.6$ inside of the windows of apparent chaotic dynamics of Figure~\ref{fig_Bif}(d).

For some parameters $\{\tau,\beta_{U},\beta_{L},\sigma,a,\sigma,\alpha\}$, there can  occur frequency locking between the perturbation period $(T_{p})$ and the limit cycle period $(T)$. For example, for the limit cycle shown in Figure~\ref{fig_PS}(a)-(b) we have $T=4.5$ and $T_{p}=0.9$ which yields a frequency locking 5:1, while for the solution of Figure~\ref{fig_PS}(c)-(d) we have $T=26.1$ and $T_{p}=0.9$, with a frequency locking of 29:1. The solutions shown in Figures~\ref{fig_PropLC4x4} and~\ref{fig_PropLC9} converge to limit cycles where $T=T_{p}$, and a frequency locking of 1:1. In fact, the frequency locking 1:1 occurs for all limit cycles with $a\geq a_{1}$ (see Proposition~\ref{prop.LC2} and~\ref{prop.LC}).

\begin{figure}[!t]
\includegraphics[width=0.49\textwidth,height=0.276\textwidth]{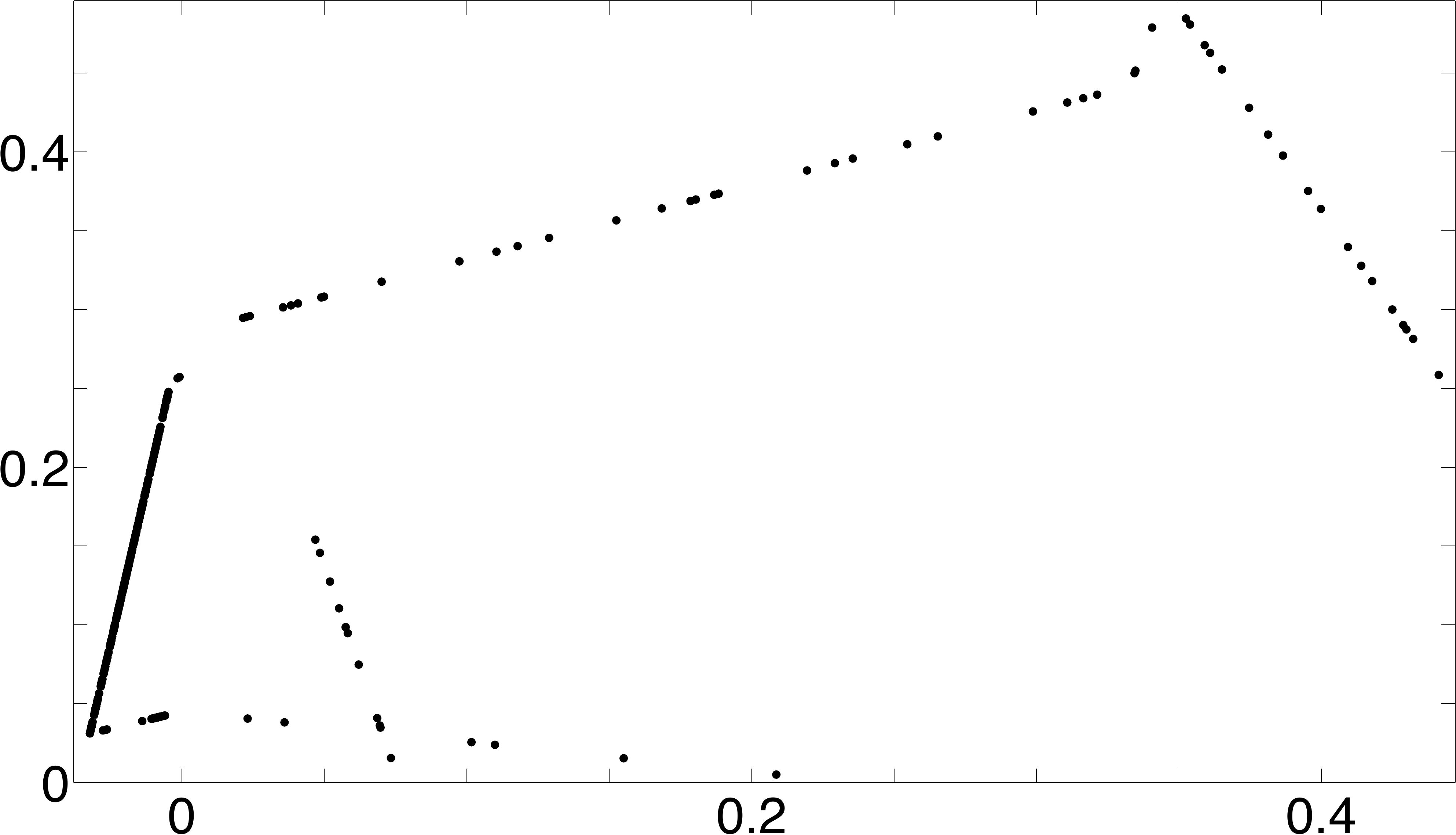}
\put(-15,120){\footnotesize{\textit{(a)}}}
\put(-80.8,1.5){\footnotesize{$x(t-\tau)$}}
\put(-228,64.5){\rotatebox{90}{\footnotesize{$x(t-2\tau)$}}}
\hspace*{1mm}
\includegraphics[width=0.49\textwidth,height=0.276\textwidth]{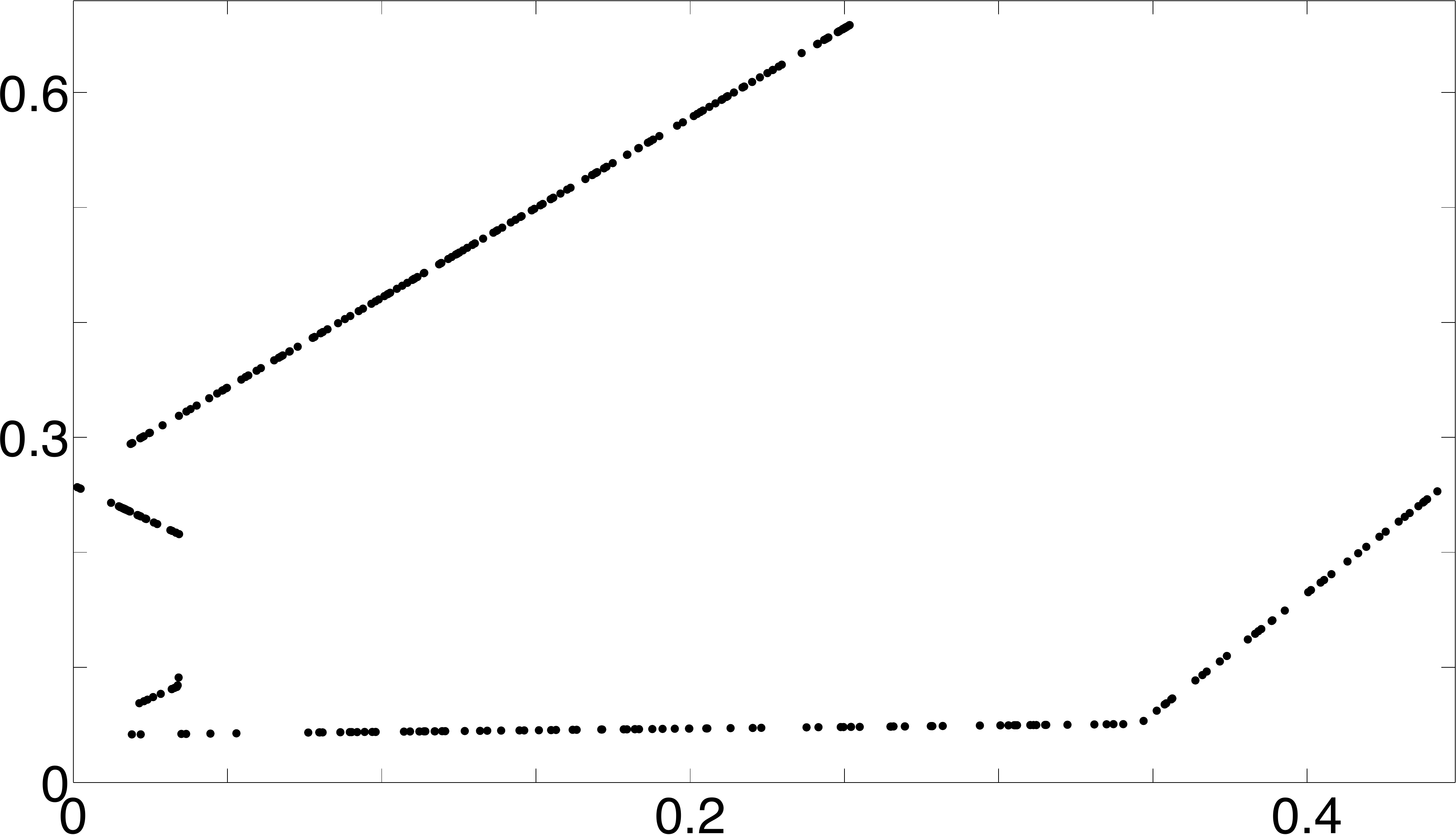}
\put(-15,120){\footnotesize{\textit{(b)}}}
\put(-86.8,1.5){\footnotesize{$x(t-\tau)$}}
\put(-228,71.2){\rotatebox{90}{\footnotesize{$x(t-2\tau)$}}}
\vspace*{-02mm}
\caption{Projected Poincar\'{e} section of the orbit of Figure~\ref{fig_PS}(e)-(f) onto the plane $(x(t-\tau),x(t-2\tau))$ for crossing of the Poincar\'{e} section $x(t)=0.14$ with $x^{\prime}(t)>0$ for panel (a) and $x^{\prime}(t)<0$ for panel (b).}
\label{fig_poincare}
\end{figure}
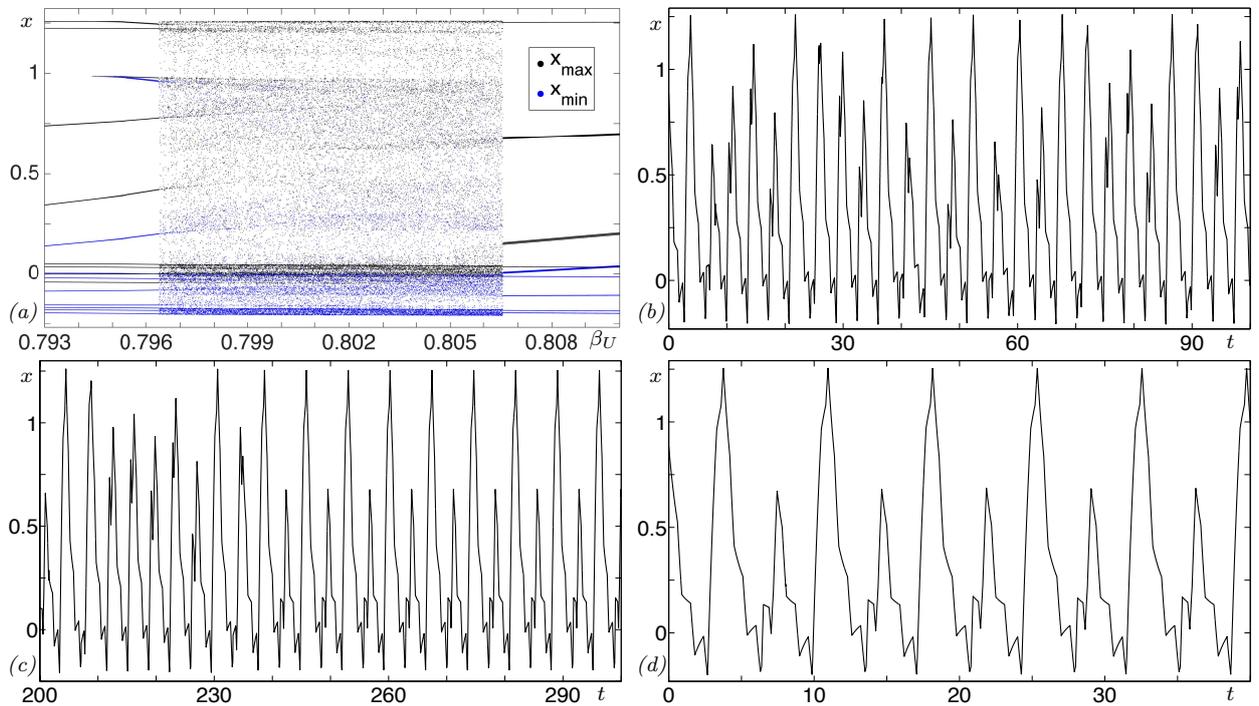
\begin{figure}[!t]
\includegraphics[width=0.49\textwidth,height=0.276\textwidth]{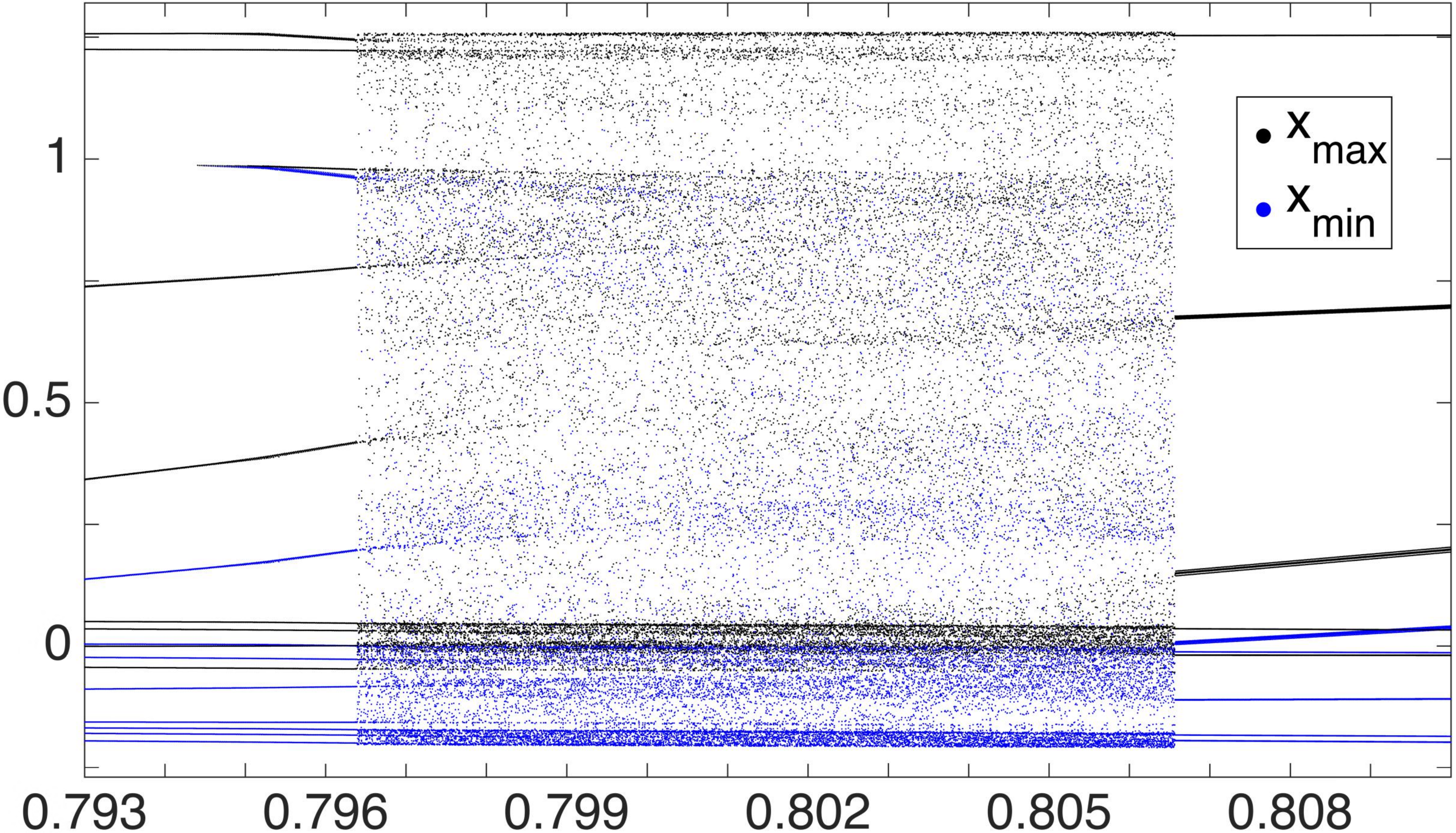}
\put(-230,12){\footnotesize{\textit{(a)}}}
\put(-12,2.1){\footnotesize{$\beta_{U}$}}
\put(-225,122){\footnotesize{$x$}}
\hspace*{1mm}
\includegraphics[width=0.49\textwidth,height=0.276\textwidth]{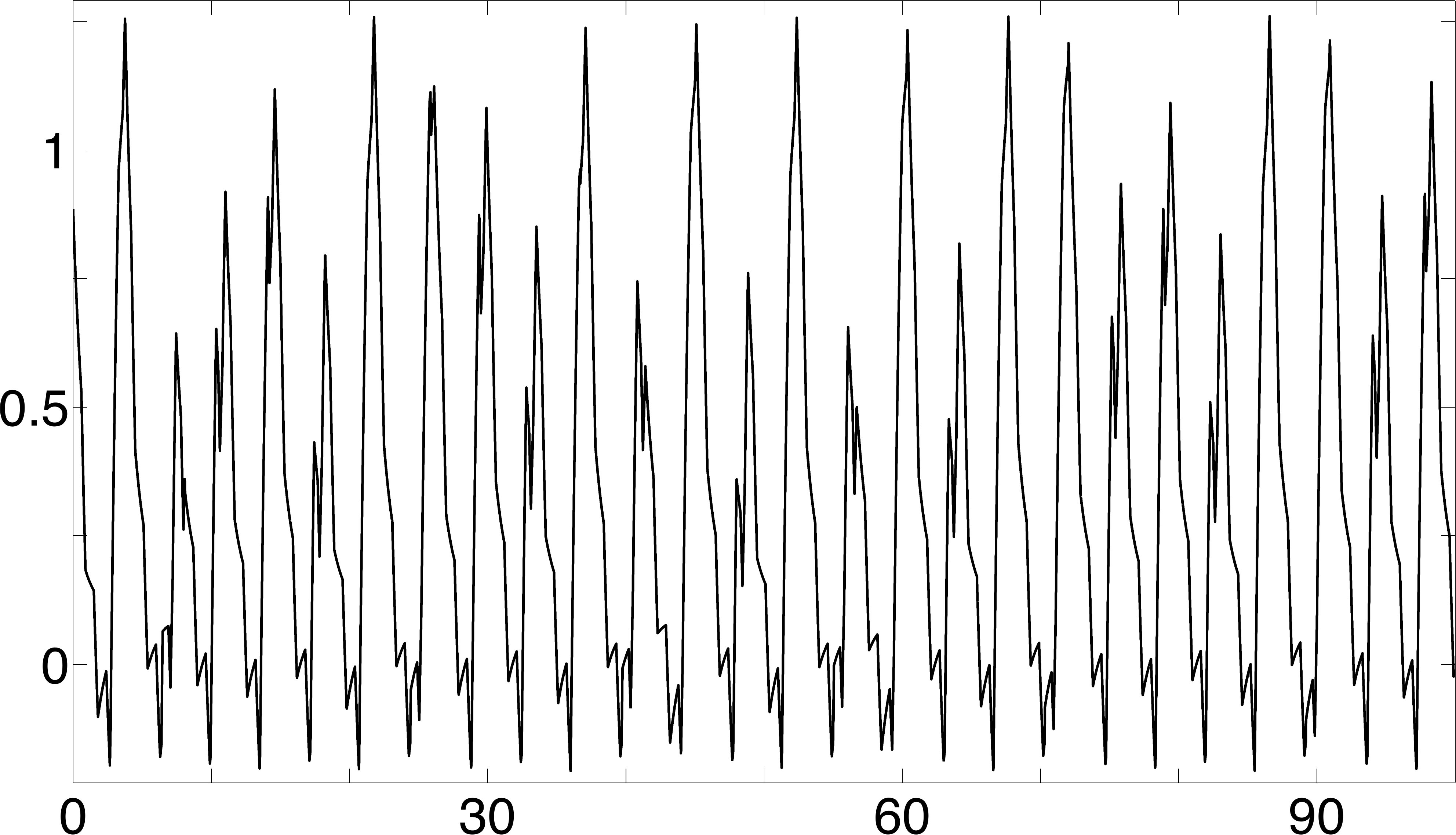}
\put(-230,12){\footnotesize{\textit{(b)}}}
\put(-9,1){\footnotesize{$t$}}
\put(-225,121){\footnotesize{$x$}}
\vspace*{1mm}\\
\includegraphics[width=0.49\textwidth,height=0.276\textwidth]{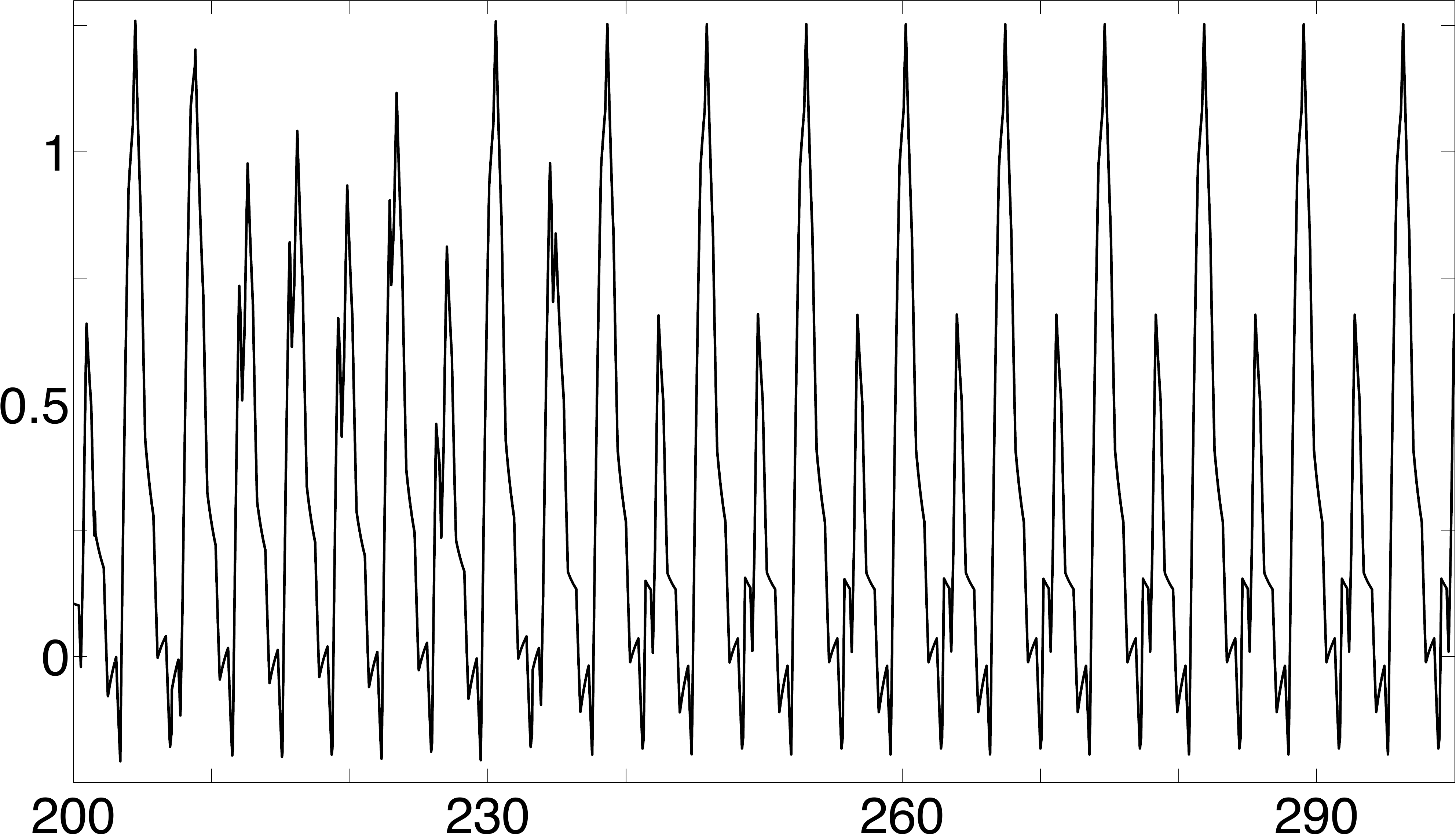}
\put(-230,12){\footnotesize{\textit{(c)}}}
\put(-9,1){\footnotesize{$t$}}
\put(-225,121){\footnotesize{$x$}}
\hspace*{1mm}
\includegraphics[width=0.49\textwidth,height=0.276\textwidth]{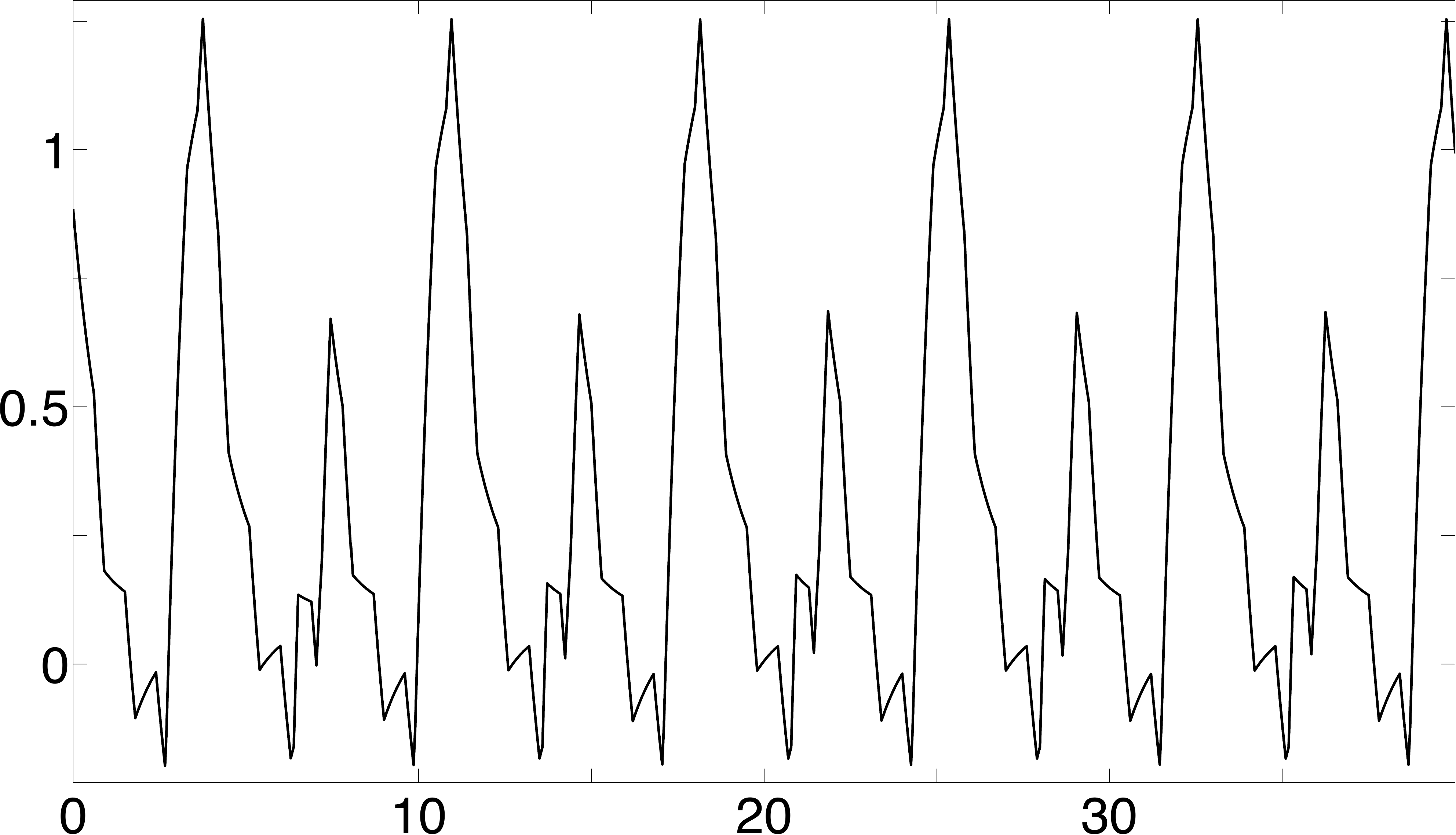}
\put(-230,12){\footnotesize{\textit{(d)}}}
\put(-9,1){\footnotesize{$t$}}
\put(-225,121){\footnotesize{$x$}}
\vspace*{-02mm}
\caption{(a) A magnification of part of the diagram of Figure~\ref{fig_Bif}(c). For the same parameters, the other three panels show solutions of~\eqref{eq.diff.pert} transitioning from irregular to regular motion with $\beta_{U}=0.805$ for panel (b), $\beta_{U}=0.807$ for panel (c) and $\beta_{U}=0.808$ for panel (d).}
\label{fig_trans}
\end{figure}
In order to solve~\eqref{eq.diff.pert} it is necessary to define a history function for $t\in[-\tau,0]$, which is an infinite-dimensional set of initial conditions. Thereby the solution space of~\eqref{eq.diff.pert} is infinite dimensional and consequently a hyperplane defined by a Poincar\'{e} section is also infinite dimensional. Although the Poincar\'{e} section is infinite dimensional, we can project it onto $\mathbb{R}^2$ by taking, for example, the solution points $x(t-\tau)$ and $x(t-2\tau)$ such that $x(t)=k$ for some constant $k\in\mathbb{R}$. In Figure~\ref{fig_poincare} we project a Poincar\'{e} section of the orbit of Figure~\ref{fig_PS}(e)-(f) onto the plane $(x(t-\tau),x(t-2\tau))$ for crossing of the Poincar\'{e} section $x(t)=0.1$ with $x^{\prime}(t)>0$ for panel (a) and $x^{\prime}(t)<0$ for panel (b). Both Poincar\'{e} sections form a locus with sparse points, indicating that the irregular motion of the orbit of Figure~\ref{fig_PS}(e)-(f) is chaotic~\citep{Nayfeh_2007}.

Figure~\ref{fig_trans}(a) shows a magnified part of the diagram of Figure~\ref{fig_Bif}(c). In the other panels we show orbits from the diagram of panel (a) transitioning from a region of regular motion to regular motion. In panel (a) we have an irregular motion with $\beta_{U}=0.805$, in panel (b) we have $\beta_{U}=0.808$ and the orbit changes from irregular to regular motion for $t\approx 240$, while in panel (c) the transition to regular motion occurs for $t\approx 2$. This scenario of an abrupt transition from periodic motion to irregular motion and vice versa is highly suggestive of  a \textit{boundary or an exterior crisis} scenario~\citep{Nayfeh_2007,Grebogi_1983}.
%
\section{Discussion and Conclusions}\label{sec.conc}
Here we have been able to exploit the relative simplicity (due to the piecewise linear nonlinearity) of our model system in an extension of the results of~\citet{Mackey_2017} examining the response to a single stimulus (Section~\ref{sec.single}), as well as examining the response of the system to a periodic stimulus in Section~\ref{sec.per}.  The insights and techniques of Section~\ref{sec.per} allowed us to draw conclusions about treatment implications in Section~\ref{sec:treatment}.  The numerical bifurcation diagrams of Figure~\ref{fig_Bif} revealed that an effective way to increase the minima of the oscillations, and hence decrease the neutropenia duration, is to increase the drug dosage by either increasing the time duration of the drug administration $\sigma$, or increasing  the drug dose administrated per unit of time, which is proportional to $a$. From Propositions~\ref{prop.LC1} and~\ref{prop.LC2} and Eq.~\eqref{eq.alpha} we obtained the condition  $a=a_{1}$ to avoid neutropenia, from which we computed the minimal interval between drug administrations $\alpha$ to avoid neutropenia as function of the time duration of the drug administration $\sigma$ by Eq.~\eqref{eq.alpha}.

The numerical results we have presented in Section~\ref{sec.num.exp.per.pert}, while not exhaustive, certainly indicate that there is a wealth of bifurcation behaviour to be understood in this relatively simple system.  However, this must remain the object of further study as it is outside the main thrust of this paper.

Finally we note that the study of periodic perturbation of limit cycle systems has been long and intensive, particularly in a biological context, c.f.~\citet{Winfree_1980,guevara1982,glass-winfree-1984,trine2004,bodnar2013b}, and the emphasis has been on an examination of the phase response curve.  However, to apply numerical methods such as the phase reduction method due to~\citet{Novicenko_2012} and~\citet{Kotani_2012} to calculate the phase response (or phase resetting) curve for a DDE and the method to compute the approximating Lyapunov exponents for DDEs due to~\citet{Breda_2014} we need to linearize the DDE~\eqref{eq.diff.pert} around a reference orbit, but the feedback function $f(x(t-\tau))$ is discontinuous at $x=0$.  An extended approach using~\eqref{eq.diff.pert} with a continuous delayed feedback $f(x(t-\tau))$ (a sigmoid type) would allow us to apply these numerical methods to study the solutions with $f(x(t-\tau))$ approaching the discontinuous delayed feedback.  This too we reserve for a future study.
%
%
\section*{Acknowledgments}
DCS was supported by National Council for Scientific and Technological Development of Brazil (CNPq) postdoctoral fellowship 201105/2014-4, and MCM is supported by a Discovery Grant from the Natural Sciences and Engineering Research Council (NSERC) of Canada.  
MCM would like to thank the Institut f\"{u}r Theoretische Neurophysik, Universit\"{a}t Bremen for their hospitality during the time in which much of the writing of this paper took place.
We are very grateful to Tony Humphries for fruitful discussions and suggestions.
%
\begin{appendices}
\section{Proof of the Results}
\label{sec.app}
Here we present the proofs of Remarks and Propositions from Sections~\ref{sec.single} and~\ref{sec.per}.
\begin{proof}[{\bf Proof of Remark~\ref{rem.FNFP}}]
First it is shown that $a>\beta_{U}$ holds. For $\Delta\leq t \leq \Delta+\sigma$ we have
\begin{equation}\label{FNFPa}
x^{(\Delta)}(t)  = -\beta_{U} + a + (x^{(\Delta)}(\Delta)+\beta_{U}-a)\mathrm{e}^{-(t-\Delta)},
\end{equation}
\noindent with $x^{(\Delta)}(\Delta) = -\beta_{U} + \beta_{U}\mathrm{e}^{\tilde{z}_{2}-\Delta}$ and
\begin{eqnarray}
x^{(\Delta)}(\Delta+\sigma) &=& -\beta_{U} + a + (x^{(\Delta)}(\Delta)+\beta_{U}-a)\mathrm{e}^{-\sigma} \label{eq.xDeltaSigma}\\
&=& \beta_{U}(\mathrm{e}^{\tilde{z}_{2}-\Delta-\sigma}-1) + a(1-\mathrm{e}^{-\sigma}). \label{eq.xDeltaSigmaB}
\end{eqnarray}
From~\eqref{eq.xDeltaSigma} the condition $x^{(\Delta)}(\Delta+\sigma)> 0$ can be written as
\begin{equation}\nonumber
x^{(\Delta)}(\Delta+\sigma)=(a-\beta_{U})(1-\mathrm{e}^{-\sigma}) + x^{(\Delta)}(\Delta)\mathrm{e}^{-\sigma},
\end{equation}
but $x^{(\Delta)}(\Delta)=-\beta_{U}(1-\mathrm{e}^{\tilde{z}_{2}-\Delta})<0$ since $\Delta>\tilde{z}_{2}$, then we must have $a>\beta_{U}$.

For $\Delta+\sigma\leq t \leq \tilde{T}$ we have
\begin{equation}\label{FNFPb}
x^{(\Delta)}(t)  = -\beta_{U} + (x^{(\Delta)}(\Delta+\sigma)+\beta_{U})\mathrm{e}^{-(t-\Delta-\sigma)}.
\end{equation}

Equation~\eqref{eq.xDeltaSigmaB} gives $x^{(\Delta)}(\Delta+\sigma) + \beta_{U} = \beta_{U}\mathrm{e}^{\tilde{z}_{2}-\Delta-\sigma} + a(1-\mathrm{e}^{-\sigma})$ and this combined with $\tilde{T}=\tilde{z}_{2}+\tau$ from~\eqref{z1z2T} in the solution~\eqref{FNFPb} computed at $t=\tilde{T}$ yields
\begin{eqnarray}
x^{(\Delta)}(\tilde{T})  &=& -\beta_{U} + (\beta_{U}\mathrm{e}^{\tilde{z}_{2}-\Delta-\sigma} + a(1-\mathrm{e}^{-\sigma}))\mathrm{e}^{-\tilde{T}+\Delta+\sigma} \notag\\
&=& -\beta_{U}(1-\mathrm{e}^{-\tau}) + a(\mathrm{e}^{\sigma}-1)\mathrm{e}^{\Delta-\tilde{T}}.\label{eq.xDeltaTtilde}
\end{eqnarray}
Hence the inequality $x^{(\Delta)}(\tilde{T})<0$ holds if and only if
\begin{equation}\nonumber
a(\mathrm{e}^{\sigma}-1)\mathrm{e}^{\Delta} < \beta_{U}(\mathrm{e}^{\tau}-1)\mathrm{e}^{\tilde{z}_{2}}.
\end{equation}
\noindent Defining $\delta_{4}$ as the constant such that $x^{(\Delta)}(\tilde{T})=0$ for $\Delta=\delta_{4}$ leads to~\eqref{delta4}.

Then $x^{(\Delta)}(\tilde{T})<0$ if, and only if, $\Delta<\delta_{4}$ while $x^{(\Delta)}(\tilde{T})\geq0$ if, and only if, $\Delta\geq\delta_{4}$.
The interval $\mathit{I_{FNFP}}=(\tilde{z}_{2},\tilde{T}-\sigma)\cap (-\infty,\delta_{2})$ combined with the condition $x^{(\Delta)}(\tilde{T})\geq 0$ (see the examples in Figure~\ref{fig_FNFP2}) gives the interval $\mathit{I_{FNFP1}}$ given by~\eqref{FNFP1}, while $\mathit{I_{FNFP}}$ combined with the condition $x^{(\Delta)}(\tilde{T})<0$ (see the examples in Figure~\ref{fig_FNFP1}) yields  the interval $\mathit{I_{FNFP2}}$ given by~\eqref{FNFP2}.
\end{proof}
%
\begin{proof}[{\bf Proof of Remark~\ref{rem.FNFPB}}]
\vspace{5mm}
\noindent\textbf{FNFP3:} for this subcase we have to assume $x^{(\Delta)}(z_{\Delta,3}+\tau)< 0$ (see the examples in Figure~\ref{fig_FNFP24}).

From $x^{(\Delta)}(\Delta)<0$ and $x^{(\Delta)}(\Delta+\sigma)>0$ we obtain a zero $z_{\Delta,3}$ of $x^{(\Delta)}$ in $(\Delta,\Delta+\sigma)$ given by
\begin{equation}\nonumber
0=x^{(\Delta)}(z_{\Delta,3})=-\beta_{U} + a + (x^{(\Delta)}(\Delta)+\beta_{U}-a)\mathrm{e}^{-(z_{\Delta,3}-\Delta)},
\end{equation}
\begin{equation}\label{eq.zDelta3}
\mathrm{e}^{z_{\Delta,3}}= \frac{(a-\beta_{U}-x^{(\Delta)}(\Delta))\mathrm{e}^{\Delta}}{a-\beta_{U}} = \frac{a\mathrm{e}^{\Delta}-\beta_{U}\mathrm{e}^{\tilde{z}_{2}}}{a-\beta_{U}}.
\end{equation}
From~\eqref{FNFPb} and $(x^{(\Delta)}(\Delta+\sigma)+\beta_{U})>0$ it follows that $x^{(\Delta)}(t)$ is strictly decreasing for $t\in[\Delta+\sigma,\tilde{T}]$.
This together with $x^{(\Delta)}(\Delta+\sigma)>0$ and $x^{(\Delta)}(\tilde{T})<0$ yields that there is a zero $z_{\Delta,4}$ of $x^{(\Delta)}$ in $(\Delta+\sigma,\tilde{T})$ given by
\begin{equation}\nonumber
0=x^{(\Delta)}(z_{\Delta,4})=-\beta_{U} + (x^{(\Delta)}(\Delta+\sigma)+\beta_{U})\mathrm{e}^{-(z_{\Delta,4}-\Delta-\sigma)},
\end{equation}
\begin{equation}\label{eq.zDelta4}
\mathrm{e}^{z_{\Delta,4}}= \frac{(x^{(\Delta)}(\Delta+\sigma)+\beta_{U})\mathrm{e}^{\Delta+\sigma}}{\beta_{U}} = \frac{\beta_{U}\mathrm{e}^{\tilde{z}_{2}} + a(\mathrm{e}^{\sigma}-1)\mathrm{e}^{\Delta}}{\beta_{U}}.
\end{equation}

For $t\in(\tilde{T},z_{\Delta,3}+\tau)$ we have $\tilde{z}_{2}<t-\tau<z_{\Delta,3}$. Hence $x^{(\Delta)}(t-\tau)<0$ and
\begin{equation}\label{zDelta3tauB}
x^{(\Delta)}(t)  = \beta_{L} + (x^{(\Delta)}(\tilde{T})-\beta_{L})\mathrm{e}^{-(t-\tilde{T})}.
\end{equation}
From~\eqref{zDelta3tauB} and $(x^{(\Delta)}(\tilde{T})-\beta_{L})<0$ we obtain that $x^{(\Delta)}(t)$ is strictly decreasing for $t\in(\tilde{T},z_{\Delta,3}+\tau)$. At $t=(z_{\Delta,3}+\tau)$ Eq.~\eqref{zDelta3tauB} gives
\begin{equation}\label{xDelta3tau}
x^{(\Delta)}(z_{\Delta,3}+\tau)  = \beta_{L} + (x^{(\Delta)}(\tilde{T})-\beta_{L})\mathrm{e}^{-(z_{\Delta,3}-\tilde{z}_{2})},
\end{equation}
which is negative by assumption.

The condition $x^{(\Delta)}(z_{\Delta,3}+\tau)< 0$ together with Eqs.~\eqref{xDelta3tau},~\eqref{eq.zDelta3},~\eqref{eq.xDeltaTtilde} and $a>\beta_{U}$ yield
\begin{eqnarray}
\beta_{L}\mathrm{e}^{z_{\Delta,3}-\tilde{z}_{2}} &<& \beta_{L}-x^{(\Delta)}(\tilde{T}), \notag\\
\beta_{L}(a\mathrm{e}^{\Delta-\tilde{z}_{2}}-\beta_{U}) &<& (a-\beta_{U})[\beta_{L} + \beta_{U}(1-\mathrm{e}^{-\tau}) - a(\mathrm{e}^{\sigma}-1)\mathrm{e}^{\Delta-\tilde{T}}], \notag\\
\Delta &<& \ln{\frac{a\beta_{L}+\beta_{U}(a-\beta_{U})(1-\mathrm{e}^{-\tau})}{a[\beta_{L}\mathrm{e}^{-\tilde{z}_{2}}+(a-\beta_{U})(\mathrm{e}^{\sigma}-1)\mathrm{e}^{-\tilde{T}}]}}\eqqcolon \hat{\delta}_{4}.\label{eq.delta4hat}
\end{eqnarray}

Thus the interval $\mathit{I_{FNFP2}}$, given by~\eqref{FNFP2}, together with the extra condition $x^{(\Delta)}(z_{\Delta,3}+\tau)<0$, written as~\eqref{eq.delta4hat}, yield the interval $\mathit{I_{FNFP3}}$ given by~\eqref{FNFP3}.

\vspace{1mm}
\noindent\textbf{FNFP4:} in this subcase we also have to assume $x^{(\Delta)}(z_{\Delta,4}+\tau)<0$ (see the examples in Figure~\ref{fig_FNFP1}).

For $t\in[z_{\Delta,3}+\tau,z_{\Delta,4}+\tau]$ we have $z_{\Delta,3}\leq t-\tau\leq z_{\Delta,4}$. Hence $x^{(\Delta)}(t-\tau)\geq 0$ and
\begin{equation}\nonumber
x^{(\Delta)}(t)  = -\beta_{U} + (x^{(\Delta)}(z_{\Delta,3}+\tau)+\beta_{U})\mathrm{e}^{-(t-z_{\Delta,3}-\tau)}.
\end{equation}

From the condition $x^{(\Delta)}(z_{\Delta,3}+\tau)< 0$ we have $x^{(\Delta)}(z_{\Delta,3}+\tau)+\beta_{U}>0$, so $x^{(\Delta)}$ is strictly decreasing on  $[z_{\Delta,3}+\tau,z_{\Delta,4}+\tau]$, and
\begin{equation}\label{eq.xDelta4tau}
x^{(\Delta)}(z_{\Delta,4}+\tau)  = -\beta_{U} + (x^{(\Delta)}(z_{\Delta,3}+\tau)+\beta_{U})\mathrm{e}^{-(z_{\Delta,4}-z_{\Delta,3})}.
\end{equation}

The condition $x^{(\Delta)}(z_{\Delta,4}+\tau)<0$ along with Eqs.~\eqref{eq.xDelta4tau},~\eqref{eq.zDelta4} and~\eqref{xDelta3tau} yield
\begin{eqnarray}
\beta_{U}\mathrm{e}^{z_{\Delta,4}}  &>& (x^{(\Delta)}(z_{\Delta,3}+\tau)+\beta_{U})\mathrm{e}^{z_{\Delta,3}}, \notag\\
\beta_{U}\mathrm{e}^{z_{\Delta,4}}  &>& (\beta_{U} + \beta_{L})\mathrm{e}^{z_{\Delta,3}} - (\beta_{L}-x^{(\Delta)}(\tilde{T}))\mathrm{e}^{\tilde{z}_{2}}, \notag
\end{eqnarray}
\vspace{-5.0mm}
\begin{equation}\nonumber
\left[\beta_{U}\frac{\beta_{U}+\beta_{L}}{a-\beta_{U}} + \beta_{L}+\beta_{U}(2-\mathrm{e}^{-\tau})\right]\mathrm{e}^{\tilde{z}_{2}} > a\left[\frac{\beta_{U} + \beta_{L}}{a-\beta_{U}} - (\mathrm{e}^{\sigma}-1)(1-\mathrm{e}^{-\tau})\right]\mathrm{e}^{\Delta},
\end{equation}
\vspace{-1.0mm}
\begin{equation}\label{eq.delta5}
\Delta < \tilde{z}_{2} + \ln{\frac{\beta_{U}(\beta_{U}+\beta_{L})+\left[\beta_{L}+\beta_{U}(2-\mathrm{e}^{-\tau})\right](a-\beta_{U})}{a(\beta_{U}+\beta_{L})-a(\mathrm{e}^{\sigma}-1)(1-\mathrm{e}^{-\tau})(a-\beta_{U})}} \eqqcolon \delta_{5}.
\end{equation}

From~\eqref{eq.delta4hat} we infer that the condition $x^{(\Delta)}(z_{\Delta,3}+\tau)\geq0$ implies $\Delta\geq\hat{\delta}_{4}$. This condition together with the interval $\mathit{I_{FNFP2}}$, given by~\eqref{FNFP2}, plus the extra condition $x^{(\Delta)}(z_{\Delta,4}+\tau)<0$, written as~\eqref{eq.delta5}, yield the interval $\mathit{I_{FNFP4}}$ given by~\eqref{FNFP4}.
\end{proof}
%
\begin{proof}[{\bf Proof of Remark~\ref{rem.calcF}}]
For each case~\eqref{cases} we compute the resetting time $F(\Delta)$ as follows:

\noindent\textbf{RNRN:} $F(\Delta)=\sigma$, since $x^{(\Delta)}(t)=\tilde{x}(t+(\tilde{z}_{1}-z_{\Delta,1}))$ for all $t\geq \Delta+\sigma$ and $x^{(\Delta)}(z_{\Delta,1}+\tau)=\tilde{x}(\tilde{z}_{1}+\tau)=\bar{x}$, where  $z_{\Delta,1}=\tilde{z}_{1}+T(\Delta)-\tilde{T}$ and $T(\Delta)$ is as in~\citet[Proposition 5.1]{Mackey_2017};

\noindent\textbf{RNRP:} $F(\Delta)=\tilde{z}_{1}+\tau+(z_{\Delta,2}-\tilde{z}_{2})-\Delta$, since $x^{(\Delta)}(t)=\tilde{x}(t-(z_{\Delta,2}-\tilde{z}_{2}))$ for all $t\geq \tilde{z}_{1}+\tau+(z_{\Delta,2}-\tilde{z}_{2})$ and $x^{(\Delta)}(z_{\Delta,2}+\tau)=\tilde{x}(\tilde{z}_{2}+\tau)=\ubar{x}$, where $z_{\Delta,2}=\tilde{z}_{2}+T(\Delta)-\tilde{T}$ and $T(\Delta)$ is as in~\citet[Proposition 5.2]{Mackey_2017};

\noindent\textbf{RPRP:} $F(\Delta)=\tilde{z}_{1}+\tau+(z_{\Delta,2}-\tilde{z}_{2})-\Delta$, since $x^{(\Delta)}(t)=\tilde{x}(t-(z_{\Delta,2}-\tilde{z}_{2}))$ for all $t\geq \tilde{z}_{1}+\tau+(z_{\Delta,2}-\tilde{z}_{2})$ and $x^{(\Delta)}(z_{\Delta,2}+\tau)=\tilde{x}(\tilde{z}_{2}+\tau)=\ubar{x}$, where $z_{\Delta,2}=\tilde{z}_{2}+T(\Delta)-\tilde{T}$ and $T(\Delta)$ is as in~\citet[Proposition 5.3]{Mackey_2017};

\noindent\textbf{RPFP:} $F(\Delta)=\sigma$, since $x^{(\Delta)}(t)=\tilde{x}(t-(z_{\Delta,2}-\tilde{z}_{2}))$ for all $t\geq \Delta+\sigma$ and $x^{(\Delta)}(z_{\Delta,2}+\tau)=\tilde{x}(\tilde{z}_{2}+\tau)=\ubar{x}$, where $z_{\Delta,2}=\tilde{z}_{2}+T(\Delta)-\tilde{T}$ and $T(\Delta)$ is as in~\citet[Proposition 5.4]{Mackey_2017};

\noindent\textbf{RPFN:} $F(\Delta)=z_{\Delta,2}+\tau-\Delta$, since $x^{(\Delta)}(t)=\tilde{x}(t-(z_{\Delta,3}-\tilde{z}_{3}))$ for all $t\geq z_{\Delta,2}+\tau$ and $x^{(\Delta)}(z_{\Delta,3}+\tau)=\tilde{x}(\tilde{z}_{3}+\tau)=\bar{x}$, where $z_{\Delta,3}=\tilde{z}_{3}+T(\Delta)-\tilde{T}$ and $T(\Delta)$ is as in~\citet[Proposition 5.5]{Mackey_2017};

\noindent\textbf{FPFP:} $F(\Delta)=\sigma$, since $x^{(\Delta)}(t)=\tilde{x}(t-(z_{\Delta,2}-\tilde{z}_{2}))$ for all $t\geq \Delta+\sigma$ and $x^{(\Delta)}(z_{\Delta,2}+\tau)=\tilde{x}(\tilde{z}_{2}+\tau)=\ubar{x}$, where $z_{\Delta,2}=\tilde{z}_{2}+T(\Delta)-\tilde{T}$ and $T(\Delta)$ is as in~\citet[Proposition 5.6]{Mackey_2017};

\noindent\textbf{FPFN:} $F(\Delta)=z_{\Delta,2}+\tau-\Delta$, since $x^{(\Delta)}(t)=\tilde{x}(t-(z_{\Delta,3}-\tilde{z}_{3}))$ for all $t\geq z_{\Delta,2}+\tau$ and $x^{(\Delta)}(z_{\Delta,3}+\tau)=\tilde{x}(\tilde{z}_{3}+\tau)=\bar{x}$, where $z_{\Delta,2}=\tilde{z}_{2}+T(\Delta)-\tilde{T}$ and $z_{\Delta,3}=\tilde{z}_{3}+T(\Delta)-\tilde{T}$ with $T(\Delta)$ and $z_{\Delta,2}$ given by~\citet[Proposition 5.7 and its proof]{Mackey_2017};

\noindent\textbf{FNFP1:} $F(\Delta)=z_{\Delta,3}+\tau-\Delta$, since $x^{(\Delta)}(t)=$$\tilde{x}(t-(z_{\Delta,4}-\tilde{z}_{4}))$ for all $t\geq z_{\Delta,3}+\tau$ and $x^{(\Delta)}(z_{\Delta,4}+\tau)=\tilde{x}(\tilde{z}_{4}+\tau)=\ubar{x}$, where $z_{\Delta,3}$ and  $z_{\Delta,4}=\tilde{z}_{4}$ are computed as follows.
From  $x^{(\Delta)}(\Delta)<0<x^{(\Delta)}(\Delta+\sigma)$ we obtain a zero $z_{\Delta,3}$ of $x^{(\Delta)}$ in $(\Delta,\Delta+\sigma)$ given by
\begin{equation}\nonumber
x^{(\Delta)}(z_{\Delta,3})=-\beta_{U}+a+(x^{(\Delta)}(\Delta)+\beta_{U}-a)\mathrm{e}^{-(z_{\Delta,3}-\Delta)}=0,
\end{equation}
where $x^{(\Delta)}(\Delta)=-\beta_{U}+\beta_{U}\mathrm{e}^{-(\tilde{z}_{2}-\Delta)}$. For $\tilde{T}<t<z_{\Delta,3}+\tau$ we have $\tilde{z}_{2}<t-\tau<z_{\Delta,3}$. Hence $x^{(\Delta)}(t-\tau)<0$, and
\begin{equation}\nonumber
x^{(\Delta)}(t)=\beta_{L}+(x^{(\Delta)}(\tilde{T})-\beta_{L})\mathrm{e}^{-(t-\tilde{T})}.
\end{equation}
Since  $\tilde{z}_{2}<z_{\Delta,3}$ and from the proof of Remark~\ref{rem.FNFP} the inequality $\Delta\geq\delta_{4}$ implies $x^{(\Delta)}(\tilde{T})\geq0$, we obtain
\begin{equation}\label{zDelta3tau}
x^{(\Delta)}(z_{\Delta,3}+\tau)=\beta_{L}+(x^{(\Delta)}(\tilde{T})-\beta_{L})\mathrm{e}^{-(z_{\Delta,3}+\tau-\tilde{T})}>\beta_{L}(1-\mathrm{e}^{\tilde{z}_{2}-z_{\Delta,3}})>0.
\end{equation}
Notice that $x^{(\Delta)}(t)>0$ on $(z_{\Delta,3},z_{\Delta,3}+\tau]$. Using this and $z_{\Delta,3}<\Delta+\sigma<\tilde{z}_{2}+\tau<z_{\Delta,3}+\tau$ we obtain that $x^{(\Delta)}(t)$ is strictly decreasing on $[z_{\Delta,3}+\tau,\infty)$ as long as $x^{(\Delta)}(t-\tau)\geq 0$. It follows that there is a first zero $z_{\Delta,4}$ of $x^{(\Delta)}(t)$ in $[z_{\Delta,3}+\tau,\infty)$ given by
\begin{equation}\nonumber
x^{(\Delta)}(z_{\Delta,4})=-\beta_{U}+(x^{(\Delta)}(z_{\Delta,3}+\tau)-\beta_{U})\mathrm{e}^{-(z_{\Delta,4}-(z_{\Delta,3}+\tau))}=0,
\end{equation}
with $x^{(\Delta)}(z_{\Delta,3}+\tau)$ given by~\eqref{zDelta3tau}, where $x^{(\Delta)}(\tilde{T})$ is given by~\eqref{eq.xDeltaTtilde};

\noindent\textbf{FNFN:} $F(\Delta)=\tilde{z}_{2}+\tau-\Delta$, since $x^{(\Delta)}(t)=\tilde{x}(t+(\tilde{z}_{3}-z_{\Delta,3}))$ for all $t\geq \tilde{z}_{2}+\tau$ and $x^{(\Delta)}(z_{\Delta,3}+\tau)=\tilde{x}(\tilde{z}_{3}+\tau)=\bar{x}$, where $z_{\Delta,3}=\tilde{z}_{3}+T(\Delta)-\tilde{T}$ and $T(\Delta)$ is as in~\citet[Proposition 5.8]{Mackey_2017};

\noindent\textbf{FNRN:} $F(\Delta)=\sigma$, since $x^{(\Delta)}(t)=\tilde{x}(t+(\tilde{z}_{3}-z_{\Delta,3}))$ for all $t\geq \Delta+\sigma$ and $x^{(\Delta)}(z_{\Delta,3}+\tau)=\tilde{x}(\tilde{z}_{3}+\tau)=\bar{x}$, where $z_{\Delta,3}=\tilde{z}_{3}+T(\Delta)-\tilde{T}$ and $T(\Delta)$ is as in~\citet[Proposition 5.10]{Mackey_2017};

\noindent\textbf{FNRP:} $F(\Delta)=\tilde{z}_{3}+\tau-(\tilde{z}_{4}-z_{\Delta,4})-\Delta$, since $x^{(\Delta)}(t)=\tilde{x}(t+(\tilde{z}_{4}-z_{\Delta,4}))$ for all $t\geq \tilde{z}_{3}+\tau-(\tilde{z}_{4}-z_{\Delta,4})$ and $x^{(\Delta)}(z_{\Delta,4}+\tau)=\tilde{x}(\tilde{z}_{4}+\tau)=\ubar{x}$, where $z_{\Delta,4}=\tilde{z}_{4}+T(\Delta)-\tilde{T}$ and $T(\Delta)$ is as in~\citet[Proposition 5.10]{Mackey_2017}.
\end{proof}
\begin{proof}[{\bf Proof of Remark~\ref{rem.FT}}]
Define the constants $\delta_{1}$ as in~\citet[Eq. (5.6)]{Mackey_2017}
\begin{equation}\label{delta1hat1}
\delta_{1} \coloneqq \tilde{z}_{1} - \sigma - \ln\left(\frac{\beta_{L}+a(1-\mathrm{e}^{-\sigma})}{\beta_{L}}\right)<\tilde{z}_{1}. 
\end{equation}
For each case~\eqref{cases} we consider the corresponding $\Delta$ interval as computed in~\citet[Section 5]{Mackey_2017} and listed in~\eqref{cases}. Recalling that $\sigma\in(0,\tau]$, $\tilde{T}=\tilde{z}_{2}+\tau$, $t_{max}=\tilde{z}_{1}+\tau$, $\tilde{z}_{n}=\tilde{T}+\tilde{z}_{n-2}$ for $n\in\{2,3,4\ldots\}$,  $\tilde{z}_{j+1}>\tilde{z}_{j}+\tau$ for all $j\in\mathbb{N}$, and from~\citet[Proposition 4.2]{Mackey_2017} we see that $J=j_\Delta \in \{0,1,2\}$ implies $z_{\Delta,J+1}>\Delta$ and $z_{\Delta,J+1}>z_{\Delta,J+1}+\tau$. Thus, we show that for each case~\eqref{cases} we have $F(\Delta)<T(\Delta)$ as follows: 

\noindent\textbf{RNRN:} $\Delta\in\mathit{I_{RNRN}}=[0,\delta_{1}]$: since $\tilde{z}_{2}>\tilde{z}_{1}+\tau>\tilde{z}_{1}$,
\begin{equation}\nonumber
T(\Delta)=\tilde{T}+z_{\Delta,1}-\tilde{z}_{1}=\tilde{z}_{2}+\tau+z_{\Delta,1}-\tilde{z}_{1}>\tau+z_{\Delta,1}>\tau\geq\sigma=F(\Delta);
\end{equation}

\noindent\textbf{RNRP and RPRP:} $\Delta\in\mathit{I_{RNRP}}=[\max\{0,\delta_{1}\},\tilde{z}_{1})$ with $\delta_{1}$ given by~\eqref{delta1hat1} and $\Delta\in\mathit{I_{RPRP}}=[\tilde{z}_{1},t_{max}-\sigma]$: using that $\tilde{z}_{2}>\tilde{z}_{1}+\tau>\tilde{z}_{1}$ we infer
\begin{equation}\nonumber
T(\Delta)=\tilde{T}+z_{\Delta,2}-\tilde{z}_{2}=\tau+z_{\Delta,2}>\tilde{z}_{1}+\tau+z_{\Delta,2}-\tilde{z}_{2}-\Delta=F(\Delta);
\end{equation}

\noindent\textbf{RPFP:} $\Delta\in\mathit{I_{RPFP}}=(t_{\max}-\sigma,t_{\max}]\cap(-\infty,\delta_{2}]$:
\begin{equation}\nonumber
T(\Delta)=\tilde{T}+z_{\Delta,2}-\tilde{z}_{2}=\tau+z_{\Delta,2}>\tau\geq\sigma=F(\Delta);
\end{equation}

\noindent\textbf{RPFN:} $\Delta\in\mathit{I_{RPFN}}=(t_{\max}-\sigma,t_{\max}]\cap(\delta_{2},\infty)$: from~\citet[Proposition 4.2]{Mackey_2017} we have $z_{\Delta,3}>z_{\Delta,2}+\tau$, and thereby
\begin{equation}\nonumber
T(\Delta)=\tilde{T}+z_{\Delta,3}-\tilde{z}_{3}=z_{\Delta,3}-\tilde{z}_{1}>z_{\Delta,2}+\tau-\tilde{z}_{1}\geq z_{\Delta,2}+\sigma-\tilde{z}_{1},
\end{equation}
\noindent and using the lower bound of $\Delta$ in $\mathit{I_{RPFN}}$ we conclude that
\begin{equation}\nonumber
T(\Delta)>z_{\Delta,2}+\sigma-\tilde{z}_{1}= z_{\Delta,2}+\tau-(t_{max}-\sigma) \geq  z_{\Delta,2}+\tau-\Delta=F(\Delta);
\end{equation}

\noindent\textbf{FPFP:} $\Delta\in\mathit{I_{FPFP}}=[t_{\max},\tilde{z}_{2}]\cap(-\infty,\delta_{2}]$:
\begin{equation}\nonumber
T(\Delta)=\tilde{T}+z_{\Delta,2}-\tilde{z}_{2}=\tau+z_{\Delta,2}>\tau\geq\sigma=F(\Delta);
\end{equation}

\noindent\textbf{FPFN:} $\Delta\in\mathit{I_{FPFN}}=[t_{\max},\tilde{z}_{2}]\cap(\delta_{2},\infty)$: using $z_{\Delta,3}>z_{\Delta,2}$ and using the lower bound of $\Delta$ in $\mathit{I_{FPFN}}$ we have
\begin{equation}\nonumber
T(\Delta)=\tilde{T}+z_{\Delta,3}-\tilde{z}_{3}=z_{\Delta,3}-\tilde{z}_{1}>z_{\Delta,2}-\tilde{z}_{1}= z_{\Delta,2}+\tau-t_{max}\geq z_{\Delta,2}+\tau-\Delta=F(\Delta);
\end{equation}

\noindent\textbf{FNFP1:} $\Delta\in\mathit{I_{FNFP1}}=\Delta\in(\tilde{z}_{2},\tilde{T}-\sigma)\cap(-\infty,\delta_{2})\cap[\delta_{4},\infty)$: here we have $x^{(\Delta)}(z_{\Delta,4}+t)=\tilde{x}(\tilde{z}_{2}+t)$ for all $t\geq0$, which gives $T(\Delta)= z_{\Delta,4}-\tilde{z}_{2}$. Noting that $z_{\Delta,4}>z_{\Delta,3}$ and using the lower bound of $\Delta$ in $\mathit{I_{FNFP1}}$ we obtain
\begin{equation}\nonumber
T(\Delta)=z_{\Delta,4}-\tilde{z}_{2}>z_{\Delta,3}+\tau-\tilde{z}_{2}>z_{\Delta,3}+\tau-\Delta=F(\Delta);
\end{equation}

\noindent\textbf{FNFP2:} Recall that $\mathit{I_{FNFP2}}=\emptyset$ since $\delta_{4}<\tilde{z}_{2}$;

\noindent\textbf{FNFN:} $\Delta\in\mathit{I_{FNFN}}=(\tilde{z}_{2},\tilde{z}_{2}+\tau-\sigma)\cap[\delta_{2},\infty)$: since $z_{\Delta,3}>z_{\Delta,2}+\tau=\tilde{z}_{2}+\tau$ we have
\begin{equation}\nonumber
T(\Delta)=\tilde{T}+z_{\Delta,3}-\tilde{z}_{3}=z_{\Delta,3}-\tilde{z}_{1}>\tilde{z}_{2}+\tau-\tilde{z}_{1} > \tau > \tilde{z}_{2}+\tau-\Delta=F(\Delta);
\end{equation}

\noindent\textbf{FNRN:} $\Delta\in\mathit{I_{FNRN}}=[\tilde{T}-\sigma,\tilde{T})\cap(-\infty,\tilde{T}+\delta_{1})$: using the fact that $z_{\Delta,3}>z_{\Delta,2}+\tau=\tilde{z}_{2}+\tau$ we obtain
\begin{equation}\nonumber
T(\Delta)=\tilde{T}+z_{\Delta,3}-\tilde{z}_{3}=z_{\Delta,3}-\tilde{z}_{1}>\tilde{z}_{2}+\tau-\tilde{z}_{1}>\tau\geq\sigma=F(\Delta);
\end{equation}

\noindent\textbf{FNRP:} $\Delta\in\mathit{I_{FNRP}}=[\tilde{T}-\sigma,\tilde{T})\cap[\tilde{T}+\delta_{1},\infty)$: using that $\tilde{z}_{2}+\tau>\tilde{z}_{1}+2\tau>\tilde{z}_{1}+\tau+\sigma$ we have $\tilde{z}_{1}+\tau-(\tilde{T}-\sigma)<0$, and thus
\begin{equation}\nonumber
T(\Delta)=\tilde{T}+z_{\Delta,4}-\tilde{z}_{4}\geq\tilde{T}+\tilde{z}_{1}+\tau-(\tilde{T}-\sigma)+z_{\Delta,4}-\tilde{z}_{4} > \tilde{z}_{3}+\tau-(\tilde{z}_{4}-z_{\Delta,4})-\Delta=F(\Delta). \tag*{\qedhere}
\end{equation}
\end{proof}
%
\begin{proof}[{\bf Proof of Remark~\ref{rem.min.ress}}]
Each set of parameters $(\tau,\beta_{U},\beta_{L},\sigma,a,\Delta)$ defines a sequence of cases~\eqref{cases} along $\Delta\in[0,\tilde{T})$. Thus we show that for each case~\eqref{cases} we have $F(\Delta)\geq\sigma$:

\noindent\textbf{RNRN}, \textbf{RPFP}, \textbf{FPFP}, \textbf{FNRN:} $F(\Delta)=\sigma$;

\noindent\textbf{RNRP:} from~\citet[Proposition 5.2]{Mackey_2017} we have $x^{(\Delta)}(z_{\Delta,1}+\tau)>\bar{x}$. Using this and $\bar{x}+\beta_{U}=\beta_{U}\mathrm{e}^{\tilde{z}_{2}-t_{max}}$  in $0=x^{(\Delta)}(z_{\Delta,2})=-\beta_{U}+(x^{(\Delta)}(z_{\Delta,1}+\tau)+\beta_{U})\mathrm{e}^{z_{\Delta,1}+\tau-z_{\Delta,2}}$ we obtain
\begin{eqnarray}
\beta_{U}\mathrm{e}^{z_{\Delta,2}} &=& (x^{(\Delta)}(z_{\Delta,1}+\tau)+\beta_{U})\mathrm{e}^{z_{\Delta,1}+\tau} \notag\\
&>& \beta_{U}\mathrm{e}^{\tilde{z}_{2}-t_{max}}\mathrm{e}^{z_{\Delta,1}+\tau}. \notag\\
&=& \beta_{U}\mathrm{e}^{\tilde{z}_{2}+z_{\Delta,1}-\tilde{z}_{1}}. \notag
\end{eqnarray}
\noindent So $(z_{\Delta,2}-\tilde{z}_{2})>(z_{\Delta,1}-\tilde{z}_{1})$ and $z_{\Delta,1}>\Delta$ lead to
\begin{equation}\nonumber
F(\Delta)=\tilde{z}_{1}+\tau+(z_{\Delta,2}-\tilde{z}_{2})-\Delta>\tilde{z}_{1}+\tau+(z_{\Delta,1}-\tilde{z}_{1})-\Delta=z_{\Delta,1}-\Delta+\tau>\tau\geq\sigma;
\end{equation}

\noindent\textbf{RPRP:} taking the upper bound of $\Delta$ in $\mathit{I_{RPRP}}=[\tilde{z}_{1},t_{max}-\sigma]$ and using that $z_{\Delta,2}>\tilde{z}_{2}$, see~\citet[Proposition 5.3]{Mackey_2017}, we find
\begin{equation}\nonumber
F(\Delta)=\tilde{z}_{1}+\tau+(z_{\Delta,2}-\tilde{z}_{2})-\Delta>\tilde{z}_{1}+\tau+z_{\Delta,2}-\tilde{z}_{2}-t_{max}+\sigma=\sigma+z_{\Delta,2}-\tilde{z}_{2}>\tau\geq\sigma;
\end{equation}

\noindent\textbf{RPFN:} taking the upper bound of $\Delta$ in $\mathit{I_{RPFN}}=(t_{\max}-\sigma,t_{\max}]\cap(\delta_{2},\infty)$ and using that $z_{\Delta,2}\geq\tilde{z}_{2}$~\citet[Proof of Proposition 5.5, Eq. (9.5)]{Mackey_2017} we obtain
\begin{equation}\nonumber
F(\Delta)=z_{\Delta,2}+\tau-\Delta\geq z_{\Delta,2}+\tau-t_{\max}>\tilde{z}_{2}+\tau-t_{\max}\geq\tau\geq\sigma;
\end{equation}

\noindent\textbf{FPFN:} from~\citet[Proof of Proposition 5.7]{Mackey_2017} we have $z_{\Delta,2}\geq\tilde{z}_{2}$, then
\begin{equation}\nonumber
F(\Delta)=z_{\Delta,2}+\tau-\Delta\geq z_{\Delta,2}+\tau-\tilde{z}_{2} \geq\tau\geq\sigma;
\end{equation}

\noindent\textbf{FNFP1:} $F(\Delta)=z_{\Delta,3}+\tau-\Delta>z_{\Delta,3}+\tau-(\tilde{z}_{2}+\tau-\sigma)=\sigma+z_{\Delta,3}-\tilde{z}_{2}>\sigma$;

\noindent\textbf{FNFP2:} Recall that $\mathit{I_{FNFP2}}=\emptyset$ since $\delta_{4}<\tilde{z}_{2}$;

\noindent\textbf{FNFN:} $F(\Delta)=\tilde{z}_{2}+\tau-\Delta>\tilde{z}_{2}+\tau-(\tilde{z}_{2}+\tau-\sigma)=\sigma$;

\noindent\textbf{FNRP:} from~\citet[Proof of Proposition 5.10]{Mackey_2017} we have $x^{(\Delta)}(z_{\Delta,3}+\tau)>\bar{x}$. Using this and $\bar{x}+\beta_{U}=\beta_{U}\mathrm{e}^{\tilde{z}_{2}-t_{max}}$  in $x^{(\Delta)}(z_{\Delta,4})=0$ we have
\begin{eqnarray}
\beta_{U}\mathrm{e}^{z_{\Delta,4}} &=& (x^{(\Delta)}(z_{\Delta,3}+\tau)+\beta_{U})\mathrm{e}^{-(z_{\Delta,4}-(z_{\Delta,3}+\tau))} \notag\\
&>& \beta_{U}\mathrm{e}^{\tilde{z}_{2}-t_{max}}\mathrm{e}^{z_{\Delta,3}+\tau}. \notag\\
&=& \beta_{U}\mathrm{e}^{\tilde{z}_{4}+z_{\Delta,3}-\tilde{z}_{3}}. \notag
\end{eqnarray}
\noindent So $(z_{\Delta,4}-\tilde{z}_{4})>(z_{\Delta,3}-\tilde{z}_{3})$ and $z_{\Delta,3}>\Delta$ lead to
\begin{equation}\nonumber
F(\Delta)=\tilde{z}_{3}+\tau+(z_{\Delta,4}-\tilde{z}_{4})-\Delta > \tilde{z}_{3}+\tau+(z_{\Delta,3}-\tilde{z}_{3})-\Delta=z_{\Delta,3}+\tau-\Delta>\tau\geq\sigma. \tag*{\qedhere}
\end{equation}
\end{proof}
%
\begin{proof}[{\bf Proof of Remark~\ref{rem.RPFN}}]
Recall that $\delta_{2}$ is given by~\eqref{delta2} and is defined for $\beta_{U}>a(1-\mathrm{e}^{-\sigma})$. Thus the definition~\eqref{delta2} also holds for $a<\beta_{U}$. For the case \textbf{RPFN} of~\citet[Table 2]{Mackey_2017} we have $\Delta\in\mathit{I_{RPFN}}=(t_{\max}-\sigma,t_{\max}]\cap(\delta_{2},\infty)$ and $\delta_{2}>t_{max}$, then it follows that $\mathit{I_{RPFN}}=\emptyset$.
\end{proof}
%
\begin{proof}[{\bf Proof of Remark~\ref{rem.max.ress}}]
From~\citet[Corollary 4.2]{Mackey_2017} it follows that the cycle length map $T(\Delta)$ is continuous for $a<\beta_{U}$. The proof is divided into two subcases,  $\delta_{2}<t_{max}$ and $\delta_{2}\geq t_{max}$ as follows.

The condition $\delta_{2}<t_{max}$ implies that $a<\beta_{U}$ and there exists a sequence of cases from~\eqref{cases} as is shown in~\citet[Table 1]{Mackey_2017}. Once the cycle length map is continuous we see from~\citet[Third row of Table 1]{Mackey_2017} that $T(\Delta)$ is strictly increasing on $[0,\delta_{2}]$ and strictly decreasing on $[\delta_{2},\tilde{T}]$. Thus the maximum of $T(\Delta)$ occurs for the case \textbf{RPFP} with $\Delta=\delta_{2}$ and we have $\bar{T}=T(\delta_{2})$ with the cycle length map given by~\citet[Proposition 5.4]{Mackey_2017}, i.e.
\begin{equation}\label{Tmax2}
\bar{T} = \tilde{T} +  \ln\left(1+\frac{a(\mathrm{e}^{\sigma}+1)}{\beta_{U}}\mathrm{e}^{\delta_{2}-\tilde{z}_{2}}\right) = \tilde{T} +  \ln\left(\frac{\beta_{U}}{\beta_{U}-a(1-\mathrm{e}^{-\sigma})}\right) > \tilde{T},
\end{equation}
\noindent where $\tilde{T}$ is defined by~\eqref{z1z2T}.

For $\delta_{2}\geq t_{max}$ there exists a sequence of cases from~\eqref{cases} as is shown in~\citet[Table 2]{Mackey_2017} and it follows from Remark~\ref{rem.RPFN} that $a<\beta_{U}$ holds.
For the case \textbf{RPFN} of~\citet[Table 2]{Mackey_2017} we have $\delta_{2}\geq t_{max}$ and from Remark~\ref{rem.RPFN} it follows that $\mathit{I_{RPFN}}=\emptyset$.
Once the cycle length map is continuous and $\mathit{I_{RPFN}}=\emptyset$, we see from~\citet[Third row of Table 2]{Mackey_2017} that again $T(\Delta)$ is strictly increasing on $[0,\delta_{2}]$ and strictly decreasing on $[\delta_{2},\tilde{T}]$.
Thus the maximum of $T(\Delta)$ occurs for the case \textbf{FPFP} with $\Delta=\delta_{2}$ and we have $\bar{T}=T(\delta_{2})$ with the cycle length map given by~\citet[Proposition 5.6]{Mackey_2017}, which is equal to~\eqref{Tmax2}.
\end{proof}
%
\begin{proof}[{\bf Proof of Remark~\ref{rem.deltainf}}]
For $\Delta\leq t \leq \Delta+\sigma$, Eq.~\eqref{FNFPa} together with $x^{(\Delta)}(\Delta) = -\beta_{U} + \beta_{U}\mathrm{e}^{\tilde{z}_{2}-\Delta}$ yields
\begin{equation}\label{FNFPc}
x^{(\Delta)}(t)  = -\beta_{U} + a + (\beta_{U}\mathrm{e}^{\tilde{z}_{2}-\Delta}-a)\mathrm{e}^{-(t-\Delta)}.
\end{equation}
The conditions $a>\beta_{U}$ and $\Delta>\tilde{z}_{2}$ (see Remark~\ref{rem.FNFP}) combined yield $a>\beta_{U}\mathrm{e}^{\tilde{z}_{2}-\Delta}$. Hence $x^{(\Delta)}(t)$ is strictly increasing on $[\Delta,\Delta+\sigma]$.
From $x^{(\Delta)}(\Delta)<0$ and $x^{(\Delta)}(\Delta+\sigma)>0$ we obtain a zero $z_{\Delta,3}$ of $x^{(\Delta)}$ in $(\Delta,\Delta+\sigma)$ given by~\eqref{eq.zDelta3}.

For $\Delta+\sigma\leq t \leq \tilde{T}$, Eq.~\eqref{FNFPb} together with $(x^{(\Delta)}(\Delta+\sigma)+\beta_{U})>0$ shows that $x^{(\Delta)}(t)$ is strictly decreasing on $[\Delta+\sigma,\tilde{T}]$.

For $\tilde{T}<t<z_{\Delta,3}+\tau$ we have $\tilde{z}_{2}<t-\tau<z_{\Delta,3}$. Hence $x^{(\Delta)}(t-\tau)<0$, and
\begin{equation}\nonumber
x^{(\Delta)}(t)  = \beta_{L} + (x^{(\Delta)}(\tilde{T})-\beta_{L})\mathrm{e}^{-(t-\tilde{T})}.
\end{equation}
\noindent Thus $x^{(\Delta)}$ is strictly increasing on  $[\tilde{T},z_{\Delta,3}+\tau]$ since $x^{(\Delta)}(\tilde{T})<0$. For $t=z_{\Delta,3}+\tau$
\begin{equation}\label{eq.xDelta3tau}
x^{(\Delta)}(z_{\Delta,3}+\tau)  = \beta_{L} + (x^{(\Delta)}(\tilde{T})-\beta_{L})\mathrm{e}^{-(z_{\Delta,3}-\tilde{z}_{2})}.
\end{equation}

A rapidly oscillating periodic solution occurs if $x^{(\Delta)}(\Delta)=x^{(\Delta)}(\tilde{T})$, $x^{(\Delta)}(\Delta+\sigma)=x^{(\Delta)}(z_{\Delta,3}+\tau)$ and if the solution $x^{(\Delta)}$ for $\Delta\leq t\leq\Delta+\sigma$ is equal to the solution $x^{(\Delta)}$ for $\tilde{T}\leq t\leq z_{\Delta,3}+\tau$, i.e, $a=\beta_{L}+\beta_{U}$. Then, the necessary conditions for the existence of a rapid oscillation are
\begin{equation}\label{eq.3cond}
\left\{\begin{array}{ll}
a = \beta_{L}+\beta_{U}, \\
x^{(\Delta)}(\Delta)=x^{(\Delta)}(\tilde{T}), \\
x^{(\Delta)}(\Delta+\sigma)=x^{(\Delta)}(z_{\Delta,3}+\tau).
\end{array}\right.
\end{equation}
Combining Eq.~\eqref{eq.3cond} with~\eqref{eq.xDelta3tau} we get
\begin{equation}\label{eq.3condB}
\left\{\begin{array}{ll}
\beta_{L} = a - \beta_{U}, \\
x^{(\Delta)}(\tilde{T}) = x^{(\Delta)}(\Delta), \\
x^{(\Delta)}(\Delta+\sigma) = \beta_{L} + (x^{(\Delta)}(\tilde{T})-\beta_{L})\mathrm{e}^{-(z_{\Delta,3}-\tilde{z}_{2})}.
\end{array}\right.
\end{equation}
Using the first and second relation in the third line of Eq.~\eqref{eq.3condB} gives
\begin{equation}\nonumber
x^{(\Delta)}(\Delta+\sigma) + \beta_{U} - a = (x^{(\Delta)}(\Delta) + \beta_{U} - a)\mathrm{e}^{-(z_{\Delta,3}-\tilde{z}_{2})},
\end{equation}
and combining this with~\eqref{eq.xDeltaSigma} yields
\begin{equation}\nonumber
(x^{(\Delta)}(\Delta) + \beta_{U} - a)\mathrm{e}^{-\sigma} = (x^{(\Delta)}(\Delta) + \beta_{U} - a)\mathrm{e}^{-(z_{\Delta,3}-\tilde{z}_{2})}.
\end{equation}
From $a>\beta_{U}$ and $x^{(\Delta)}(\Delta) = -\beta_{U} + \beta_{U}\mathrm{e}^{\tilde{z}_{2}-\Delta}$ it follows that $(x^{(\Delta)}(\Delta) + \beta_{U} - a)<0$.
Hence the conditions~\eqref{eq.3cond} are reduced to
\begin{equation}\label{delta3z2}
\mathrm{e}^{z_{\Delta,3}}  = \mathrm{e}^{\tilde{z}_{2}+\sigma}.
\end{equation}
The relation~\eqref{delta3z2} combined with~\eqref{eq.zDelta3} yields
\begin{equation}\label{eq.deltaInf3}
a\mathrm{e}^{\Delta} = \beta_{U}\mathrm{e}^{\tilde{z}_{2}} + (a-\beta_{U})\mathrm{e}^{\tilde{z}_{2}+\sigma},
\end{equation}
where $\tilde{z}_{2}$ is given by~\eqref{z1z2T}. Substituting $\Delta=\delta_{\infty}$ in~\eqref{eq.deltaInf3} gives the constant defined $\delta_{\infty}$ in~\eqref{eq.deltaInf2} and
\begin{eqnarray}
\delta_{\infty} &=& \tilde{z}_{2} + \sigma + \ln{\frac{a-\beta_{U}(1-\mathrm{e}^{-\sigma})}{a}} < \tilde{z}_{2} + \sigma \notag\\[1mm]
&=& \tilde{z}_{2} + \ln{\frac{a\mathrm{e}^{\sigma}-\beta_{U}(\mathrm{e}^{\sigma}-1)}{a\mathrm{e}^{\sigma}-a(\mathrm{e}^{\sigma}-1)}} > \tilde{z}_{2}. \notag
\end{eqnarray}
So the conditions~\eqref{eq.3cond} yield $\Delta=\delta_{\infty}$ with $\tilde{z}_{2}<\delta_{\infty}<\tilde{z}_{2}+\sigma$.

The period of the unstable periodic solution is given by $\tilde{T}^{(\infty)}=z_{\Delta,3}+\tau-(\Delta+\sigma)$ (see the example from Figure~\ref{fig_FNFP2_23}). Computing $z_{\Delta,3}$ from~\eqref{delta3z2} and using $\Delta=\delta_{\infty}$ gives $\tilde{T}^{(\infty)}=\tilde{T}-\delta_{\infty}$, and this together with $\tilde{z}_{2}<\delta_{\infty}<\tilde{z}_{2}+\sigma$ yields $\tau-\sigma<\tilde{T}^{(\infty)}<\tau$.

Recalling that $x^{(\Delta)}(t)$ is strictly increasing on $[\Delta,\Delta+\sigma]$, strictly decreasing on $[\Delta+\sigma,\tilde{T}]$ and strictly increasing on  $[\tilde{T},z_{\Delta,3}+\tau]$, we infer
that for $\Delta=\delta_{\infty}$ the minimum and maximum of the rapid limit cycle are respectively given by $\ubar{x}^{(\Delta)}=x^{(\Delta)}(\delta_{\infty})$ and $\bar{x}^{(\Delta)}=x^{(\Delta)}(\delta_{\infty}+\sigma)$. From  $x^{(\Delta)}(\Delta) = -\beta_{U} + \beta_{U}\mathrm{e}^{\tilde{z}_{2}-\Delta}$ with $\Delta=\delta_{\infty}$ and $0<\delta_{\infty}-\tilde{z}_{2}<\sigma<\tau$ it follows that
\begin{eqnarray}
\ubar{x}^{(\Delta)} &=& -\beta_{U}(1-\mathrm{e}^{-(\delta_{\infty}-\tilde{z}_{2})})  \notag\\
&>& -\beta_{U}(1-\mathrm{e}^{-\sigma}) \notag\\
&>& -\beta_{U}(1-\mathrm{e}^{-\tau}) = \ubar{x}. \notag
\end{eqnarray}
From~\eqref{FNFPc} with $t=\delta_{\infty}+\sigma$ and $\Delta=\delta_{\infty}$ together with $a=\beta_{L}+\beta_{U}$ and $0<\delta_{\infty}-\tilde{z}_{2}<\sigma<\tau$ it follows that
\begin{eqnarray}
\bar{x}^{(\Delta)} &=&  a(1-\mathrm{e}^{-\sigma}) + \beta_{U}(\mathrm{e}^{\tilde{z}_{2}-\delta_{\infty}-\sigma} - 1) \notag\\
&=&  (\beta_{L}+\beta_{U})(1-\mathrm{e}^{-\sigma}) + \beta_{U}(\mathrm{e}^{\tilde{z}_{2}-\delta_{\infty}-\sigma} - 1) \notag\\
&=&  \beta_{L}(1-\mathrm{e}^{-\sigma}) - \beta_{U}(1-\mathrm{e}^{-(\delta_{\infty}-\tilde{z}_{2})})\mathrm{e}^{-\sigma} \notag\\
&<&  \beta_{L}(1-\mathrm{e}^{-\sigma}) \notag\\
&<&  \beta_{L}(1-\mathrm{e}^{-\tau}) = \bar{x}. \notag
\end{eqnarray}

\vspace{-7.3mm}
\qedhere
\vspace{2mm}
\end{proof}

\begin{proof}[{\bf Proof of Proposition~\ref{prop.LC1}}]
Using the fact that $x^{(p)}(\Delta_{0})=x^{(p)}(z_{2})=0$, for $t\in[\Delta_{0},\Delta_{0}+\sigma]$ we have
\begin{equation}\nonumber
x^{(p)}(t) = -\beta_{U}+a+(\beta_{U}-a)\mathrm{e}^{-(t-\Delta_{0})}.
\end{equation}
\noindent The condition $a\geq a_{1}$ implies $a>\beta_{U}$, so $x^{(p)}(t)$ is increasing and
\begin{eqnarray}\nonumber
x^{(p)}(\Delta_{0}+\sigma) &=& -\beta_{U}+a+(\beta_{U}-a)\mathrm{e}^{-\sigma} \notag\\
&=& \beta_{U}(\mathrm{e}^{-\sigma}-1)+a(1-\mathrm{e}^{-\sigma}), \notag
\end{eqnarray}
\noindent thus $x^{(p)}(\Delta_{0}+\sigma) = (a-\beta_{U})(1-\mathrm{e}^{-\sigma})>0$.
\par For $t\in[\Delta_{0}+\sigma,\Delta_{1}]$ we have
\begin{equation}\nonumber
x^{(p)}(t) = -\beta_{U}+(x^{(p)}(\Delta_{0}+\sigma) + \beta_{U})\mathrm{e}^{-(t-\Delta_{0}-\sigma)},
\end{equation}
\noindent so $x^{(p)}(t)$ is decreasing, since $(x^{(p)}(\Delta_{0}+\sigma) + \beta_{U}) = \beta_{U}\mathrm{e}^{-\sigma}+a(1-\mathrm{e}^{-\sigma})>0$, and
\begin{eqnarray}\nonumber
x^{(p)}(\Delta_{1}) &=& -\beta_{U}+(x^{(p)}(\Delta_{0}+\sigma) + \beta_{U})\mathrm{e}^{-\alpha} \notag\\
&=& \beta_{U}(\mathrm{e}^{-T_{p}}-1)+a(1-\mathrm{e}^{-\sigma})\mathrm{e}^{-\alpha}. \notag
\end{eqnarray}
\noindent The condition $a\geq a_{1}$ implies $x^{(p)}(\Delta_{1}) \geq0$.
\par For $t\in[\Delta_{1},\Delta_{1}+\sigma]$ it follows that
\begin{equation}\nonumber
x^{(p)}(t) = -\beta_{U}+a+(x^{(p)}(\Delta_{1})+\beta_{U}-a)\mathrm{e}^{-(t-\Delta_{1})},
\end{equation}
\noindent so $x^{(p)}(t)$ is increasing, $(x^{(p)}(\Delta_{1})+\beta_{U}-a) =  (\beta_{U}-a)\mathrm{e}^{-T_{p}}+a(\mathrm{e}^{-\alpha}-1)<0$, and
\begin{eqnarray}
x^{(p)}(\Delta_{1}+\sigma) &=& -\beta_{U}+a+(x^{(p)}(\Delta_{1})+\beta_{U}-a)\mathrm{e}^{-\sigma} \notag\\
&=& \beta_{U}(\mathrm{e}^{-T_{p}-\sigma}-1)+a(1-\mathrm{e}^{-\sigma})(1+\mathrm{e}^{-T_{p}}). \notag
\end{eqnarray}
\noindent Since $x^{(p)}(t)$ is increasing for $t\in[\Delta_{1},\Delta_{1}+\sigma]$ and $x^{(p)}(\Delta_{1}) \geq0$, then $x^{(p)}(\Delta_{1}+\sigma)>x^{(p)}(\Delta_{1})\geq0$.
\par For $t\in[\Delta_{1}+\sigma,\Delta_{2}]$ the solution is given by
\begin{equation}\nonumber
x^{(p)}(t) = -\beta_{U}+(x^{(p)}(\Delta_{1}+\sigma)+\beta_{U})\mathrm{e}^{-(t-\Delta_{1}-\sigma)},
\end{equation}
\noindent so $x^{(p)}(t)$ is decreasing, $(x^{(p)}(\Delta_{1}+\sigma)+\beta_{U}) = \beta_{U}\mathrm{e}^{-T_{p}-\sigma}+a(1-\mathrm{e}^{-\sigma})(1+\mathrm{e}^{-T_{p}})>0$, and
\begin{equation}\nonumber
x^{(p)}(\Delta_{2}) = -\beta_{U}+(x^{(p)}(\Delta_{1}+\sigma)+\beta_{U})\mathrm{e}^{-\alpha} = \beta_{U}(\mathrm{e}^{-2T_{p}}-1)+a(1-\mathrm{e}^{-\sigma})(1+\mathrm{e}^{-T_{p}})\mathrm{e}^{-\alpha}.
\end{equation}
\noindent Thus $x^{(p)}(\Delta_{2}) > \beta_{U}(\mathrm{e}^{-T_{p}}-1)+a(1-\mathrm{e}^{-\sigma})\mathrm{e}^{-\alpha} = x^{(p)}(\Delta_{1})\geq0$.
\par For $t\in[\Delta_{2},\Delta_{2}+\sigma]$ we have
\begin{equation}\nonumber
x^{(p)}(t) = -\beta_{U}+a+(x^{(p)}(\Delta_{2})+\beta_{U}-a)\mathrm{e}^{-(t-\Delta_{2})},
\end{equation}
\noindent so $x^{(p)}(t)$ is increasing, $(x^{(p)}(\Delta_{2})+\beta_{U}-a) = (\beta_{U}-a)\mathrm{e}^{-2T_{p}}+a(\mathrm{e}^{-\alpha}-1)+(\mathrm{e}^{-\alpha}-1)\mathrm{e}^{-2\alpha}<0$, and
\begin{eqnarray}
x^{(p)}(\Delta_{2}+\sigma) &=& -\beta_{U}+a+(x^{(p)}(\Delta_{2})+\beta_{U}-a)\mathrm{e}^{-\sigma} \notag\\
&=& \beta_{U}(\mathrm{e}^{-2T_{p}-\sigma}-1)+a(1-\mathrm{e}^{-\sigma})[1+\mathrm{e}^{-T_{p}}+\mathrm{e}^{-2T_{p}}]. \notag
\end{eqnarray}
\noindent Since $x^{(p)}(t)$ is increasing for $t\in[\Delta_{2},\Delta_{2}+\sigma]$ and $x^{(p)}(\Delta_{2})\geq 0$, then $x^{(p)}(\Delta_{2}+\sigma)>x^{(p)}(\Delta_{2})>0$.
\par For $t\in[\Delta_{2}+\sigma,\Delta_{3}]$ it follows that
\begin{equation}\nonumber
x^{(p)}(t) = -\beta_{U}+(x^{(p)}(\Delta_{2}+\sigma)+\beta_{U})\mathrm{e}^{-(t-\Delta_{2}-\sigma)},
\end{equation}
\noindent so $x^{(p)}(t)$ is decreasing, $(x^{(p)}(\Delta_{2}+\sigma)+\beta_{U}) = \beta_{U}\mathrm{e}^{-2T_{p}-\sigma}+a(1-\mathrm{e}^{-\sigma})[1+\mathrm{e}^{-T_{p}}+\mathrm{e}^{-2T_{p}}]>0$, and
\begin{eqnarray}
x^{(p)}(\Delta_{3}) &=& -\beta_{U}+(x^{(p)}(\Delta_{2}+\sigma)+\beta_{U})\mathrm{e}^{-\alpha} \notag\\
&=& \beta_{U}(\mathrm{e}^{-3T_{p}}-1)+a(1-\mathrm{e}^{-\sigma})[1+\mathrm{e}^{-T_{p}}+\mathrm{e}^{-2T_{p}}]\mathrm{e}^{-\alpha}. \notag
\end{eqnarray}
\noindent Thus $x^{(p)}(\Delta_{3}) > \beta_{U}(\mathrm{e}^{-T_{p}}-1)+a(1-\mathrm{e}^{-\sigma})\mathrm{e}^{-\alpha} = x^{(p)}(\Delta_{2})>0$.
\par For $t\in[\Delta_{3},\Delta_{3}+\sigma]$ the solution is given by
\begin{equation}\nonumber
x^{(p)}(t) = -\beta_{U}+a+(x^{(p)}(\Delta_{3})+\beta_{U}-a)\mathrm{e}^{-(t-\Delta_{3})},
\end{equation}
\noindent so $x^{(p)}(t)$ is increasing, $(x^{(p)}(\Delta_{3})+\beta_{U}-a) = (\beta_{U}-a)\mathrm{e}^{-3T_{p}}+a(\mathrm{e}^{-\alpha}-1)+(\mathrm{e}^{-\alpha}-1)\mathrm{e}^{-3\alpha}<0$, and
\begin{eqnarray}
x^{(p)}(\Delta_{3}+\sigma) &=& -\beta_{U}+a+(x^{(p)}(\Delta_{3})+\beta_{U}-a)\mathrm{e}^{-\sigma} \notag\\
&=& \beta_{U}(\mathrm{e}^{-3T_{p}-\sigma}-1)+a(1-\mathrm{e}^{-\sigma})[1+\mathrm{e}^{-T_{p}}+\mathrm{e}^{-3T_{p}}]. \notag
\end{eqnarray}
\noindent Since $x^{(p)}(t)$ is increasing for $t\in[\Delta_{3},\Delta_{3}+\sigma]$ and $x^{(p)}(\Delta_{3})\geq 0$, then $x^{(p)}(\Delta_{3}+\sigma)>x^{(p)}(\Delta_{3})>0$.
\par Generalizing this procedure, we see that $x^{(p)}(\Delta_{n})$ and $x^{(p)}(\Delta_{n}+\sigma)$ can be written according to Eqs.~\eqref{eq.xpd} and~\eqref{eq.xpdsigma}. The proof is completed inductively for~\eqref{eq.xpd} and~\eqref{eq.xpdsigma} with $n=1,2,3$ and $\Delta_{n}=nT_{p}+\Delta_{0}$. Hence, we see that for $t\in[\Delta_{n},\Delta_{n}+\sigma]$ the solution $x^{(p)}(t)$ is given by Eq.~\eqref{eq.xpt1} and for $t\in[\Delta_{n}+\sigma,\Delta_{n+1}]$ the solution $x^{(p)}(t)$ is given by Eq.~\eqref{eq.xpt2}.
\end{proof}
\begin{proof}[{\bf Proof of Proposition~\ref{prop.LC2}}]
\textbf{(i):} Since $a\geq a_{1}$ we can take the limit $n\longrightarrow\infty$ in Eqs.~\eqref{eq.xpd} and~\eqref{eq.xpdsigma}. We know that $\sum_{k=0}^{\infty}y^{k}=1/(1-y)$ for $|y|<1$, then
\begin{equation}\nonumber
\lim\limits_{n\longrightarrow\infty} x^{(p)}(\Delta_{n}) = -\beta_{U}+a\frac{(1-\mathrm{e}^{-\sigma})}{(1-\mathrm{e}^{-T_{p}})}\mathrm{e}^{-\alpha}= \ubar{x}^{(p)},
\end{equation}
\noindent and
\begin{equation}\nonumber
\lim\limits_{n\longrightarrow\infty} x^{(p)}(\Delta_{n}+\sigma)  = -\beta_{U}+a\frac{(1-\mathrm{e}^{-\sigma})}{(1-\mathrm{e}^{-T_{p}})} = \bar{x}^{(p)}.
\end{equation}

\noindent\textbf{(ii):} Using that $x(0)=\varphi(0)=\ubar{x}^{(p)}$ and $\Delta_{0}=0$, for $t\in[0,\sigma]$ the solution is given by
\begin{equation}\nonumber
x^{(p)}(t) = -\beta_{U}+a+(\ubar{x}^{(p)}+\beta_{U}-a)\mathrm{e}^{-t},
\end{equation}
\noindent once that $x(t-\tau)=\varphi(t)\geq 0$. So $x^{(p)}(t)$ is increasing, since $(\ubar{x}^{(p)}+\beta_{U}-a)=a\mathrm{e}^{-\alpha}(1-\mathrm{e}^{-\sigma})/(1-\mathrm{e}^{-\sigma-\alpha})-a<0$, and
\begin{eqnarray}
x^{(p)}(\sigma) &=& -\beta_{U}+a+(\ubar{x}^{(p)}+\beta_{U}-a)\mathrm{e}^{-\sigma}, \notag\\
&=& -\beta_{U}+a(1-\mathrm{e}^{-\sigma})+a\dfrac{(1-\mathrm{e}^{-\sigma})}{(1-\mathrm{e}^{-T_{p}})}\mathrm{e}^{-T_{p}}= \bar{x}^{(p)}. \notag
\end{eqnarray}
Using that $x^{(p)}(\sigma) = \bar{x}^{(p)}$, for $t\in[\sigma,\Delta_{1}]$ we have
\begin{equation}\nonumber
x^{(p)}(t) = -\beta_{U}+(\bar{x}^{(p)} + \beta_{U})\mathrm{e}^{-(t-\sigma)},
\end{equation}
\noindent so $x^{(p)}(t)$ is decreasing, since $(\bar{x}^{(p)} + \beta_{U}) = a(1-\mathrm{e}^{-\sigma})/(1-\mathrm{e}^{-T_{p}})>0$, and
\begin{eqnarray}
x^{(p)}(\Delta_{1}) &=& -\beta_{U}+(\bar{x}^{(p)} + \beta_{U})\mathrm{e}^{-\alpha}, \notag\\
&=& -\beta_{U}+a\dfrac{(1-\mathrm{e}^{-\sigma})}{(1-\mathrm{e}^{-T_{p}})}\mathrm{e}^{-\alpha} = \ubar{x}^{(p)}. \notag
\end{eqnarray}
Furthermore, once $x^{(p)}(\Delta_{1}) = \ubar{x}^{(p)}$, for $t\in[\Delta_{1},\Delta_{1}+\sigma]$ it follows that
\begin{equation}\nonumber
x^{(p)}(t) = -\beta_{U}+a+(\ubar{x}^{(p)}+\beta_{U}-a)\mathrm{e}^{-(t-\Delta_{1})},
\end{equation}
and thus $x^{(p)}(\Delta_{1}+\sigma) = \bar{x}^{(p)}$. For $t\in[\Delta_{1}+\sigma,\Delta_{2}]$ we have
\begin{equation}\nonumber
x^{(p)}(t) = -\beta_{U}+(\bar{x}^{(p)} + \beta_{U})\mathrm{e}^{-(t-\Delta_{1}-\sigma)}.
\end{equation}
Hence, repeating this process for $t\in[\Delta_{2},\Delta_{2}+\sigma]$, $[\Delta_{2}+\sigma,\Delta_{2}]$, $[\Delta_{3},\Delta_{3}+\sigma]$ and so forth we see that the solution $x^{(p)}(t)$ is given by~\eqref{eq.xp}.
The proof is completed by checking that the \textit{Principle of Mathematical Induction} holds for~\eqref{eq.xp} with $n=1,2,3$ and $\Delta_{n}=nT_{p}+\Delta_{0}$.
\end{proof}

\begin{proof}[{\bf Proof of Proposition~\ref{prop.LC}}]
Consider the maximum~\eqref{eq.xpn2} and define an initial perturbation phase $\Delta_{l}$ such that $x^{(p)}(\Delta_{l})=\ubar{x}^{(p)}$. The solution for the initial pulse is given by~\eqref{FNFPa}~\citep[Section 5.3]{Mackey_2017}
\noindent with $x^{(\Delta)}(\Delta)=-\beta_{U}+\beta_{U}\mathrm{e}^{\tilde{z}_{2}-\Delta}$. Thus $x^{(p)}(\Delta_{l})  = -\beta_{U} + \beta_{U}\mathrm{e}^{\tilde{z}_{2}-\Delta_{l}}=\ubar{x}^{(p)}$. This with $\ubar{x}^{(p)}\geq 0$, since $a\geq a_{1}$, gives
\begin{equation}\label{DeltaL}
\Delta_{l}  = \tilde{z}_{2} + \ln\left(\frac{\beta_{U}}{\beta_{U}+\ubar{x}^{(p)}}\right)\leq \tilde{z}_{2}.
\end{equation}
The proof is divided between the four cases shown in Figure~\ref{fig_PropLC4x4}, where each $\Delta_{0}$ interval is given by: (a) $(\tilde{z}_{1},t_{max})$, (b) $[t_{max},\Delta_{l})$, (c) $[\Delta_{l},\tilde{z}_{2})$, (d) $[\tilde{z}_{2},\tilde{T}+\tilde{z}_{1}]$.

First, we prove case (b) by showing that for $\Delta_{0}\in[t_{max},\Delta_{l})$ the points $x^{(p)}(\Delta_{n})$ converge to~\eqref{eq.xpn1} and the points $x^{(p)}(\Delta_{n}+\sigma)$ converge to~\eqref{eq.xpn2}.
For this case $x^{(p)}(t_{max})\geq x^{(p)}(\Delta_{0})>\ubar{x}^{(p)}\geq0$ and $x^{(p)}(\Delta_{0}-\tau)\geq 0$. Thus for $t=\Delta_{0}$ we have
\begin{equation}\nonumber
x^{(p)}(\Delta_{0})=-\beta_{U}+\beta_{U}\mathrm{e}^{\tilde{z}_{2}-\Delta_{0}}.
\end{equation}
For $t\in[\Delta_{0},\Delta_{0}+\sigma]$ the solution is given by $x^{(p)}(t)=-\beta_{U}+a+(x^{(p)}(\Delta_{0})+\beta_{U}-a)\mathrm{e}^{-(t-\Delta_{0})}$, so
\begin{equation}\nonumber
x^{(p)}(\Delta_{0}+\sigma)=-\beta_{U}+a+(x^{(p)}(\Delta_{0})+\beta_{U}-a)\mathrm{e}^{-\sigma}.
\end{equation}
For $t\in[\Delta_{0}+\sigma,\Delta_{1}]$ we have $x^{(p)}(t)=-\beta_{U}+(x^{(p)}(\Delta_{0}+\sigma)+\beta_{U})\mathrm{e}^{-(t-\Delta_{0}-\sigma)}$, then
\begin{equation}\nonumber
x^{(p)}(\Delta_{1})=-\beta_{U}+(x^{(p)}(\Delta_{0}+\sigma)+\beta_{U})\mathrm{e}^{-\alpha}.
\end{equation}

Repeating the previous steps for $t\in[\Delta_{1},\Delta_{1}+\sigma]$, $t\in[\Delta_{1}+\sigma,\Delta_{2}]$ and so forth, we see that for $t\in[\Delta_{n},\Delta_{n}+\sigma]$ and $n\in\mathbb{N}_{>0}$ the solution is given by $x^{(p)}(t)=-\beta_{U}+a+(x^{(p)}(\Delta_{n})+\beta_{U}-a)\mathrm{e}^{-(t-\Delta_{n})}$, so
\begin{equation}\label{lcn1}
x^{(p)}(\Delta_{n}+\sigma)=-\beta_{U}+a+(x^{(p)}(\Delta_{n})+\beta_{U}-a)\mathrm{e}^{-\sigma},
\end{equation}
and for $t\in[\Delta_{n}+\sigma,\Delta_{n+1}]$ and $n\in\mathbb{N}_{>0}$ we have $x^{(p)}(t)=-\beta_{U}+(x^{(p)}(\Delta_{n}+\sigma)+\beta_{U})\mathrm{e}^{-(t-\Delta_{n}-\sigma)}$, then
\begin{equation}\label{lcn2}
x^{(p)}(\Delta_{n+1})=-\beta_{U}+(x^{(p)}(\Delta_{n}+\sigma)+\beta_{U})\mathrm{e}^{-\alpha}.
\end{equation}
From~\eqref{lcn1},~\eqref{lcn2} and $T_{p}=\alpha+\sigma$ it follows that
\begin{equation}\label{lcn3}
x^{(p)}(\Delta_{n+1})=-\beta_{U}+a\mathrm{e}^{-\alpha}(1-\mathrm{e}^{-\sigma})+(x^{(p)}(\Delta_{n})+\beta_{U})\mathrm{e}^{-T_{p}}.
\end{equation}
Equation~\eqref{lcn3} is recursive and can be rewritten as
\begin{equation}\label{lcn4}
x^{(p)}(\Delta_{n})=-\beta_{U}+a\mathrm{e}^{-\alpha}(1-\mathrm{e}^{-\sigma})\sum_{k=0}^{n-1}\mathrm{e}^{-kT_{p}}+(x^{(p)}(\Delta_{0})+\beta_{U})\mathrm{e}^{-nT_{p}}.
\end{equation}
It is known that $\sum_{k=0}^{\infty}y^{k}=1/(1-y)$ for $|y|<1$, thus taking the limit $n\longrightarrow\infty$ of~\eqref{lcn4} gives
\begin{equation}\label{lcn5}
\lim\limits_{n\longrightarrow\infty} x^{(p)}(\Delta_{n}) = -\beta_{U}+a\frac{(1-\mathrm{e}^{-\sigma})}{(1-\mathrm{e}^{-T_{p}})}\mathrm{e}^{-\alpha}= \ubar{x}^{(p)},
\end{equation}
\noindent and taking the limit $n\longrightarrow\infty$ of~\eqref{lcn1} and using~\eqref{lcn5} yields
\begin{equation}\nonumber
\lim\limits_{n\longrightarrow\infty} x^{(p)}(\Delta_{n}+\sigma)  = -\beta_{U}+a\frac{(1-\mathrm{e}^{-\sigma})}{(1-\mathrm{e}^{-T_{p}})} = \bar{x}^{(p)}.
\end{equation}
Then in the limit $n\longrightarrow\infty$ we see that the solution of~\eqref{eq.diff.pert} converges to the limit cycle given by~\eqref{eq.xp}.

For case (c) the proof is the same as the case (b), but with $\Delta_{0}\in[\Delta_{l},\tilde{z}_{2})$ and hence $x^{(p)}(\Delta_{l})\geq x^{(p)}(\Delta_{0})>x^{(p)}(\tilde{z}_{2})\geq0$.

To prove case (a) we first note that $x^{(p)}(t-\tau)<0$ for $t\in(\tilde{z}_{1},t_{max})$. Given an initial perturbation phase $\Delta_{0}\in(\tilde{z}_{1},t_{max})$, the solution on $(\tilde{z}_{1},t_{max})$ alternates between $x^{(p)}(t)=\beta_{L}+a+(x^{(p)}(\Delta_{n})-\beta_{L}-a)\mathrm{e}^{-(t-\Delta_{n})}$ for $[\Delta_{n},\Delta_{n}+\sigma]$ and $x^{(p)}(t)=\beta_{L}+(x^{(p)}(\Delta_{0}+\sigma)-\beta_{L})\mathrm{e}^{-(t-\Delta_{0}-\sigma)}$ for $[\Delta_{n}+\sigma,\Delta_{n}]$ with $n=0,1,2,\ldots,k$, where the $k$-th index is such that $\Delta_{k-1}< t_{max} \leq \Delta_{k}$. 
Along the interval $(\tilde{z}_{1},t_{max})$ the solution $x^{(p)}(t)$ might have a maximum point if it reaches the value $(\beta_{L}+a)$ in the intervals $[\Delta_{n},\Delta_{n}+\sigma]$ or if it reaches the value $\beta_{L}$ in the intervals $[\Delta_{n}+\sigma,\Delta_{n}]$, otherwise, it will be strictly increasing for $t\in(\tilde{z}_{1},t_{max})$, as is in the examples of Figure~\ref{fig_PropLC4x4}(a). For $t\geq t_{max}$ the proof is the same as the case (b), but using $x^{(p)}(t_{max})$ as initial point $x^{(p)}(\Delta_{0})$.

For case (d) we need to distinguish it between two subcases, the interval $t\in[\Delta_{0},\tilde{T}]$, for which $x^{(p)}(t-\tau)\geq 0$, and the interval $t\in(\tilde{T},z_{p,1}+\tau)$, where $x^{(p)}(t-\tau)<0$ and $z_{p,1}$ is the first zero of $x^{(p)}(t)$ with $t>\Delta_{0}$. For the first interval the solution alternates between $x^{(p)}(t)=-\beta_{U}+a+(x^{(p)}(\Delta_{n})+\beta_{U}-a)\mathrm{e}^{-(t-\Delta_{n})}$ for $[\Delta_{n},\Delta_{n}+\sigma]$ and $x^{(p)}(t)=-\beta_{U}+(x^{(p)}(\Delta_{0}+\sigma)+\beta_{U})\mathrm{e}^{-(t-\Delta_{0}-\sigma)}$ for $[\Delta_{n}+\sigma,\Delta_{n}]$ with $n=0,1,2,\ldots,k$, where the $k$-th interval is such that  $\Delta_{k-1}<(z_{p,1}+\tau)\leq\Delta_{k}$. For the second interval, the solution alternates between $x^{(p)}(t)=\beta_{L}+a+(x^{(p)}(\Delta_{n})-\beta_{L}-a)\mathrm{e}^{-(t-\Delta_{n})}$ for $[\Delta_{n},\Delta_{n}+\sigma]$ and $x^{(p)}(t)=\beta_{L}+(x^{(p)}(\Delta_{0}+\sigma)-\beta_{L})\mathrm{e}^{-(t-\Delta_{0}-\sigma)}$ for $[\Delta_{n}+\sigma,\Delta_{n}]$ with $n=0,1,2,\ldots,k$, where the $k$-th index is such that $\Delta_{k-1}<(\tilde{T}+\tilde{z}_{1})\leq\Delta_{k}$. 
Along the interval $[\Delta_{0},\tilde{T}]$ the solution $x^{(p)}(t)$ might oscillate if it reaches the value $(a-\beta_{U})$ in the intervals $[\Delta_{n},\Delta_{n}+\sigma]$ or if it reaches the value $-\beta_{U}$ in the intervals $[\Delta_{n}+\sigma,\Delta_{n}]$, otherwise, it will be strictly increasing for $t\in[\Delta_{0},\tilde{T}]$, as is in the examples of Figure~\ref{fig_PropLC4x4}(d). For $t>(z_{p,1}+\tau)$ the proof is the same as the case (a), but using $x^{(p)}(z_{p,1}+\tau)$ as initial point $x^{(p)}(\Delta_{0})$.
\end{proof}
\end{appendices}
\bibliography{PeriodicPerturbationsArxivV2}
\end{document}